\documentclass[a4paper,11pt,reqno]{amsart}
\usepackage[utf8]{inputenc}
\usepackage[T1]{fontenc} 
\usepackage[dvipsnames]{xcolor}
\usepackage{textcomp}
\usepackage[utf8]{inputenc}
\usepackage[T1]{fontenc}
\usepackage{amsmath,amsthm,amssymb}
\usepackage{newtxtext,newtxmath}
\usepackage{mathtools}
\usepackage{bm}
\usepackage{wrapfig}
\usepackage{subcaption}
\usepackage[square,numbers,sort&compress]{natbib}
\usepackage{fullpage}
\usepackage{algorithm}
\usepackage{algpseudocode}
\usepackage[margin=1in]{geometry}
\usepackage[backref,colorlinks=true,linktoc=page]{hyperref}
\usepackage{etoolbox}
\usepackage{enumitem}
\usepackage{tikz}
\usetikzlibrary{positioning,arrows.meta}
\usepackage{tabularx}
\usepackage{pifont}
\usepackage{booktabs}
\usepackage{enumitem}
\usepackage{graphicx}
\usepackage{color}
\usepackage{microtype}
\usepackage{url}
\usepackage{booktabs}
\usepackage{microtype}
\usepackage{float}
\usepackage{pgfplots}
\usepackage{caption}

\definecolor{dark2green}{rgb}{0.1, 0.65, 0.3}
\definecolor{dark2orange}{rgb}{0.9, 0.4, 0.}
\definecolor{dark2purple}{rgb}{0.4, 0.4, 0.8}

\definecolor{sns_green}{HTML}{2CA02C}
\definecolor{sns_orange}{HTML}{FF7F0E}
\definecolor{sns_blue}{HTML}{1F77B4}

\colorlet{first_color}{sns_green!40}
\colorlet{second_color}{sns_blue!50}
\colorlet{third_color}{sns_orange!40}

\makeatletter
\def\@footnotecolor{red}
\define@key{Hyp}{footnotecolor}{%
 \HyColor@HyperrefColor{#1}\@footnotecolor%
}
\patchcmd{\@footnotemark}{\hyper@linkstart{link}}{\hyper@linkstart{footnote}}{}{}
\makeatother

\definecolor{intcolor}{HTML}{CA0020}
\definecolor{extcolor}{HTML}{0571B0}
\hypersetup{
linkcolor=intcolor
,citecolor=intcolor
,filecolor=extcolor
,urlcolor=extcolor
,menucolor=intcolor
,runcolor=extcolor
,linkbordercolor=intcolor
,citebordercolor=intcolor
,filebordercolor=extcolor
,urlbordercolor=extcolor
,menubordercolor=intcolor
,runbordercolor=extcolor
,footnotecolor=intcolor
}

\newtheorem{theorem}{Theorem}[section]

\theoremstyle{definition}

\newtheorem{remark}[theorem]{Remark}

\numberwithin{equation}{section}

\newcommand{\R}{\mathbb{R}}
\newcommand{\E}{\mathbb{E}}
\newcommand{\X}{\mathcal{X}}

\renewcommand{\H}{\mathcal{H}}
\newcommand{\calL}{\mathcal{L}}
\newcommand{\vol}{\operatorname{vol}}

\newcommand{\FMBoost}{\textsc{FlowBoost}}

\title{Flow-based Extremal Mathematical Structure Discovery}
\author[G. B\'erczi]{Gergely B\'erczi}
\address{Aarhus University, Denmark}
\email{gergely.berczi@math.au.dk}
\author[B. Hashemi]{Baran Hashemi}
\address{Max Planck Institute for Mathematics in the Sciences, Leipzig, Germany}
\email{baran.hashemi@mis.mpg.de}
\author[J. Kl\"uver]{Jonas Kl\"uver}
\address{Aarhus University, Denmark}
\email{jonas.kluever@gmail.com}

\date{\today}
\pgfplotsset{compat=1.18}

\begin{document}

\begin{abstract}
    The discovery of extremal structures in mathematics requires navigating vast and nonconvex landscapes where analytical methods offer little guidance and brute-force search becomes intractable. We introduce \FMBoost{}, a closed-loop generative framework that learns to discover rare and extremal combinatorial structures by combining three components: (i) a geometry-aware conditional flow-matching model that learns to sample high-quality configurations, (ii) reward-guided policy optimization with action exploration that directly optimizes the generation process toward the objective while maintaining diversity, and (iii) stochastic local search (Stochastic Relaxation with Perturbations) for both training-data generation and final refinement. 
    Unlike prior open-loop approaches, such as PatternBoost \cite{charton2024patternboost} that either retrains on filtered discrete samples, or AlphaEvolve \cite{novikov2025alphaevolve} where they rely on frozen Large Language Models (LLMs) as evolutionary mutation operators, \FMBoost{} enforces geometric feasibility during sampling, and propagates reward signal directly into the generative model, closing the optimization loop and requiring much smaller training sets and shorter training times, and reducing the required outer-loop iterations by orders of magnitude, while eliminating dependence on LLMs. 
    We demonstrate the framework on four geometric optimization problems: sphere packing in hypercubes, circle packing maximizing sum of radii, the Heilbronn triangle problem, and star discrepancy minimization. In several cases, \FMBoost{} discovers configurations that match or exceed the best known results. For circle packings, we improve the best known lower bounds, surpassing the LLM-based system AlphaEvolve while using substantially fewer computational resources. For the Heilbronn problem, \FMBoost{} improves the minimum triangle area approaching the best known numerical values. For sphere packing in dimension $12$, our method finds configurations denser than those produced by classical heuristics. 
    To our knowledge, \FMBoost{} is the first systematic application of flow-based generative models with Reinforcement Learning to extremal structure discovery in pure mathematics. We term this paradigm \emph{de novo mathematical structure design}. Our results demonstrate that lightweight, domain-specific generative models, when equipped with geometric inductive biases and reward-guided learning, can match or exceed the performance of LLM-based systems at a fraction of the computational cost. The code is available at \url{https://github.com/berczig/FlowBoost}.
\end{abstract}

\maketitle

\tableofcontents

\section{Introduction}

AI-driven mathematics \cite{he2024ai} has advanced in two main directions: (i) formal reasoning and verification \cite{de2015lean,bansal2019holist,agrawal2211towards,azerbayev2023proofnet,yang2025position,schmitt2025improofbench}, where machine learning is integrated with proof assistants to propose proof steps and complete formal proofs, and (ii) discovery via interpretability analysis \cite{davies2021advancing,dong2024machine,hashemican,he2025murmurations,schmitt2025extremal} in a human-AI collaboration, and search \cite{wagner2021constructions,berczi2023ml,gukov2023searching,charton2024patternboost,georgiev2025mathematical,chervov2025cayleypy,alfarano2024global,shehper2025what,BercziWagner2024PercolatingSets}, where neural networks either propose counter-examples, relations, and candidates (programs or constructions) that are evaluated by human mathematicians. While these approaches have achieved notable successes in symbolic and discrete domains, a distinct class of problems lies at their boundary, namely continuous optimization over geometric configurations, where objectives are not normalized (there is no access to likelihood), admit no closed-form gradient, constraints encode hard geometric feasibility, and the landscape admits exponentially many local optima. Such problems are extremely common in combinatorial geometry and algebra, but they resist brute-force search because the search domain is continuous. This is precisely where deep generative models, trained to capture complex distributions over high-dimensional spaces, may offer a distinctive advantage.

Many structure discovery problems in extremal geometry can be re-formulated as finite-dimensional continuous Simulation Based Optimization (SBO) \cite{amaran2016simulation} problem. These problems share a common setup: a configuration space $\mathcal{X} \subseteq \mathbb{R}^{d \times N}$ of geometric objects, an objective functional $J: \mathcal{X} \to \mathbb{R}$ to be optimized, and hard constraints (geometric or algebraic) defining feasibility. Typical examples include packing problems (configurations of spheres or circles maximizing density or coverage), discrepancy problems (point sets minimizing irregularity of distribution), and extremal point configurations (maximizing minimum distances or triangle areas). More broadly, our pipeline can be efficiently used to discover rare candidates in constrained algebraic–combinatorial search spaces, such as graph, or polytope encoded generators.
Analytically, these problems are often intractable beyond low dimensions or small instance sizes. Computationally, they exhibit rugged energy landscapes with exponentially many local optima separated by high barriers. Progress has historically relied on a combination of mathematical insight and carefully designed heuristic search.

Modern deep generative models, such as Transformer \cite{vaswani2017attention}, Diffusion \cite{ho2020denoising}, Normalizing Flows \cite{rezende2015variational}, Generative Adversarial networks \cite{goodfellow2020generative}, Variational Auto Encoders (VAEs) \cite{kingma2013auto} and Flow Matching \cite{lipman2022flow} models, have shown remarkable capacity at learning complex distributions in high-dimensional spaces. These methods have transformed domains from image synthesis \cite{esser2024scaling} to protein structure prediction\cite{watson2023novo}, and physics simulation \cite{hashemi2024deep}, demonstrating an ability to capture intricate structural regularities from data. Hence, a natural question arises: 

\noindent \emph{Can deep generative models accelerate mathematical discovery and exploration by learning to generate high-quality Out-of-Distribution (OOD) configurations directly?} We call this paradigm \textbf{\emph{De Novo Mathematical Structure Design}}. 

Recent work has begun to explore this direction by combining generative discovery with explicit search-and-evaluate loops. These evolutionary search models are used to iteratively boost the quality of mathematical constructions, progressively reshaping the search process toward better solutions.

The \textsc{PatternBoost} framework~\cite{charton2024patternboost} alternates between a \emph{local} classical optimization phase and a global generative modeling phase. A problem-specific local search heuristic generates high-scoring objects; a sequence-to-sequence transformer model \cite{vaswani2017attention} is then trained on the best objects and sampled to provide \emph{seeds} that restart the local search, iterating this local-global search to discover improved constructions.
\textsc{FunSearch}~\cite{romera2024mathematical} replaces local search of PatternBoost by explicit program evolution. It searches directly in the space of functions/programs where a pretrained LLM proposes variants of a python script  (by editing a specific function), while an automatic evaluator executes the code on a set of instances to assign a deterministic score. Candidates scripts are retained, clustered (e.g.\ by behavioral signatures on the evaluation set), and iteratively recombined/mutated through further LLM proposals, yielding an evolutionary procedure in which the LLM acts as a powerful mutation operator and the evaluator enforces correctness and selection pressure. More recently, \textsc{AlphaEvolve}~\cite{novikov2025alphaevolve,georgiev2025mathematical} scales this paradigm into a more general-purpose evolutionary coding agent. 
Rather than evolving a single short function, AlphaEvolve provides a pipeline in which LLMs propose code edits, the fixed problem-specific scoring system in the evaluation harness measures fitness, and population-based selection retains and refines the most promising candidates across many generations, enabling algorithm discovery and optimization across diverse domains. 

In short, FunSearch and AlphaEvolve are evolutionary search models with states represented by python scripts, values are given by a fixed scoring system and actions are LLM-suggested modifications of scripts. 
However, these models can face limitations when applied for continuous geometric optimization problems. Here we list four type of problem which arises.

\begin{enumerate}
	\item \emph{Discrete representations for continuous problems.} PatternBoost operates on tokenized sequences, requiring geometric configurations to be discretized. As a result, it loses the smooth and high-precision structure of the optimization landscape and complicates constraint satisfaction.
	
	\item \emph{Open-loop iteration without convergence guarantees.} PatternBoost and AlphaEvolve iterate a generate-evaluate-select system, but the generative model receives no direct feedback from the objective. They learn to match the distribution of current best solutions rather than to optimize or rely on evolutionary selection with a frozen LLM. This can require significantly amount of  iterations without any guarantee of improvement.
	
	\item \emph{Dependence on Large Language Models.} AlphaEvolve's power derives from frontier LLMs (Gemini Pro, $\sim$10B+ parameters), requiring substantial computational infrastructure and API access beyond the reach of most researchers.
	
	\item \emph{Post-hoc constraint handling.} Geometric constraints (non-overlap, boundary containment, symmetries, algebraic relations, etc.) are typically enforced through rejection sampling, repair heuristics after generation, or prompt engineering, rather than being incorporated into the generative process itself.
\end{enumerate}

In this work we introduce \FMBoost, a framework that addresses each of these limitations through three innovations:

\begin{itemize}
	\item \emph{Continuous generation via CFM.} We operate directly in configuration space $\mathbb{R}^{d \times N}$, learning a time-dependent vector field that transports a prior to the distribution of high-quality configurations. This preserves the geometric structure of the problem and enables gradient-based reasoning about constraints.
	
	\item \emph{Geometry-Aware sampling (GAS).} Rather than generating configurations and filtering invalid ones, we interleave flow integration with projection onto the constraint manifold. This generates feasible samples throughout the generative process, dramatically improving sample efficiency.
	
	\item \emph{Closed-loop Reward-Guided optimization.} We close the optimization loop by fine-tuning the flow model with direct feedback from the objective. Using online reward-guided flow matching with problem-specific action exploration, we upweight high-reward rare samples in the training objective while a trust-region self-distillation term prevents distribution collapse. This transforms open-loop iteration into closed-loop optimization, reducing required iterations from $\mathcal{O}(100)$ to $\mathcal{O}(10)$.
\end{itemize}

\begin{table}[!ht]
\centering
\small
\renewcommand{\arraystretch}{1.15}
\setlength{\tabcolsep}{6pt}

\begin{tabularx}{\linewidth}{@{}p{3.1cm}XXX@{}}
\toprule
 & \textsc{PatternBoost} & \textsc{AlphaEvolve} & \textsc{FlowBoost} \\
\midrule
\textbf{Generative Model} & 
Autoregressive Transformer (Makemore) & 
Frontier LLMs & 
Conditional Flow Matching (CFM) \\
\textbf{Model Scale} & 
$\sim$10M parameters & 
$\sim$10B+ parameters (frozen) & 
$\sim$2M parameters \\
\textbf{Search Space} & 
Discrete token sequences & 
Program space (code) & 
Continuous configuration space \\
\textbf{Inductive Bias} & 
Sequential structure &  
LLM priors & 
Permutation Equivariance, Geometric Consraints \\
\textbf{Exploration mechanism} & 
Random Sampling & 
LLM-guided mutation, MAP-Elites & 
Online Reward-guidance + Action Exploration \\
\textbf{LLM Required?} & 
\textcolor{green}{\ding{55}} No & 
\textcolor{red}{\ding{51}} Yes (essential) & 
\textcolor{green}{\ding{55}} No \\
\textbf{Constraint Handling} & 
Greedy repair heuristics & 
Post-hoc verification & 
Geometry-Aware Sampling \\
\textbf{Iterations to converge} & 
$\mathcal{O}(10\text{--}100)$ & 
$\mathcal{O}(10\text{--}10^6)$ & 
$\mathcal{O}(1\text{--}10)$ \\
\textbf{Loop structure} & 
Open (No Direct Feedback) & 
Open (Evolutionary Selection) & 
Closed (Reward-Weighted Updates) \\
\textbf{Compute} & 
Single GPU ($\sim$10–100\,h) & 
Cluster + API ($\gg$1000\,h) & 
Single GPU ($\sim$1–10\,h) \\
\textbf{Requirements} & 
Objective + Tokenization & 
Evaluator + Prompt engineering & 
Objective \\
\bottomrule
\end{tabularx}
\caption{\textbf{Comparison of AI-driven mathematical structure discovery systems.} We compare architectural choices, computational requirements, and demonstrated capabilities. \FMBoost{} achieves competitive or superior results on continuous geometric problems with 2-5 orders of magnitude less compute, and no reliance on LLMs.}
\label{tab:method_comparison}
\end{table}

The fundamental distinction between our approach and LLM-based systems lies in how structure is represented and exploited. AlphaEvolve-family of models \cite{novikov2025alphaevolve,openevolve,lange2025shinkaevolve} operate in \emph{program space}: the LLM proposes syntactic modifications to code, and an evaluator checks whether the resulting program produces valid, high-scoring solutions. This approach inherits the power of language models where they can generate human-readable algorithms. 
However, for many mathematical problems where the moduli of solutions admits natural algebraic and geometric structures and follow certain set of sharp optimization objectives, this indirection introduces a huge complexity and over-reliance on SOTA LLMs. 
While AlphaEvolve requires access to frontier LLMs and substantial computational infrastructure, \FMBoost{} can be trained and deployed by individual researchers on commodity hardware. When this framework is instantiated with domain-appropriate generative models with few simulated data, comparable or superior results can be achieved at a fraction of the computational cost. Our pipeline is largely reusable across problems, with only modest adjustments to the penalty structure, conditioning variables, and reward definition. In short, when the problem domain admits natural geometric structure within continuous or discrete optimization problems, encoding this structure directly into the generative model as an inductive bias, yields superior results with far less computation. This parallels developments in computational biology, where domain-specific architectures such as AlphaFold \cite{jumper2021highly} and ESM \cite{lin2023evolutionary}, and equivariant diffusion models for molecular generation \cite{watson2023novo} have largely supplanted general-purpose language models for protein structure prediction and synthesis. 

The contributions of this paper are threefold:
\begin{enumerate}
	\item \textbf{A closed-loop deep generative framework for mathematical structure discovery.} We introduce \FMBoost, combining conditional flow matching, Transformers, geometry-aware sampling, and reward-guided policy optimization into a unified pipeline for rare structure discovery. The framework is problem-agnostic: the same architecture applies across problem domains with only changes to the constraint structure and reward function.
	
	\item \textbf{First systematic application of flow-based models to extremal mathematics.} We demonstrate that modern generative models, when properly adapted, can tackle problems in combinatorial geometry, not as benchmarks, but as genuine research tools capable of discovering new constructions.
	
	\item \textbf{New records and near-optimal constructions across multiple problems:}
	\begin{itemize}
		\item \emph{Circle packing (sum of radii):} New best constructions for $n=26$ and $n=32$ circles, surpassing AlphaEvolve.
		\item \emph{Sphere packing in hypercubes:} High-density configurations in dimension $3$ ($N=50$-$200$) and dimension $12$ ($N=31$), matching or improving known constructions.
		\item \emph{Heilbronn problem:} For $n=13$ points, improvement from $A_{\min} = 0.025236$ to $0.025727$ in four iterations, approaching the best known value.
	\end{itemize}
    
    \item \textbf{Bridging two separate research areas.} We identify a unifying perspective showing that, despite superficial differences in problem formulation, rare structure discovery in mathematics and OOD generative design in domains such as molecular and materials discovery share the same underlying SBO blueprint, constrained and likelihood-free generative OOD sampling. This bridge enables direct transfer of methods, theoretical tools, and practical heuristics between the two communities.
\end{enumerate}

The paper is organized as follows. Section~\ref{sec:method} formalizes the simulation-based optimization framework and our closed-loop extension, details the \FMBoost{} architecture: conditional flow matching, geometry-constrained sampling, and reward-guided policy optimization. Sections~\ref{sec:results} present applications to sphere packing, circle packing, the Heilbronn problem, and star discrepancy. Section~\ref{sec:discussion} discusses implications and future directions.

\section{Methods}
\label{sec:method} 
\subsection{Simulation-Based Optimization (SBO): \textit{Learning a Stochastic Policy}}

The discovery of rare and extremal mathematical structures, from dense sphere packings to point configurations maximizing geometric functionals, requires navigating high-dimensional, non-convex landscapes with rugged local optima. We introduce \FMBoost{}, a unified framework that combines conditional flow matching with geometry-aware sampling and reward-guided fine-tuning. This methodology goes under the hood a broader paradigm, namely \emph{Simulation-Based Optimization} (SBO) \cite{gosavi2015simulation,amaran2016simulation,bengio2021flow,jain2022biological}, the iterative alternation between local refinement and global generative exploration.

SBO is closely related to the paradigm of \emph{Simulation-Based Inference} (SBI) \cite{cranmer2020frontier,deistler2025simulation}, but the two differ in what is being learned and what the simulator is used for. In SBI, one assumes access to a simulator (or implicit model) $p(x\mid \theta)$ and aims to infer latent parameters $\theta$ from observations $x$ by learning an approximation to the posterior $p(\theta\mid x)$ or a surrogate for the likelihood $p(x\mid \theta)$. Hence, the simulator is treated as a generative process that produces data, and the learned model is used for statistical inference under a fixed data generating mechanism. SBO instead treats simulation and evaluation as an \emph{optimization oracle}. Here the latent variable is the candidate solution itself, and the objective induces a target distribution over solutions. Hence one can define an \emph{energy-based} (Boltzmann/Gibbs) target
\[
\pi_\beta(x)
=
\frac{1}{Z_\beta}\, p_0(x)\,\exp \big(\beta\,J(x)\big) 
\qquad
Z_\beta = \int_{\mathcal{X}} p_0(x) \exp \big(\beta\,J(x)\big) dx,
\]

where $p_0(x)$ is a \emph{prior} over configurations (e.g.\ a simple base distribution), and $\exp(\beta J(x))$ is an (unnormalized) objective-induced measure that upweights high-scoring configurations. Equivalently, defining an energy $E(x)=-J(x)$ yields $\pi_\beta(x)\propto p_0(x)\exp(-\beta E(x))$.
In this view, SBO is not posterior inference over parameters; it is \emph{posterior sampling over solutions} under an objective-induced measure. The central task is therefore to learn a generative model $p_\theta(x)$ (amortized sampler) so that it concentrates mass near $\pi_\beta$ (with the amortized sampler playing the role that variational posteriors play in SBI), with larger $\beta$ corresponding to stronger optimization pressure. This problem, converting an energy function into a generative distribution, admits several classical solutions. Markov Chain Monte Carlo (MCMC) methods can sample from $\pi_\beta$ given only the unnormalized density, but they are computationally expensive (requiring many iterations per sample) and fundamentally local in their exploration: mixing between distant modes is slow, and discovering new modes requires the chain to traverse low-probability regions~\cite{bengio2021flow}. When the modes occupy a tiny volume of the configuration space (extremal regions), the probability of initializing a chain near an unknown mode becomes negligible, rendering MCMC unsatisfactory for rare structure discovery.

Deep generative models offer an alternative to \emph{amortize} the cost of sampling by training a neural network to map noise to high-quality configurations directly. Sampling then becomes a single forward pass, the model can generalize to propose candidates in regions it has never explicitly visited, provided there is learnable structure relating known modes to undiscovered ones. This amortization fundamentally changes the exploration-exploitation trade-off: whereas MCMC or heuristic algorithms must mix between modes (a slow, iterative process), a generative model can jump to new modes if its inductive biases and training procedure enable appropriate generalization.

Both SBI and SBO rely on repeated simulation and amortized learning, and both may use the same modern generative machinery (flows, diffusion/flow matching) to represent complex distributions. The key distinction is the feedback signal. SBI updates the model to match the simulator-implied data distribution (likelihood or posterior consistency), whereas SBO updates the model to concentrate probability mass on rare, high-scoring regions defined by the objective function and the constraints. 
This distinction also explains why neither local search nor global generation alone suffices for hard extremal problems. local search provides high-precision refinement but is basin-limited, while pure generation provides exploration but struggles with tight feasibility and extrapolative OOD discovery. By embedding both within a boosting loop where a generative model proposes diverse candidates and local refinement sharpens them, SBO can escape local minima while progressively shifting the proposal toward the objective-induced posterior over solutions. From a statistical perspective, \FMBoost\ is related and complimentary to Boltzmann-generator approaches in molecular and materials modeling \cite{noe2019boltzmann}, which learn flow-based transports from a simple prior to a Gibbs target $\pi(x)\propto \exp(-\beta E(x))$; in our setting we induce an analogous target by taking $E(x)=-J(x)$ (along with geometric feasibility constraints), and extend the paradigm with geometry-aware sampling and closed-loop reward-guided updates for extrapolative discovery.

\subsubsection*{Comparison to prior work.}
Our pipeline differs from prior classical SBO method in mathematical structure discovery, PatternBoost, in two critical respects. First, we replace discrete autoregressive sequence models with continuous conditional flow matching, a natural choice for geometric problems where the configuration space admits both continuous coordinates and geometric constraint functions. Second, we close the optimization loop through online reward-guided policy updates, transforming an open iterative scheme into a convergent closed-loop control system. For example, PatternBoost is \emph{open-loop} where the generator is retrained on selected solutions, but does not receive direct objective feedback during the update. Therefore, there is no mechanism guaranteeing convergence to optimal or near-optimal solutions. In PatternBoost, one only hopes that an extremal and rare construction would emerge after $\mathcal{O}(100)$ iterations of the local-global cycle. In practice, however, the distribution of generated solutions can plateau, oscillate, or drift without systematic improvement. The generative model learns to mimic the current best solutions but receives no direct signal pushing it toward better solutions, it matches a distribution rather than optimizing an objective.
In contrast, \FMBoost\ introduces a closed-loop update in which generated candidates are scored and immediately used to fine-tune the generator via self-distillation reward-guidance with a trust region.
As a result, we formalize \FMBoost\ as a closed-loop instance of SBO for rare mathematical structure discovery in generic configuration spaces.  

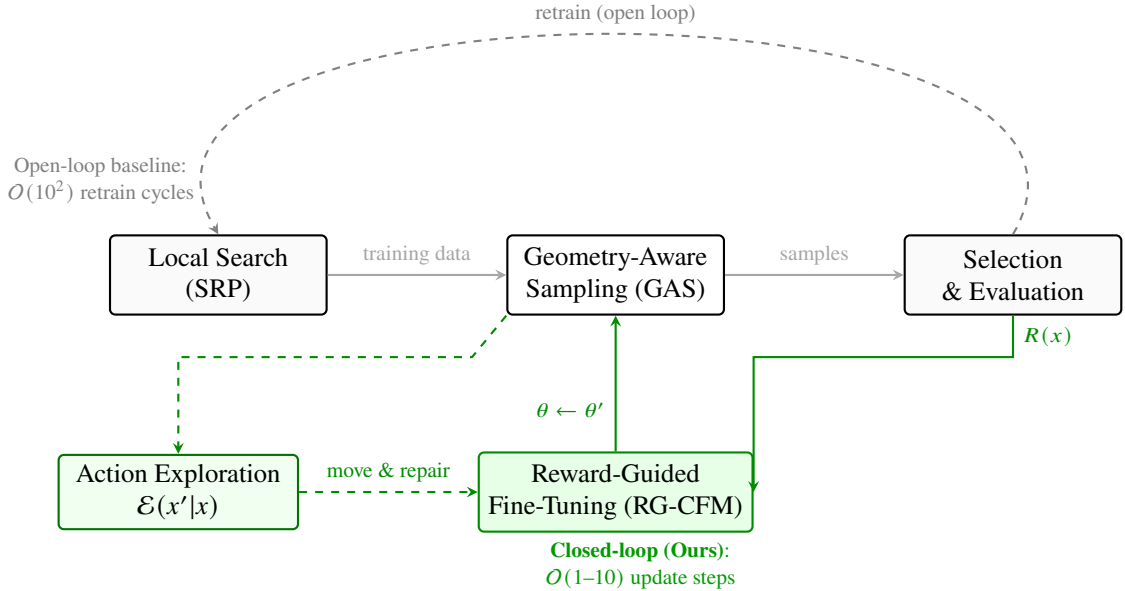
\begin{figure}[ht]
\centering
\begin{tikzpicture}[
    node distance=2.35cm,
    box/.style={rectangle, draw, rounded corners=2pt, minimum width=2.85cm, minimum height=1.05cm, align=center, font=\small, thick},
    boxsoft/.style={rectangle, draw, rounded corners=2pt, minimum width=2.85cm, minimum height=1.05cm, align=center, font=\small, thick, fill=gray!4},
    reward/.style={rectangle, draw, rounded corners=2pt, minimum width=3.15cm, minimum height=1.05cm, align=center,
                   fill=green!10, draw=green!55!black, thick, font=\small},
    expl/.style={rectangle, draw, rounded corners=2pt, minimum width=3.15cm, minimum height=0.95cm, align=center,
                 fill=green!6, draw=green!45!black, thick, font=\small},
    arrow/.style={->, >=stealth, thick},
    openarrow/.style={->, >=stealth, thick, dashed, gray},
    closedarrow/.style={->, >=stealth, thick, draw=green!55!black},
    closedarrow2/.style={->, >=stealth, thick, draw=green!55!black, dashed}
]

\node[boxsoft] (local) {Local Search\\(SRP)};
\node[box] (gen) [right=of local] {Geometry-Aware\\Sampling (GAS)};
\node[boxsoft] (select) [right=of gen] {Selection\\\& Evaluation};

\node[reward] (rw) [below=1.8cm of gen] {Reward-Guided\\Fine-Tuning (RG-CFM)};
\node[expl] (explore) [left=of rw] {Action Exploration\\$\mathcal{E}(x'|x)$};

\draw[arrow, draw=gray!70] (local) -- node[above, font=\scriptsize, text=gray!70] {training data} (gen);
\draw[arrow, draw=gray!70] (gen) -- node[above, font=\scriptsize, text=gray!70] {samples} (select);
\draw[openarrow] (select.north) to[out=60, in=120]
    node[above, font=\scriptsize, gray] {retrain (open loop)} (local.north);

\draw[closedarrow] (select.south) -- ++(0,-0.55)
    node[midway, right, font=\scriptsize, text=green!55!black] {$R(x)$}
    -| (rw.east);
\draw[closedarrow] (rw.north) -- ++(0,0.55) -| 
    node[pos=0.25, left, font=\scriptsize, text=green!55!black] {$\theta \leftarrow \theta'$}
    (gen.south);

\draw[closedarrow2] (gen.south west) -- ++(-0.45,-0.55) -| (explore.north);
\draw[closedarrow2] (explore.east) -- node[midway, above, font=\scriptsize, text=green!55!black] {move \& repair} (rw.west);

\node[font=\scriptsize, gray, align=center] at (-1.55, 1.25)
    {Open-loop baseline:\\$\mathcal{O}(10^2)$ retrain cycles};
\node[font=\scriptsize, green!60!black, align=center] at (5.55, -3.85)
    {\textbf{Closed-loop (Ours)}:\\$\mathcal{O}(1\text{--}10)$ update steps};

\end{tikzpicture}


\[
\theta \leftarrow \arg\min_\theta 
\mathbb{E}_{x \sim p_\theta}\!\Big[w\!\big(R(x)\big)\,\|v_\theta(x_t,t)-v_t^*\|^2\Big]
+
\alpha \,\|v_\theta - v_{\mathrm{ref}}\|^2
\]

\caption{\textbf{Simulation-based optimization with closed-loop updates and action exploration.}
Gray dashed arrows indicate the open-loop baseline: generate candidates, select, and retrain on the selected set over many cycles.
Green arrows indicate our closed-loop pipeline: samples produced by geometry-aware sampling are perturbed by an action exploration operator $\mathcal{E}(x'|x)$, scored by rewards $R(x)$, and used for online reward-guided fine-tuning of the flow model.
The reward-weighted objective upweights high-reward samples via $w(R(x))$, while the teacher, student consistency term $\|v_\theta - v_{\mathrm{ref}}\|^2$ limits distribution shift and mitigates generative collapse.}
\label{fig:sbo_closure}
\end{figure}

%
Empirically, this converts slow, stochastic improvement across many outer cycles into a small number of effective update rounds, because the generator update is explicitly biased by the objective functional rather than implicitly by imitation of the current elite set.

\subsection{Conditional Flow Matching for Configuration Generation}

Let $\X \subseteq \R^{d \times N}$ denote a configuration space of $N$ objects in $d$-dimensional Euclidean space, equipped with an objective functional $J \colon \X \to \R$ to be maximized. Rather than directly optimizing $J$, we learn a generative model whose samples concentrate on high-quality (exceptional) regions of $\X$.

We adopt the conditional flow matching (CFM) \cite{lipman2022flow} framework, which learns a time-dependent vector field $v_\theta \colon \X \times [0,1] \to T\X$ that transports a simple prior distribution $p_0$ on $\X$ (e.g., uniform points in a bounding box) to a data distribution $\mu_{\mathrm{data}}$ of high-quality configurations. The flow is defined by the ordinary differential equation
\begin{equation}
\frac{\mathrm{d} x_t}{\mathrm{d} t} = v_\theta(x_t, t), \quad x_0 \sim p_0,
\label{eq:flow_ode}
\end{equation}
with the objective that the push-forward of $p_0$ under this flow approximates $\mu_{\mathrm{data}}$ at $t=1$. Training proceeds via regression onto the conditional velocity field of an optimal transport interpolation. Given paired samples $(x_0, x_1)$ with $x_0 \sim p_0$ and $x_1 \sim \mu_{\mathrm{data}}$, and the linear interpolant $x_t = (1-t)x_0 + t\,x_1$, the target velocity is simply $v^*_t = x_1 - x_0$. The CFM training objective takes the form
\begin{equation}
\calL_{\mathrm{CFM}}(\theta) = \E_{t \sim \mathcal{U}(0,1)}\, \E_{x_0 \sim p_0,\, x_1 \sim \mu_{\mathrm{data}}} \bigl\| v_\theta(x_t, t) - v^*_t \bigr\|^2.
\label{eq:cfm_loss}
\end{equation}
Unlike diffusion models that rely on stochastic noise schedules and iterative denoising, flow matching learns a deterministic vector field enabling direct, efficient sampling via ODE integration.

\subsubsection*{Auxiliary Geometric Penalties.}
Pure distribution matching can under-penalize rare but catastrophic constraint violations. Following the principle of fusing geometric inductive bias into training \cite{battaglia2018relational}, we augment the CFM loss with a geometric-aware penalty. For packing problems, we define an overlap energy on the projected data endpoint:
\begin{equation}
E_{\mathrm{overlap}}(x) = \sum_{i < j} \psi\bigl( 2r - \|x_i - x_j\| \bigr), \quad \psi(u) = \frac{1}{\beta} \log\bigl(1 + e^{\beta u}\bigr),
\label{eq:overlap_penalty}
\end{equation}
where $r$ is the object radius and $\psi$ is a \textsc{Softplus} function penalizing violations. Given the predicted velocity $v_\theta(x_t, t)$ and the affine path coefficients $(\alpha_t, \sigma_t)$ defining the interpolation $x_t = \alpha_t x_1 + \sigma_t x_0$ (with $\alpha_t = t$, $\sigma_t = 1-t$ for optimal transport), we project back to the data endpoint via $\hat{x}1 = (v\theta(x_t, t) - \dot{\sigma}_t x_0)/\dot{\alpha}t = v\theta(x_t, t) + x_0$.
The total training loss becomes
\begin{equation}
\calL(\theta) = \calL_{\mathrm{CFM}}(\theta) + \lambda(t, \mathrm{epoch}) \cdot E_{\mathrm{overlap}}(\hat{x}_1),
\label{eq:total_loss}
\end{equation}
with $\lambda$ following a cosine warm-up schedule that prioritises distribution matching early in training and constraint satisfaction thereafter.

\subsection{Geometry-Aware Sampling (GAS)}
Standard ODE integration of the learned flow~\eqref{eq:flow_ode} can drift off the feasible set and produce samples violating hard geometric constraints \cite{wang2022and,chen2023flow}, particularly for packing problems where non-overlap is essential. We introduce a \emph{Geometry-Aware Sampling} (GAS) algorithm that interleaves flow integration with constraint projection, analogous to the use of molecular dynamics corrections and energy minimization in protein structure refinement \cite{zhou2025energy}. GAS is a tangent-space projection-corrected sampler that pushes learned transport with explicit constraint correction during inference. Conceptually, GAS can be viewed as a geometric specialization of physics-constrained sampling \cite{utkarsh2025physics,cheng2024gradient}, where in our setting we target a wide range of geometric feasibility hypersurfaces.

In this section, we describe the GAS model specifically for the sphere packing problem. Details for the other problems studied in this paper are provided in Section~\ref{sec:results}.

Let $\H \subset \X$ denote the feasible set defined by non-overlap constraints. Rather than fixing boundary margins to the sphere radius $r$, we employ an adaptive wall constraint that depends on the configuration's minimum separation:
\begin{equation}
\H = \Bigl\{ x \in \X : \|x_i - x_j\| \geq 2r \forall\, i < j, \quad m(x) \leq x_i^{(k)} \leq L - m(x) \forall\, i, k \Bigr\},
\label{eq:feasible_set}
\end{equation}
where $m(x) = \tfrac{1}{2} \min_{i < j} \|x_i - x_j\|$ is half the minimum pairwise distance. This adaptive formulation ensures that wall margins scale with the effective radius supported by the current configuration, enabling the sampling process to discover tighter packings without imposing a fixed radius constraint during sampling. GAS maintains approximate feasibility throughout the sampling trajectory via iterative projection, without requiring gradient information for the constraints during training.

The algorithm is the following. We discretize the flow time $\tau \in [0,1]$ into $K$ steps with $\tau_k = k/K$, where $\tau = 0$ corresponds to the prior distribution and $\tau = 1$ to the data distribution. The model's internal time coordinate is $t = 1 - \tau$, so integration proceeds from $t=1$ toward $t=0$. At each step:

\begin{enumerate}[label=(\roman*), nosep, leftmargin=*]
\item \emph{Flow Integration:} Advance the state via ODE integration using an adaptive midpoint or higher-order solver:
\begin{equation}
\tilde{x}_{k+1} = x_k + \int_{\tau_k}^{\tau_{k+1}} v_\theta(x, 1 - \tau)  \mathrm{d}\tau.
\label{eq:ode_step}
\end{equation}
After each integration step, we enforce the adaptive box constraint $x \in [m(x), L - m(x)]^d$ via iterative reflection, which preserves the minimum pairwise separation while preventing boundary violations.

\item \emph{Gauss--Newton Projection:} Project onto the constraint manifold via linearised constraint enforcement. For active constraint pairs $(i,j)$ with overlap $h_{ij} = 2r - \|x_i - x_j\| > 0$, let $J$ denote the Jacobian of the constraint residuals $h$. The projection step computes
\begin{equation}
\Delta x = J^\top(J J^\top)^{-1} h,
\label{eq:gn_step}
\end{equation}
which we approximate efficiently via degree-normalized pairwise updates:
\begin{equation}
\Delta x_i = \frac{1}{\deg(i)} \sum_{j : h_{ij} > 0} \frac{h_{ij}}{2} \cdot \hat{n}_{ij}, \quad \hat{n}_{ij} = \frac{x_i - x_j}{\|x_i - x_j\|}.
\label{eq:contact_update}
\end{equation}
Boundary violations are handled analogously with wall-normal corrections.

\item \emph{Proximal Relaxation:} To prevent projection from deviating excessively from the learned flow manifold, we apply a proximal correction that balances constraint satisfaction with fidelity to the flow trajectory. Given the anchor point $\hat{x} = (1 - \tau') x_0 + \tau' x_{\mathrm{proj}}$ on the optimal transport path, we iteratively minimise
\begin{equation}
x \leftarrow x - \eta \Bigl( (x - \hat{x}) + \lambda_{\mathrm{prox}}\, J^\top h(x_{\mathrm{next}}) \Bigr),
\label{eq:prox_step}
\end{equation}
where $x_{\mathrm{next}} = x + (1-\tau') v_\theta(x, 1-\tau')$ is the forward-projected terminal state, $h$ is the constraint residual, and $J$ is its Jacobian. This couples the proximal anchor with a lookahead penalty on constraint violations at the predicted endpoint.

\item \emph{Terminal Refinement:} After completing flow integration, we apply additional projection passes with tightened tolerances to ensure strict feasibility. The algorithm terminates when the maximum overlap residual falls below a prescribed tolerance:
\begin{equation}
\max_{i < j} \bigl( 2r - \|x_i - x_j\| \bigr)^+ \leq \epsilon_{\mathrm{tol}},
\label{eq:terminal_tol}
\end{equation}
where $(\cdot)^+ = \max(\cdot, 0)$ denotes the positive part. This criterion is independent of any fixed radius assumption, yielding samples that are both distributionally faithful and geometrically valid for the effective radius $r_{\mathrm{eff}} = \tfrac{1}{2}\min_{i<j}\|x_i - x_j\|$ supported by the configuration.
\end{enumerate}

This projection-relaxation scheme is conceptually similar to constraint satisfaction in molecular dynamics \cite{nam2024flow}, where forces from bonded interactions (analogous to our flow field) are balanced against non-bonded repulsions (our overlap penalties).

\subsection{Reward-Guided Fine-Tuning}
While the base CFM model learns to sample from the distribution of training configurations, the discovery of extremal structures requires \emph{improving} upon this distribution towards OOD regions, generating configurations with higher objective values than any in the training set. In other words, rare structure discovery, is intrinsically extrapolative. This challenge is very similar to de novo molecular design, where one seeks to extrapolatively generate molecules with properties exceeding those in known databases. As a result, we introduce \emph{Reward-Guided CFM} (RG-CFM), an online fine-tuning procedure that synthesizes: (i) importance-weighted policy optimization from reinforcement learning, and (ii) consistency regularization from self-supervised learning (SSL) \cite{balestriero2023cookbook}. This synthesis addresses a fundamental in pure reward-guided generation \cite{bengio2021flow}. Pure reward maximization drives the model toward high-scoring configurations, but without regularization it collapses to a degenerate policy that produces only a single output, losing the diversity necessary for continued exploration and discovery.

\subsubsection*{Reward Function.}
For each problem, we define a reward $R \colon \X \to \R$ aligned with the optimization objective. For sphere packing, we use the effective radius, half the minimal separation distance, normalized by the box length:
\begin{equation}
R(x) = \frac{r_{\mathrm{eff}}(x)}{L}, \quad r_{\mathrm{eff}}(x) = \frac{1}{2} \min_{i < j} \|x_i - x_j\|,
\label{eq:reward}
\end{equation}
which directly measures packing quality: larger $R$ supports larger sphere radii. This formulation is agnostic to any predetermined target radius, allowing the model to discover configurations that exceed the quality of training data.

\subsubsection*{Online Reward-Guided Flow Matching}
RG-CFM adopts a teacher-student framework inspired by self-supervised methods such as DINO \cite{caron2021emerging} and Mean Teacher \cite{tarvainen2017mean}. We maintain two copies of the velocity field, a \emph{student} model $v_\theta$ that is actively fine-tuned, and a \emph{teacher} model $v_{\mathrm{ref}}$ that remains frozen at the pretrained weights. The teacher serves as an anchor, providing consistency targets that prevent the student from drifting into degenerate solutions, a phenomenon we term \emph{generative collapse}, analogous to representation collapse in contrastive learning \cite{jing2021understanding}. To push generation toward high-reward regions, we adopt importance weighting from offline RL \cite{peng2019advantage,nair2020awac}. If $q(x)$ denotes the distribution of training endpoints and $w(x) \geq 0$ is a weighting function, then minimizing a weighted flow matching loss yields a learned terminal distribution proportional to $w(x)\,q(x)$. Repeated cycles of sampling from the current model and retraining with the same weighting amplify mass by successive powers of $w$, progressively concentrating the distribution around reward-maximizing modes. This observation motivates both (i) reward-weighted updates for systematic improvement, and (ii) explicit regularization to prevent degeneracy.

Concretely, given a batch of samples $\{x^{(b)}\}_{b=1}^B$ generated by the current student model via GAS, along with their reward scores $\{R^{(b)}\}$, we compute importance weights through $z$-score normalization followed by exponentiation:
\begin{equation}
w^{(b)} = \frac{1}{Z} \exp\Bigl( \tau \cdot \frac{R^{(b)} - \bar{R}}{\sigma_R + \epsilon} \Bigr), \quad Z = \frac{1}{B}\sum_{b'} \exp\Bigl( \tau \cdot \frac{R^{(b')} - \bar{R}}{\sigma_R + \epsilon} \Bigr),
\label{eq:reward_weights}
\end{equation}
where $\bar{R}$ and $\sigma_R$ are the batch mean and standard deviation of rewards, $\tau > 0$ is a temperature controlling selection pressure, and weights are clipped to $[0, w_{\max}]$ for numerical stability. This scheme ensures that samples with above-average rewards receive amplified influence during training, while below-average samples are downweighted but not discarded entirely, reminiscent of Advantage-Weighted Regression (AWR), where policy updates take the form of weighted maximum likelihood with weights exponential in the advantage.

Without regularization, pure reward maximization leads to a well-known failure mode. Given sufficient training, the student model converges to a degenerate policy that produces only a single high-reward output, losing the diversity necessary for continued exploration. Each retraining round amplifies probability mass on the current best modes, and the model eventually converges to a delta distribution concentrated on one or a few configurations. This generative collapse is catastrophic for optimization: once the model produces only one configuration, it cannot explore alternative basins, and the search terminates prematurely at a local optimum. We prevent collapse through \emph{consistency regularization}, penalising deviation of the student's velocity field from the frozen teacher's predictions. At each training step, both models are evaluated on the same interpolated states $x_t$, and we minimize their squared discrepancy:
\begin{equation}
\calL_{\mathrm{consist}} = \frac{1}{B} \sum_{b=1}^B \bigl\| v_\theta(x_t^{(b)}, t) - v_{\mathrm{ref}}(x_t^{(b)}, t) \bigr\|^2.
\label{eq:consistency_reg}
\end{equation}
This term acts as a soft constraint in function space, bounding how far the fine-tuned velocity field can deviate from the pretrained baseline. The mechanism mirrors the role of the exponential moving average (EMA) teacher in DINO \cite{caron2021emerging} and BYOL \cite{grill2020bootstrap}, which provides stable targets that anchor the student's learning and maintain output diversity through implicit regularization.

As a result, the total RG-CFM loss combines reward-weighted flow matching with consistency regularization:
\begin{equation}
\calL_{\mathrm{RG}} = \calL_{\mathrm{FM}}^{\mathrm{w}} + \alpha \cdot \calL_{\mathrm{consist}},
\label{eq:rg_total}
\end{equation}
where $\alpha > 0$ controls the exploration-exploitation trade-off. Small $\alpha$ permits aggressive reward optimization at the risk of collapse and large $\alpha$ preserves diversity but limits improvement over the baseline. In practice, we find that $\alpha \in [0.1, 1.0]$ provides a robust operating regime, analogous to the KL penalty coefficient in proximal policy Optimisation (PPO) \cite{schulman2017proximal} or the trust region in TRPO \cite{schulman2015trust} methods where the new policy is the old policy reweighted by exponentiated advantages, regularized to stay close to a reference and RLHF \cite{christiano2017deep} for LLMs.

\subsubsection*{Geometry-Aware Action Exploration.}
An important limitation of pure reward-guidance is support sub-optimality, meaning that weighted training cannot assign probability mass to configurations that are essentially absent from the distribution used to generate training endpoints. In other words, if the current policy (or current elite dataset) has negligible density in a region containing better solutions, reweighting alone cannot reliably discover it. This motivated us to have an explicit exploration operators that propose structured perturbations beyond the current support. Moreover, naive continuous-action online RL methods encourage exploration by sampling actions from a normal distribution when collecting episodes, but simply adding noise to the action trajectories produces non-smooth trajectories and does not efficiently explore the action space. Thus, to further enhance sample diversity during fine-tuning, during on-policy data collection for RG-CFM  we introduce an exploration operator $\mathcal{E}(x'|x)$ that proposes smooth, constraint-informed moves:
\begin{equation}
\label{eq:explore}
x' = x + \Delta x,
\qquad
\Delta x = (M r)\, s(x)\, \widehat{d}(x),
\qquad
\widehat{d}(x) = \frac{c\,\Delta_{\mathrm{contact}}(x) + (1-c)\,\Delta_{\mathrm{wall}}(x)}{\|c\,\Delta_{\mathrm{contact}}(x) + (1-c)\,\Delta_{\mathrm{wall}}(x)\|_{\mathrm{F}}+\varepsilon}.
\end{equation}
Here $\|\cdot\|_{\mathrm{F}}$ is the Frobenius norm on $\R^{d\times N}$, $M$ is a dimensionless exploration magnitude, $c\in[0,1]$ mixes contact- and wall-driven directions, and $s(x)$ is an overlap-severity scalar (set to $1$ by default, and increased when overlaps are large in order to push more strongly out of infeasible tangles). $\Delta_{\mathrm{contact}}$ is computed from the contact graph (directions that relieve overlaps), $\Delta_{\mathrm{wall}}$ from boundary proximity, and $M$ controls the exploration magnitude. After exploration, a repair step projects the perturbed configuration back onto the feasible manifold, ensuring that rewards are evaluated on geometrically valid samples. 
This exploration mechanism plays the same conceptual role as guided MCMC moves in molecular sampling \cite{plainer2023transition} that prevent uniform collapse, it explicitly expands the effective support of the policy so that reward-weighted updates can amplify genuinely novel, higher-reward structures once they are encountered.

\subsection{Closing the Loop: Convergence Properties}
A critical limitation of prior SBO system for mathematical structure discovery, notably PatternBoost \cite{charton2024patternboost}, is that the boosting loop remains \emph{open}: there is no mechanism guaranteeing convergence to optimal or near-optimal solutions. In open-loop methods, the generative model learns to mimic the current best solutions but receives no direct signal pushing it toward \emph{better} and OOD solutions, it matches a distribution rather than optimizing an objective. We have addressed this fundamental gap through online reward-guided policy optimization, where rather than treating the generative model as a passive density estimator retrained on filtered samples, we fine-tune it with direct feedback from the optimization objective. This transforms the open iterative scheme into a closed-loop control system:

\begin{itemize}[nosep, leftmargin=*]
\item \emph{Open loop:} The generative model approximates $p_{\mathrm{data}}$, the distribution of current best solutions. Improvement relies on sampling and local search occasionally escaping to better basins, an indirect, stochastic process with no gradient signal toward the objective that prompts to a lengthy and blind search.

\item \emph{Closed loop (ours):} The generative model is fine-tuned to maximize expected reward $\E_{x \sim p_\theta}[R(x)]$ while remaining close to a reference policy. This provides a direct optimization signal, where configurations with higher objective values receive higher weight in the policy update, systematically guiding sampling toward better and OOD solutions.
\end{itemize}

The closure mechanism operates through the reward-guided flow matching objective where samples from the current policy are evaluated, and the model is updated to increase the likelihood of high-reward trajectories while the consistency regularization and action exploration prevents mode collapse. Combined with geometry-aware sampling, which maintains geometric feasibility throughout the generative process, this closed-loop formulation enables convergence to near-optimal solutions in much fewer iterations, $\mathcal{O}(1\text{--}10)$ boosting rounds rather than the $\mathcal{O}(100)$ iterations characteristic of open-loop methods.

\subsection{The \FMBoost{} Pipeline}

The complete \FMBoost{} pipeline integrates the above components in an iterative boosting loop:

\begin{enumerate}[label=(\roman*), nosep, leftmargin=*]
\item \emph{Initialization.} Generate an initial dataset $\mathcal{D}_0$ of configurations by running stochastic relaxation with perturbations (SRP) from random initial conditions. Retain the top fraction (typically 25--50\%) according to the objective $J$.

\item \emph{Training.} Train a conditional flow matching model $v_\theta^{(0)}$ on $\mathcal{D}_0$ using the combined flow matching and geometric penalty objectives.

\item \emph{Sampling.} Sample a batch of configurations from $v_\theta^{(0)}$ using geometry-aware sampling with projection-corrected integration.

\item \emph{Reward-Guided Fine-Tuning.} Evaluate rewards on the refined samples. Apply reward-guided fine-tuning with consistency regularization and geometry-aware action exploration to update the model parameters.

\item \emph{Refinement \& Selection.} Apply local search (SRP followed by L-BFGS optimization) to each generated configuration, obtaining a refined batch.

\item \emph{Iteration.} Repeat steps 1--5, obtaining successive model updates $v_\theta^{(1)}, v_\theta^{(2)}, \ldots$ with progressively improving sample quality.
\end{enumerate}

Empirically, we observe that in each iteration the distribution of $J$ over the samples shifts upward, and the best configurations improve over time. The closed-loop structure ensures systematic improvement rather than stochastic drift, with convergence typically achieved within 1--3 boosting rounds. 

\remark{Note that with the closed-loop pipeline, we \emph{do not} require the extremely long outer boosting iterations of PatternBoost. In such open-loop schemes, the generator is retrained only to imitate the current elite set; improvement is therefore indirect and essentially stochastic as one must repeatedly sample and re-run local search in the hope of occasionally landing in a better basin. This makes progress heavily dependent on hundreds to thousands of outer iterations and offers no principled mechanism for systematic ascent. By contrast, our \FMBoost\ closes the loop, meaning that rewards computed on newly generated (and repaired) configurations directly drive an online update of the generator via reward-guided fine-tuning. As a result, the optimization signal is injected where it matters, namely into the parameters of the generative policy, so the dominant driver of improvement becomes a small number of targeted inner updates rather than a large number of blind outer retrain cycles.}

\subsection{Local Search: Stochastic Relaxation with Perturbations}

The local search component employs a variant of stochastic relaxation with perturbations (SRP). Given an initial configuration $x \in \X$, we consider a smooth surrogate objective $\widetilde{J}(x) \approx J(x)$ that is differentiable in the ambient Euclidean coordinates. For packing problems, this comprises an overlap energy penalizing sphere intersections and boundary violations, combined with terms proportional to the negative sum of radii.

SRP alternates between two phases: (i) \emph{random perturbation}, adding small random displacements to the configuration, and (ii) \emph{gradient relaxation}, performing several steps of normalized gradient descent on $\widetilde{J}$, often with an annealing schedule for the step size. For packing problems, SRP is followed by L-BFGS-B optimization over the same surrogate to polish the configuration while respecting box constraints.

The key properties of SRP are that it is computationally inexpensive and parallel, robust to initial conditions due to the perturbation and annealing, and captures fine geometric features that may be difficult for a global model to learn purely from data. In \FMBoost{}, we use SRP in two modes: as a training-data generator starting from random configurations, and as a final refinement step on samples produced by the flow model.

In continuous optimization settings, \emph{local search} simply refers to a fast \emph{local fixing} (or local improvement) procedure: starting from an initial configuration, it repeatedly applies small, inexpensive updates that reduce constraint violations and improve the objective, until it reaches a locally stable configuration. In practice, one often has several plausible local search heuristics that all produce near-feasible, near-best-known configurations, but with markedly different reliability and final quality. This choice is crucial: since our goal is to surpass the best known constructions, the local search mechanism must produce samples that are as strong as possible, both to supply high-quality training data and to provide an effective refinement step for generative samples.

\begin{figure}[t]
    \centering
    \includegraphics[width=0.8\linewidth]{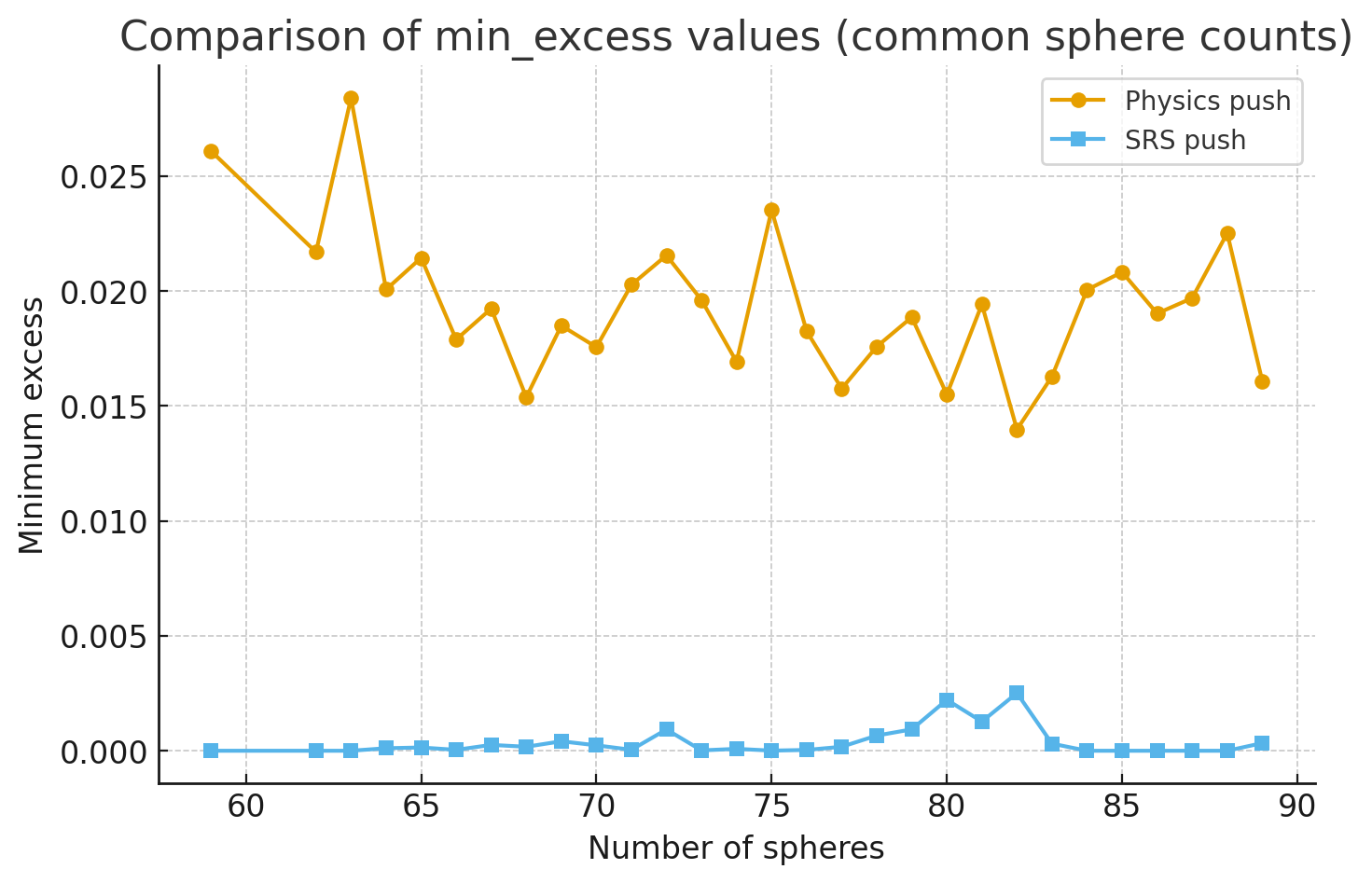}
    \caption{Comparison of the minimum-excess metric for two local search heuristics (common circle counts). The physics-push heuristic yields consistently higher minimum-excess values than SRP/SRS, indicating worse configurations across the tested regime.}
    \label{fig:physics_push_vs_srp}
\end{figure}

Figure \ref{fig:physics_push_vs_srp} illustrates how different geometric heuristics can explore different regions of configuration space and yield radically different solution quality. In this example (circle packing in the unit square), a simple \emph{physics-push} baseline (where each center is iteratively displaced along a repulsion direction given by a weighted sum of vectors pointing away from neighbouring centers, with weights proportional to the reciprocal of the inter-center distances) consistently produces substantially worse configurations than SRP across common circle counts.

\section{Results}\label{sec:results} 
This section contains a detailed description of the results of our \FMBoost{} experiments on geometric optimization problems. 

\subsection{Overview}

\subsubsection*{Sphere Packing}
We evaluate \FMBoost{} on the classical problem of packing $N$ non-overlapping spheres of radius $r$ inside a $d$-dimensional unit hypercube, a problem that combines geometric rigidity with combinatorial complexity. The objective is to maximize $r$ (equivalently, the packing fraction $\phi = N \cdot \vol_d(B_r)$) subject to non-overlap constraints $\|x_i - x_j\| \geq 2r$ for all pairs and containment constraints $x_i \in [r, 1-r]^d$ for all centres. This problem admits no closed-form solution for general $N$; even verifying local optimality is NP-hard, and the best known packings for most $N$ have been found through extensive computational search. Our investigation of high-dimensional packings is motivated by recent breakthroughs in the asymptotic theory. Viazovska~\cite{viazovska2017sphere} proved that the $E_8$ lattice achieves the optimal sphere packing density in $\mathbb{R}^8$, and together with collaborators~\cite{cohn2017sphere} established optimality of the Leech lattice in $\mathbb{R}^{24}$. These remain the only dimensions $d > 3$ where the densest packing is known. Dimension $d = 12$ is particularly interesting as it lies between the solved cases, admits rich lattice structure \cite{conway1983coxeter} (e.g., the Coxeter-Todd lattice $K_{12}$), yet the optimal packing density remains unknown, making it a natural testbed for computational exploration.
We report experiments in two regimes: (i) three-dimensional packings for $N \in \{50, \ldots, 89\}$ and $N \in \{191, 200\}$, where high-quality reference packings exist from decades of prior work; and (ii) twelve-dimensional packings for $N \in 31$, where the geometry is far from lattice-like and classical heuristics struggle. In several instances, \FMBoost{} matches or exceeds the best previously reported packing fractions, and in high-dimensional settings it consistently discovers configurations denser than those produced by simulated annealing or basin-hopping alone.

\subsubsection*{The Heilbronn Problem} 
We study the classical Heilbronn triangle problem in the unit square: for a point set
$X=\{p_1,\dots,p_n\}\subset[0,1]^2$ one considers the minimum area over all triangles spanned by triples,
\[
A_{\min}(X) := \min_{1\le i<j<k\le n}\frac12\bigl|\det(p_j-p_i,p_k-p_i)\bigr|,
\]
and the goal is to maximize $A_{\min}(X)$ over all $n$-point configurations.
In our \FMBoost{} installation, SRP produces an initial elite set by maximizing a smooth soft-min surrogate of $A_{\min}$, followed by a deterministic L-BFGS-B polish and an active max-min push (SLSQP on the currently smallest-area triangles).
A conditional flow-matching model is then trained on the top portion of these SRP-refined configurations and used to propose new point sets; each generated candidate is again refined by the same SRP final push before selection.
Empirically, the raw generator underperforms the SRP elite distribution, but the final push reliably recovers near-elite configurations, and the closed-loop iteration can exceed the previous best.
In our runs we obtain $A_{\min}=0.0259285$ for $n=13$ (beating the training maximum in that experiment) and $A_{\min}=0.0187494$ for $n=15$, with clear right-shifts in the distributions of $A_{\min}$ across \FMBoost{} iterations. These results approach the best known constructions in \cite{FriedmanHeilbronnSquaresBenchmark} $A_{\min} = 0.0270$ for $n=13$ and $A_{\min} = 0.0211$ for $n=15$.

\subsubsection*{Circle packing with maximal sum of radii}

We study circle packings in the unit square where the objective is to maximize the \emph{sum of radii} for a fixed number $N$ of circles, subject to non-overlap and boundary constraints.  This problem has been recently attacked by AlphaEvolve and its open-source variants using LLM-based genetic programming to evolve code that produces good packings. In our approach, SRP generates an initial dataset of circle arrangements, and a flow-matching model on centers learns to propose new center configurations.  The final push stage solves a convex optimization problem (via linear programming) to maximize the sum of radii subject to constraints for the given centres, and combines this with an SRP-based refinement. For $N=26$ and $N=32$ circles, we obtain configurations whose sum of radii strictly exceeds the best values reported for AlphaEvolve-type methods, while using a considerably smaller computational budget.

\subsubsection*{Star discrepancy problem}

We also apply \FMBoost{} to constructing low star-discrepancy point sets in $[0,1]^2$, a central objective in quasi-Monte Carlo integration \cite{ClementDoerrKlamrothPaquete2025OptimalStarDiscrepancy,ALGOLabStarDiscrepancyBenchmark,ClementDoerrPaquete2022StarSubsetSelection}. For $P=\{p_1,\dots,p_N\}\subset[0,1]^2$, the (anchored) star discrepancy is the minimax quantity
\begin{equation}\label{def:star_discrepancy}
D^\ast(P)
:=
\sup_{(a,b)\in[0,1]^2}
\left|
\frac{1}{N}\#\bigl(P\cap([0,a)\times[0,b))\bigr)-ab
\right|,
\qquad \text{(smaller is better),}
\end{equation}
i.e.\ the worst deviation between empirical mass and volume over anchored axis-aligned boxes.

Our local SRP search minimizes a differentiable soft-max surrogate of $D^\ast$ based on sigmoid-smoothed box indicators and a log-sum-exp over a grid of anchored boxes; this is followed by L-BFGS-B and exact evaluation of $D^\ast$ on the critical grid.
A conditional flow-matching model (conditioned on $(N, D^\ast)$) is then trained on the best half of the pushed samples, and sampling interleaves ODE integration with short projection/proximal steps that directly decrease the smooth discrepancy surrogate before the final push.
As with Heilbronn, generation alone is not sufficient, but the final push recovers and slightly improves the best tail.
In our experiments we obtain $D^\ast=0.06290897$ for $N=20$ and $D^\ast=0.02943972$ for $N=60$, and the pushed distributions consistently shift left relative to the initial SRP training set.

\subsection{Sphere Packing in Hypercube}
\label{sec:sphere_packing}
\subsubsection*{Model architecture.}
The velocity field $v_\theta(x, t)$ is parameterized by a permutation-equivariant Transformer \cite{vaswani2017attention,lee2019set} operating on the point-cloud $x = (x_1, \ldots, x_N) \in (\R^d)^N$, where each configuration $x \in \R^{d \times N}$ is treated as a set of $N$ tokens with $d$-dimensional coordinates. Time information is embedded via random Fourier features combined with polynomial basis functions \cite{rahimi2007random}, $\phi(t) = \bigl[\sin(2\pi \omega_j t),\, \cos(2\pi \omega_j t),\, t,\, t^2,\, t^3,\, \log(1+t)\bigr]_{j=1}^{F},$
where $\{\omega_j\}$ are fixed random frequencies. This embedding is projected to a conditioning vector and injected into transformer layers via Feature-wise Linear Modulation (FiLM) \cite{perez2018film}: $h \mapsto h \odot (1 + \gamma) + \beta$, where $(\gamma, \beta)$ are predicted from the time-condition embedding. Our architecture further incorporates problem-specific conditioning variable $y=(r/L, N, \text{face-contact ratio}, \text{min-separation}/L)$, enabling a single to generalize across problem instances by learning a conditional probability density $p(\X|y)$. For packing in a box $[r,L-r]^d$, the base distribution is uniform in the feasible box but does not enforce non-overlap. As a result, we have a geometry-aware prior that improves coverage of boundary-touching structures by placing each point on a randomly chosen box, then adding small jitter. This prior is conditional on $y$ and improves sample efficiency in regimes where high-quality solutions exhibit frequent face contacts. Key hyperparameters are: model dimension $512$, depth $2$ layers, $8$ attention heads, and GLU \cite{shazeer2020glu} feed-forward blocks with RMSNorm \cite{zhang2019root}. The output layer is zero-initialized to ensure the initial velocity field is near-constant, stabilizing early training. Total parameter count is approximately $2$M.

\subsubsection*{Training.}
Initial datasets are generated via Stochastic Repulsion Push (SRP), a physics-inspired local search that eliminates overlaps through annealed gradient descent on a soft-overlap energy. We retain the top $5$--$10\%$ of configurations ranked by packing fraction. The flow-matching model is trained for $200$--$500$ epochs using the AdamW variant with schedule-free learning rate adaptation, batch size $128$, and gradient clipping at norm $1.0$. An auxiliary penalty encouraging large pairwise separations is ramped in over the first half of training via a cosine schedule. Time sampling is biased toward small $t$ (near-data regime) using a mixture distribution with weight $0.5$ on a power-law component ($\gamma = 2$).

\subsubsection*{Sampling.}
Geometry-Aware Sampling (GAS) integrates the learned velocity field from $\tau = 0$ (prior) to $\tau = 1$ (data) using a midpoint ODE solver with $40-60$ steps. At each step, a Gauss--Newton projection enforces non-overlap and wall constraints, followed by proximal relaxation that anchors the trajectory to the optimal-transport path. Terminal refinement applies additional projection passes until the maximum overlap residual falls below $10^{-8}$. The prior distribution places centers uniformly on box faces with small jitter, providing geometric diversity while respecting approximate containment. We then fine-tune the pretrained model using reward guidance. The teacher model (frozen at pretrained weights) provides consistency targets while the student model is updated via importance-weighted flow matching with temperature $\tau \in [0.5, 2.0]$ and consistency coefficient $\alpha \in [0.1, 0.5]$. Each fine-tuning epoch samples a batch of $128$--$256$ configurations from the current student policy via GAS, evaluates their effective radii $r_{\mathrm{eff}} = \tfrac{1}{2}\min_{i<j}\|x_i - x_j\|$, computes $z$-scored importance weights, and performs a ten gradient step. Fine-tuning typically converges within $2$--$10$ epochs, after which the best configurations are passed to SRP for final polishing. 
The required steps for a single GAS sample plus projection overhead, is comparable in wall-clock time to $\sim\!100$ L-BFGS-B iterations. However, the important things is that the generative model produces \emph{diverse} samples that explore multiple basins simultaneously, whereas local search must be restarted from scratch to escape a given basin. In practice, $10^3$ GAS samples explore the configuration space more thoroughly than $10^4$ SRP restarts, yielding both higher peak quality and better coverage of near-optimal configurations.

\subsubsection{Results}
We benchmark against the Packomania database \cite{packomania} , which compiles the best known packings for $N \leq 200$ spheres in the unit hypercube, obtained through decades of simulated annealing, genetic algorithms, and manual refinement.

\begin{figure}[t]
\centering
\includegraphics[width=0.92\linewidth]{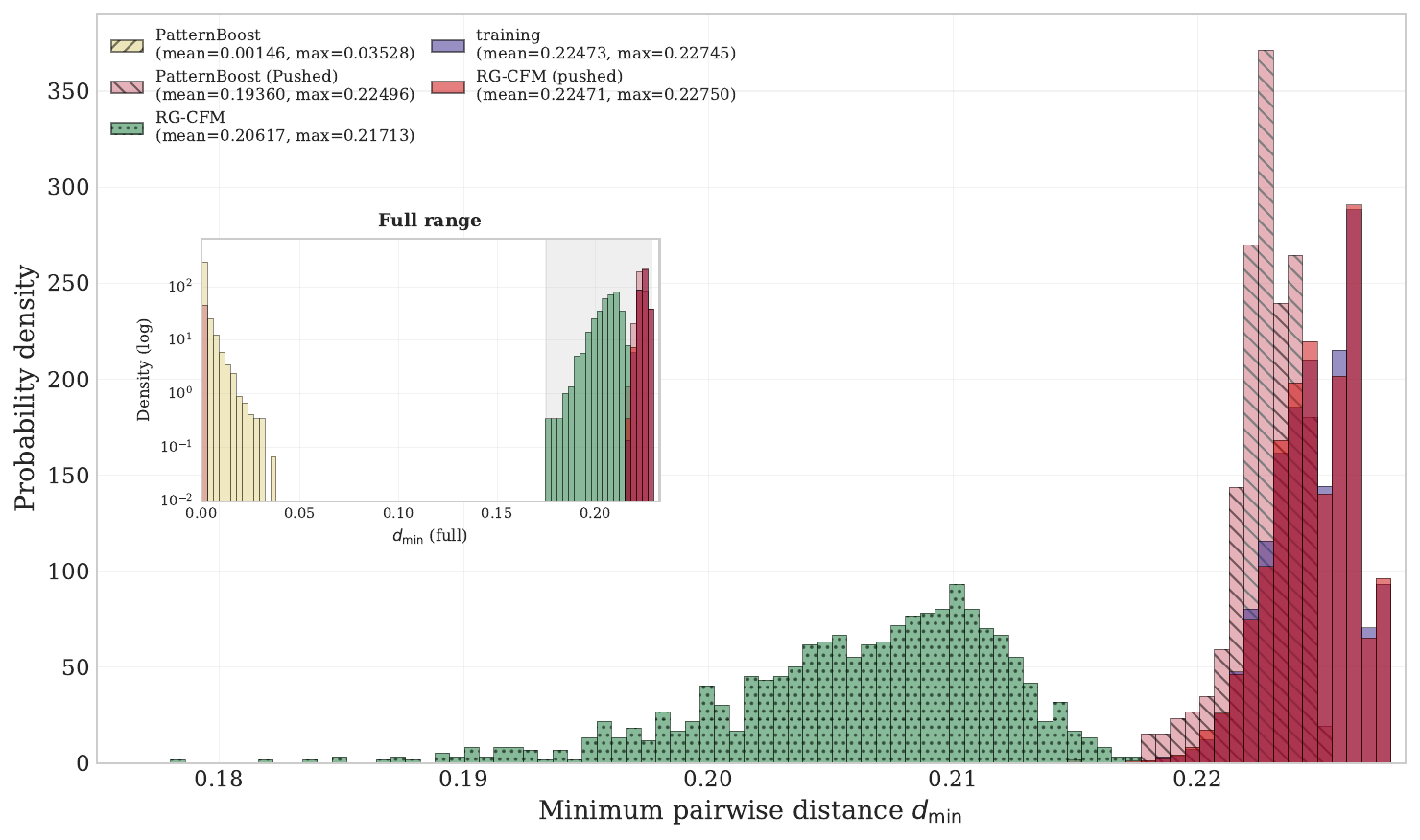}
\caption{\textbf{3D Sphere Packing in Cube, $n=89$.} The comparison between PatternBoost, RG-CFM, and their pushed version with the training data. Our Reward-guided fine tuning alone could provide the same performance as the expensive PatternBoost and SRP Push combined, while the pushed version of RG-CFM could improve the max over the training data.}
\label{fig:3d-89}
\end{figure}

For $N = 50--200$, we apply the \FMBoost{} pipeline: train on SRP-generated data, sample with GAS, fine-tune with reward guidance, and refine with SRP, and iterate. Across $1$--$5$ boosting rounds, the generative model progressively shifts probability mass toward higher-quality basins. We also provide a comparision against PatternBoost. Figure~\ref{fig:3d-89} directly compares \FMBoost{} against PatternBoost for $N=89$ spheres. Two findings stand out. First, RG-CFM \emph{without} local refinement achieves comparable peak quality to PatternBoost even \emph{with} SRP push demonstrating that reward-guided flow matching partially internalizes the refinement that PatternBoost delegates entirely to post-hoc search. Second, RG-CFM with push exceeds both the training maximum and PatternBoost's best result, while requiring less outer iterations and sampling time. This reduction in convergence time, combined with the elimination of LLM dependencies, constitutes the primary computational advantage of the closed-loop formulation.

\begin{table}[htbp]
\centering
\caption{\textbf{Computational comparison for $N=89$, $d=3$.} Wall-clock times on a single A100 GPU. FLOPs exclude local search.}
\label{tab:compute-comparison}
\smallskip
\begin{tabular}{@{}lcccc@{}}
\toprule
\textbf{Method} & \textbf{Parameters} & \textbf{Iterations} & \textbf{FLOPs} & \textbf{Wall-clock} \\
\midrule
PatternBoost & $\sim$8M & $\mathcal{O}(10\text{--}100)$ & $\sim 10^{12}$ & 10--100\,h \\
FlowBoost & $\sim$2M & $\mathcal{O}(1\text{--}10)$ & $\sim 10^{8}$ & 0.5--2\,h \\
\bottomrule
\end{tabular}
\end{table}

\subsubsection*{12D Sphere Packing} High-dimensional sphere packing presents qualitatively different challenges: the configuration space grows exponentially, local search becomes trapped in shallow basins, and there is no lattice structure to exploit for small $N$. We test \FMBoost{} on $N = 31$ spheres in the $12$-dimensional hypercube cube, a regime where even generating \emph{feasible} packings is nontrivial. In high-dimensional packing, any improvement over extensive local search is noteworthy. The training dataset comprises $\mathcal{O}(10^3)$ configurations generated by SRP with random restarts that is a strong baseline that has been refined over many iterations. That a single round of RG-CFM fine-tuning yields a configuration exceeding this entire dataset's maximum suggests that the learned surrogate captures geometric structure invisible to coordinate-wise local search.

\begin{figure}[htbp]
    \centering
    \begin{subfigure}[t]{0.48\textwidth}
        \centering
        \includegraphics[width=\textwidth]{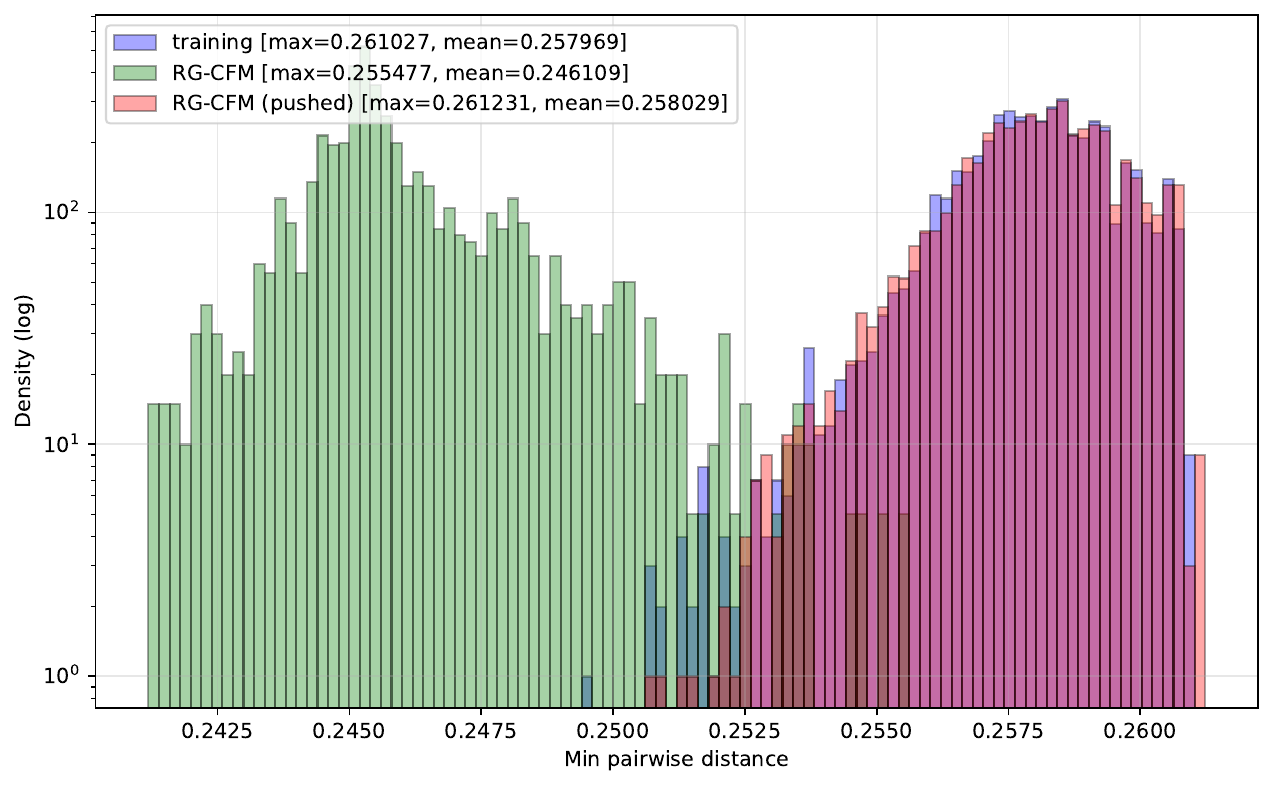}
        \caption{$N = 55$: training vs.\ RG-CFM (pushed).}
        \label{fig:3d-55}
    \end{subfigure}
    \hfill
    \begin{subfigure}[t]{0.48\textwidth}
        \centering
        \includegraphics[width=\textwidth]{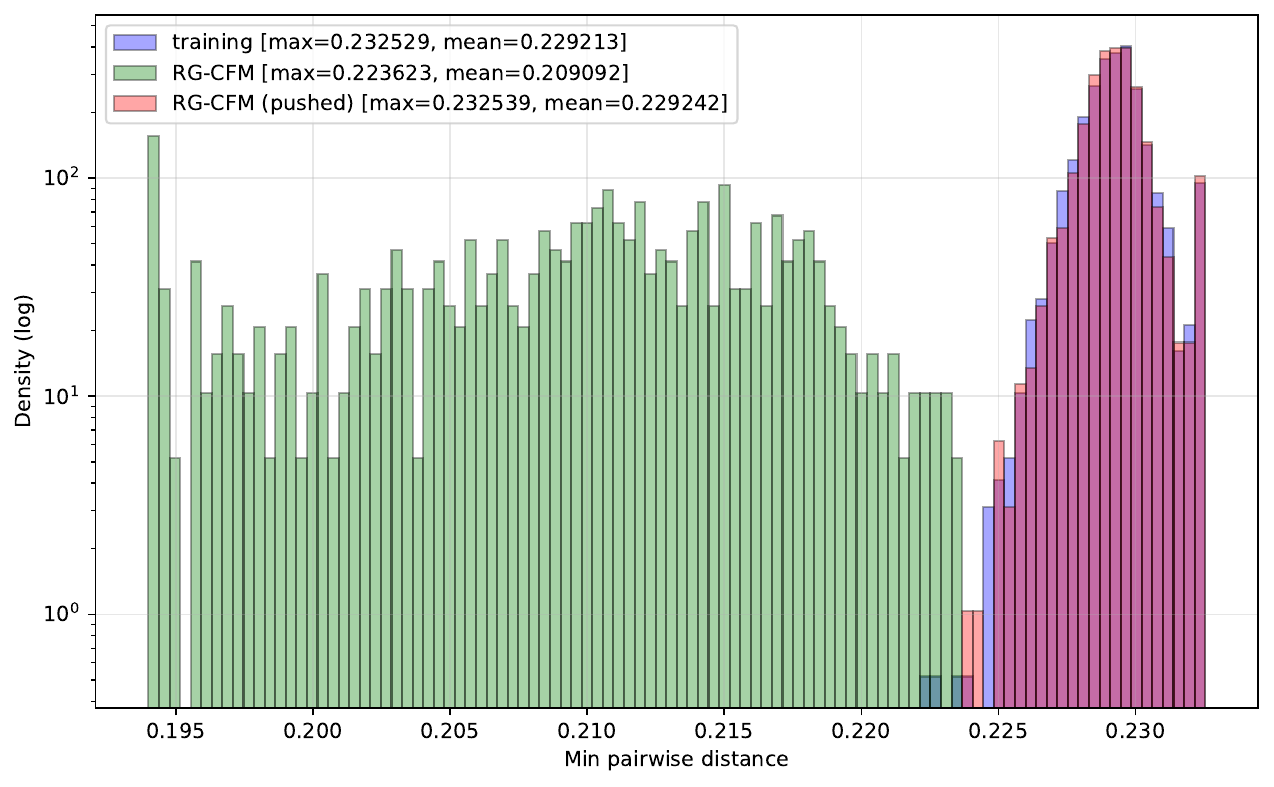}
        \caption{$N = 83$: training vs.\ RG-CFM (pushed).}
        \label{fig:3d-83}
    \end{subfigure}
    
    \vspace{0.5em}
    
    \begin{subfigure}[t]{0.48\textwidth}
        \centering
        \includegraphics[width=\textwidth]{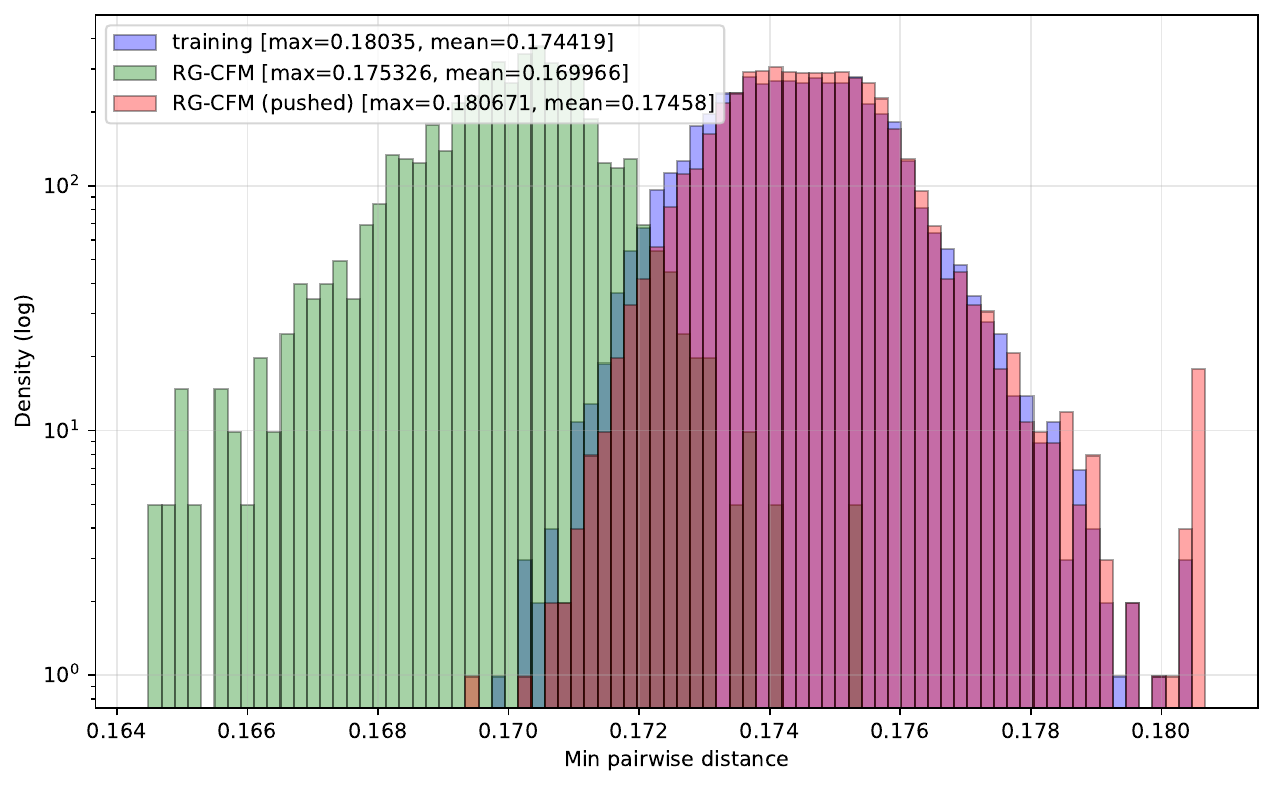}
        \caption{$N = 191$: training vs.\ RG-CFM (pushed).}
        \label{fig:3d-191}
    \end{subfigure}
    \hfill
    \begin{subfigure}[t]{0.48\textwidth}
        \centering
        \includegraphics[width=\textwidth]{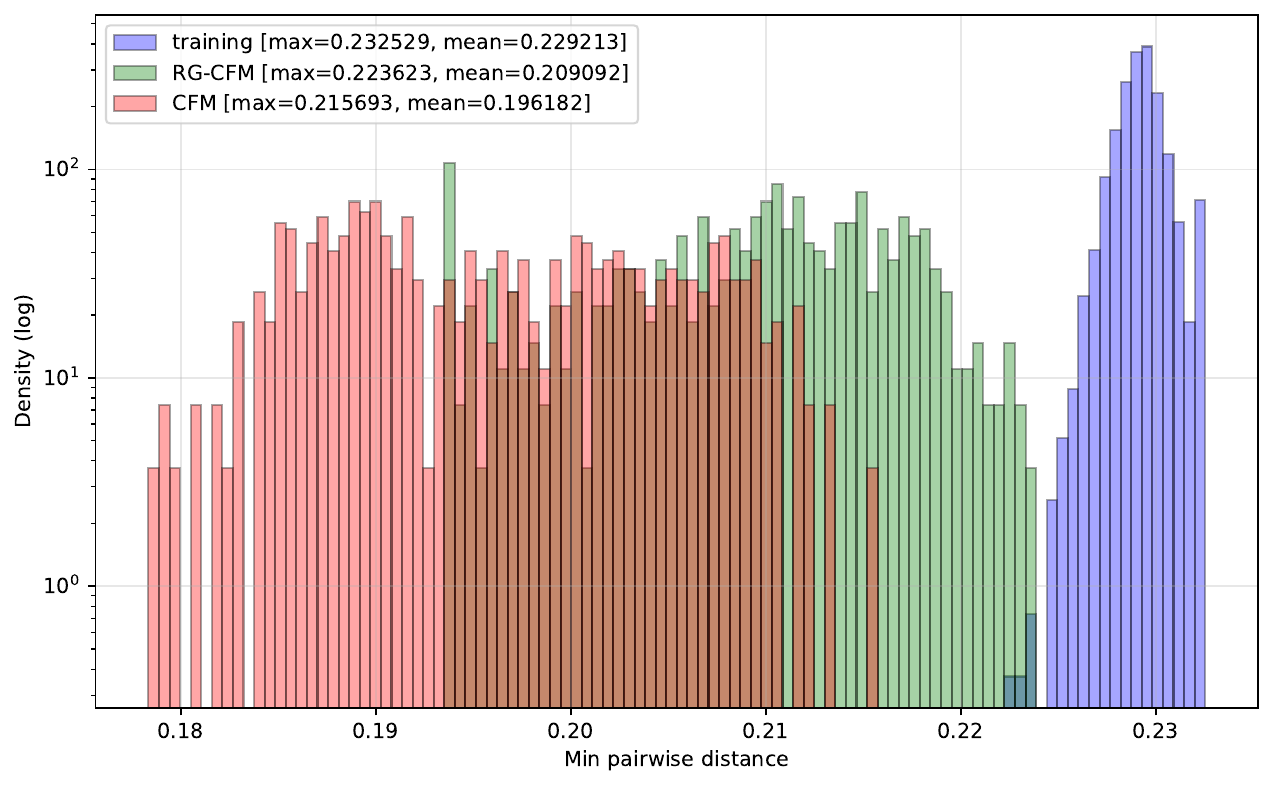}
        \caption{$N = 83$: ablation of reward guidance.}
        \label{fig:3d-83-ablation}
    \end{subfigure}
    \caption{\textbf{Sphere packing in $d=3$.} (a–c) Normalized histograms of minimum pairwise distance (log scale). In all cases, RG-CFM with final push (blue) recovers the training distribution and yields configurations exceeding the training maximum: $d_{\min} = 0.261231$ vs.\ $0.261027$ for $N=55$; $0.232539$ vs.\ $0.232529$ for $N=83$; $0.180671$ vs.\ $0.180350$ for $N=191$. (d)~Ablation comparing vanilla CFM (orange) against reward-guided CFM (green) for $N=83$: reward guidance with action exploration shifts the distribution toward higher-quality samples ($d_{\min}^{\max} = 0.2236$ vs.\ $0.2157$), demonstrating that direct objective feedback improves generation quality prior to local refinement.}
    \label{fig:3d-packing}
\end{figure}

\begin{figure}[htbp]
    \centering
    \begin{subfigure}[t]{0.48\textwidth}
        \centering
        \includegraphics[width=\textwidth]{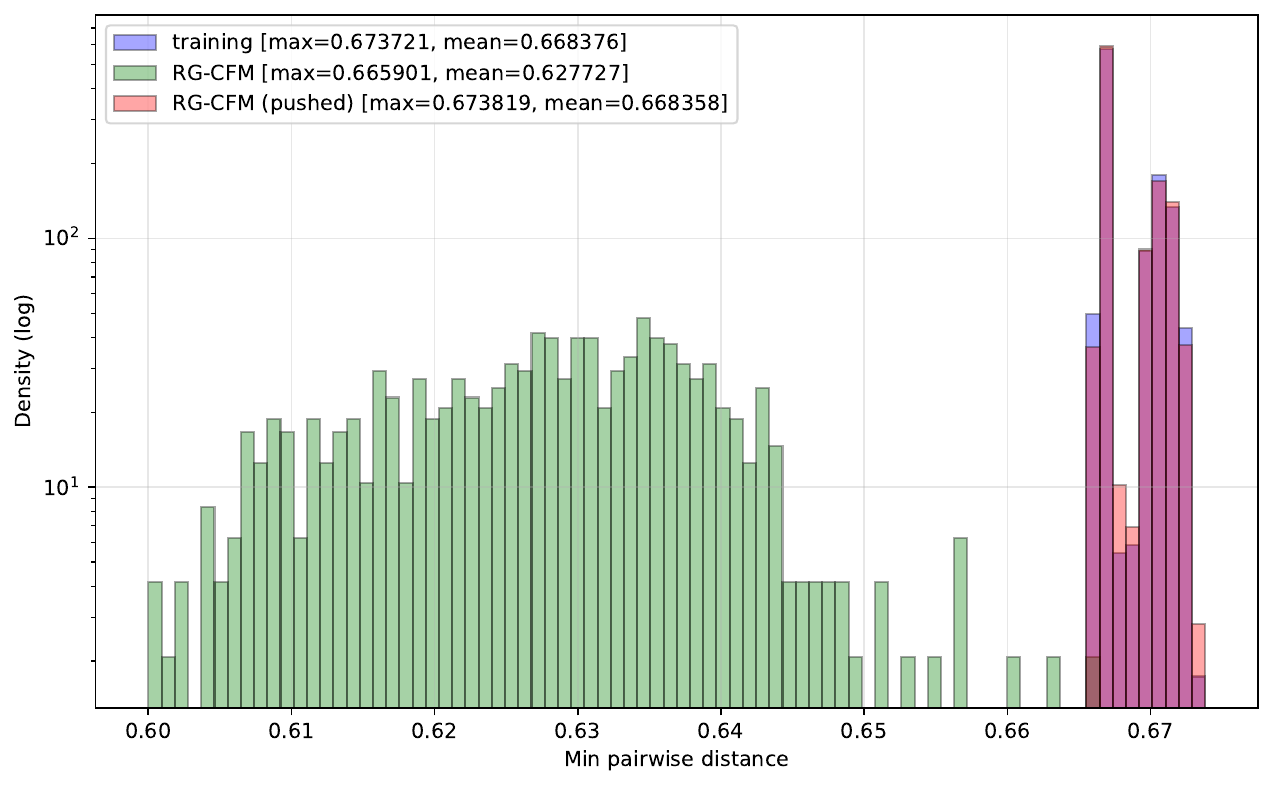}
        \caption{Training vs.\ generated vs.\ pushed.}
        \label{fig:12d-31-compare}
    \end{subfigure}
    \hfill
    \begin{subfigure}[t]{0.48\textwidth}
        \centering
        \includegraphics[width=\textwidth]{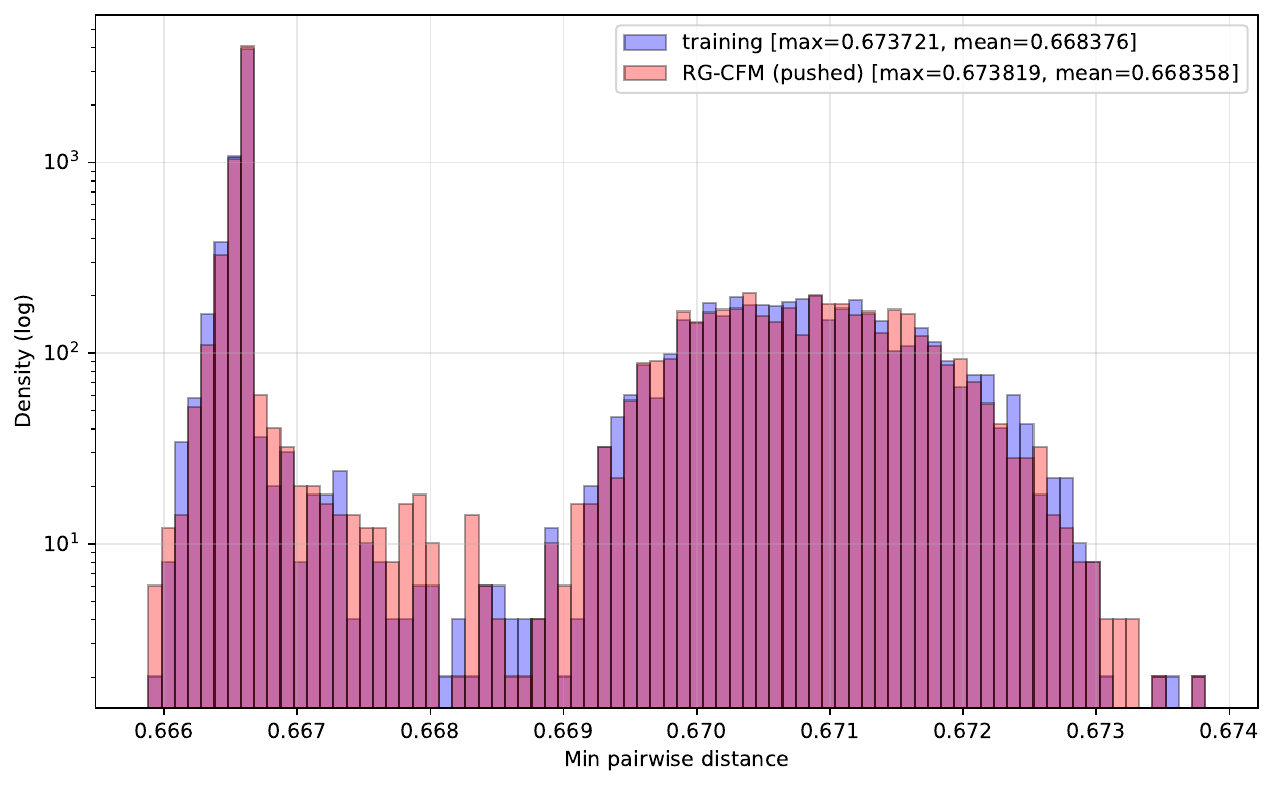}
        \caption{Training vs.\ pushed (zoomed).}
        \label{fig:12d-31-zoomed}
    \end{subfigure}
    \caption{\textbf{Sphere packing in $d=12$, $N=31$.} Normalized histograms of minimum pairwise distance (log scale).}
    \label{fig:12d-31}
\end{figure}

The pattern mirrors lower-dimensional results: raw flow samples require geometric repair, but the pushed distribution matches the training regime while the closed-loop update enables discovery of configurations exceeding the training maximum. For $N=31$, starting from SRP-generated packings with best minimum separation $d_{\min} = 0.673721$, a single round of RG-CFM fine-tuning yields $d_{\min} = 0.673819$, an improvement in a regime where any progress is nontrivial.

\subsubsection*{Large-training-set runs (long per-iteration training).}
\begin{figure}[t]
\centering
\begin{minipage}[t]{0.49\textwidth}
  \centering
  \includegraphics[width=\linewidth]{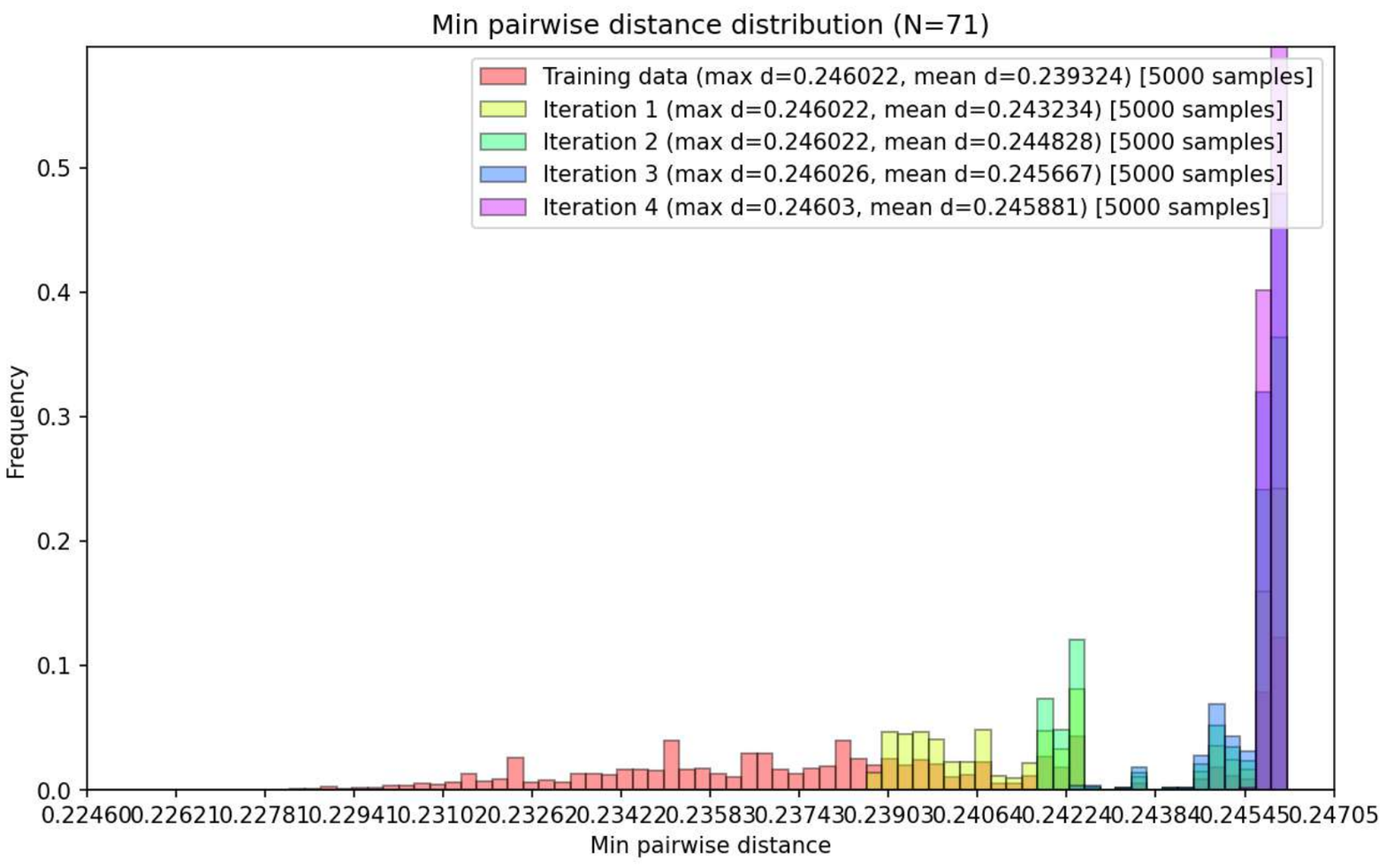}
  \caption*{$N=71$}
\end{minipage}\hfill
\begin{minipage}[t]{0.49\textwidth}
  \centering
  \includegraphics[width=\linewidth]{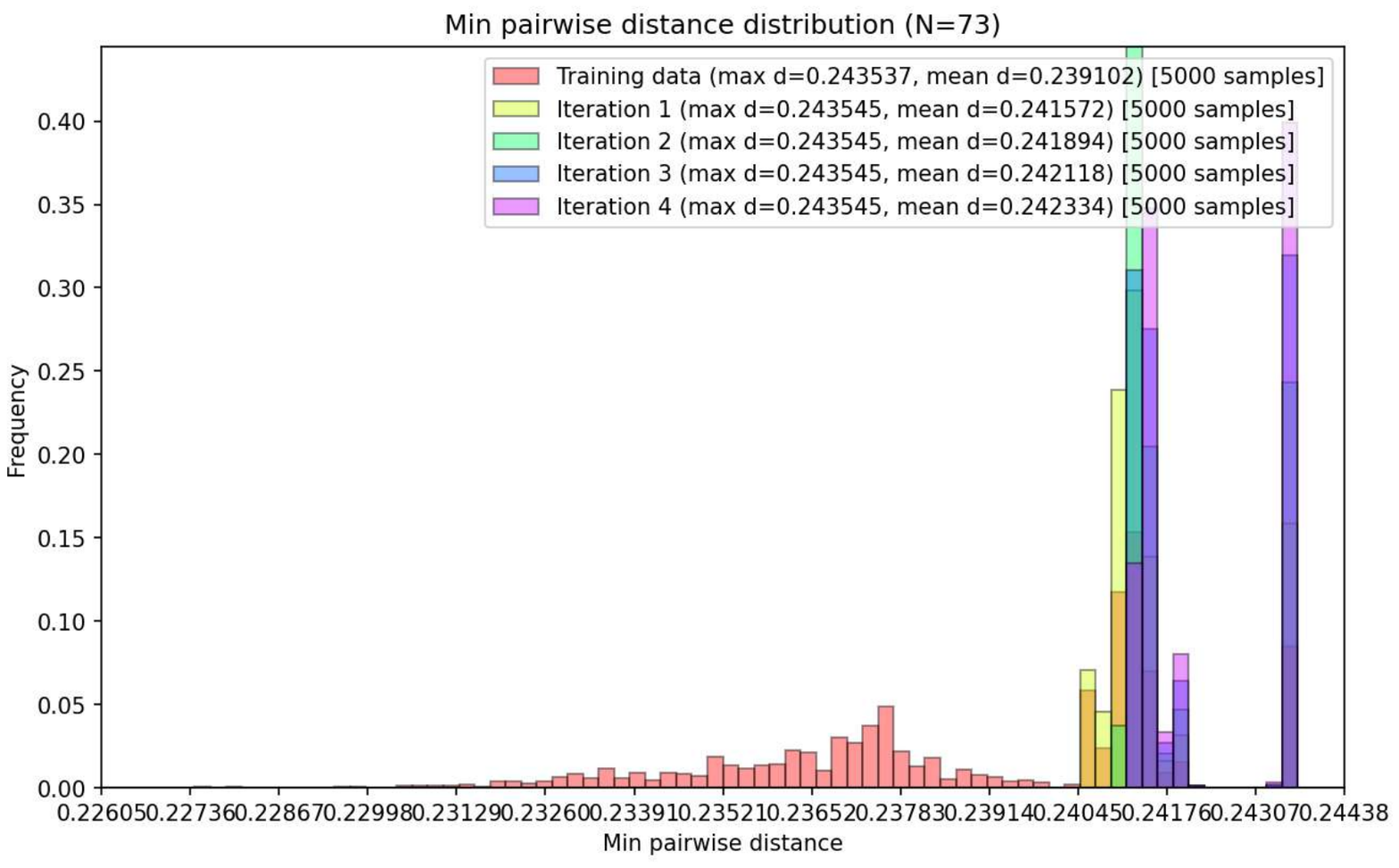}
  \caption*{$N=73$}
\end{minipage}
\begin{minipage}[t]{0.49\textwidth}
  \centering
  \includegraphics[width=\linewidth]{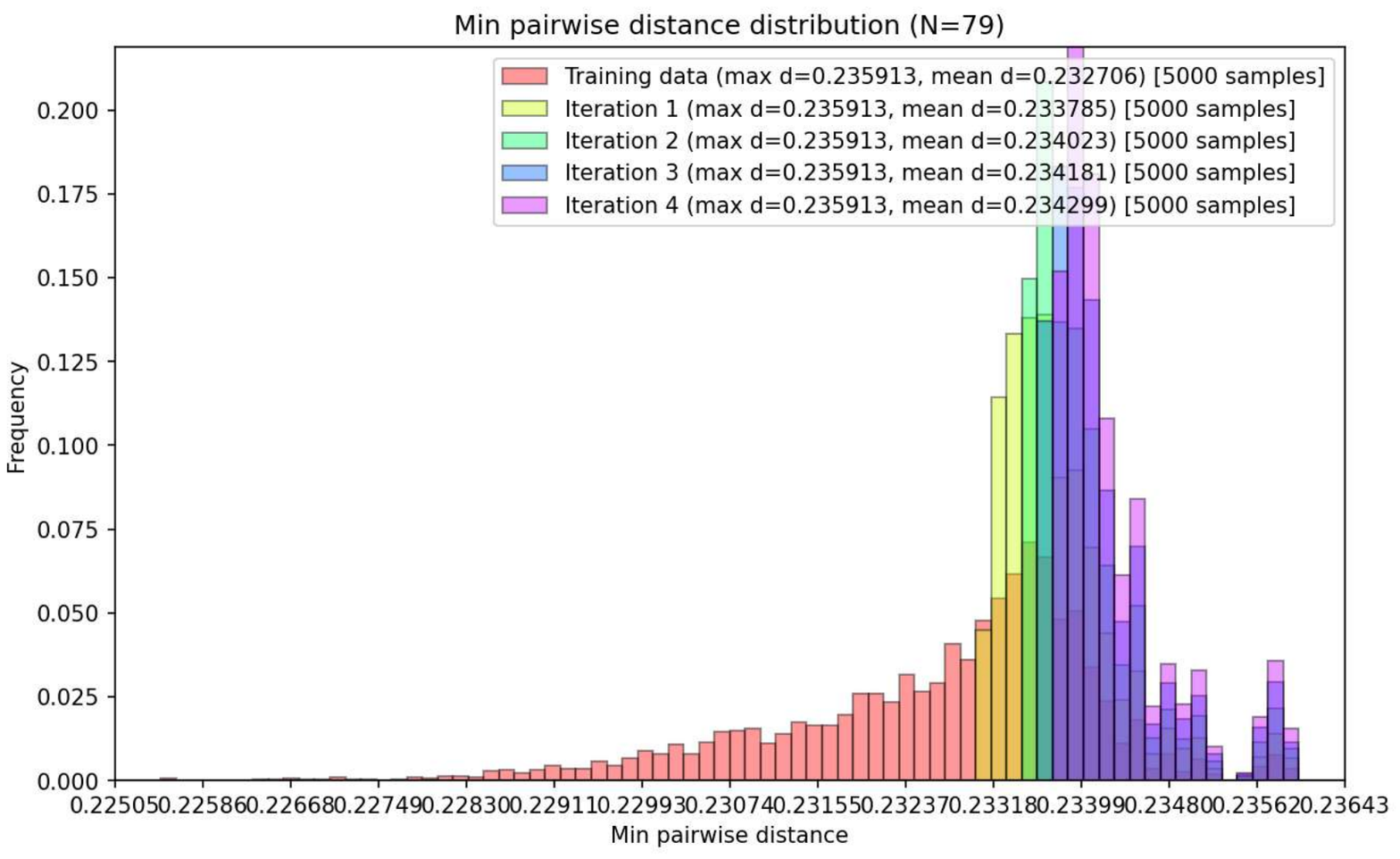}
  \caption*{$N=79$}
\end{minipage}\hfill
\begin{minipage}[t]{0.49\textwidth}
  \centering
  \includegraphics[width=\linewidth]{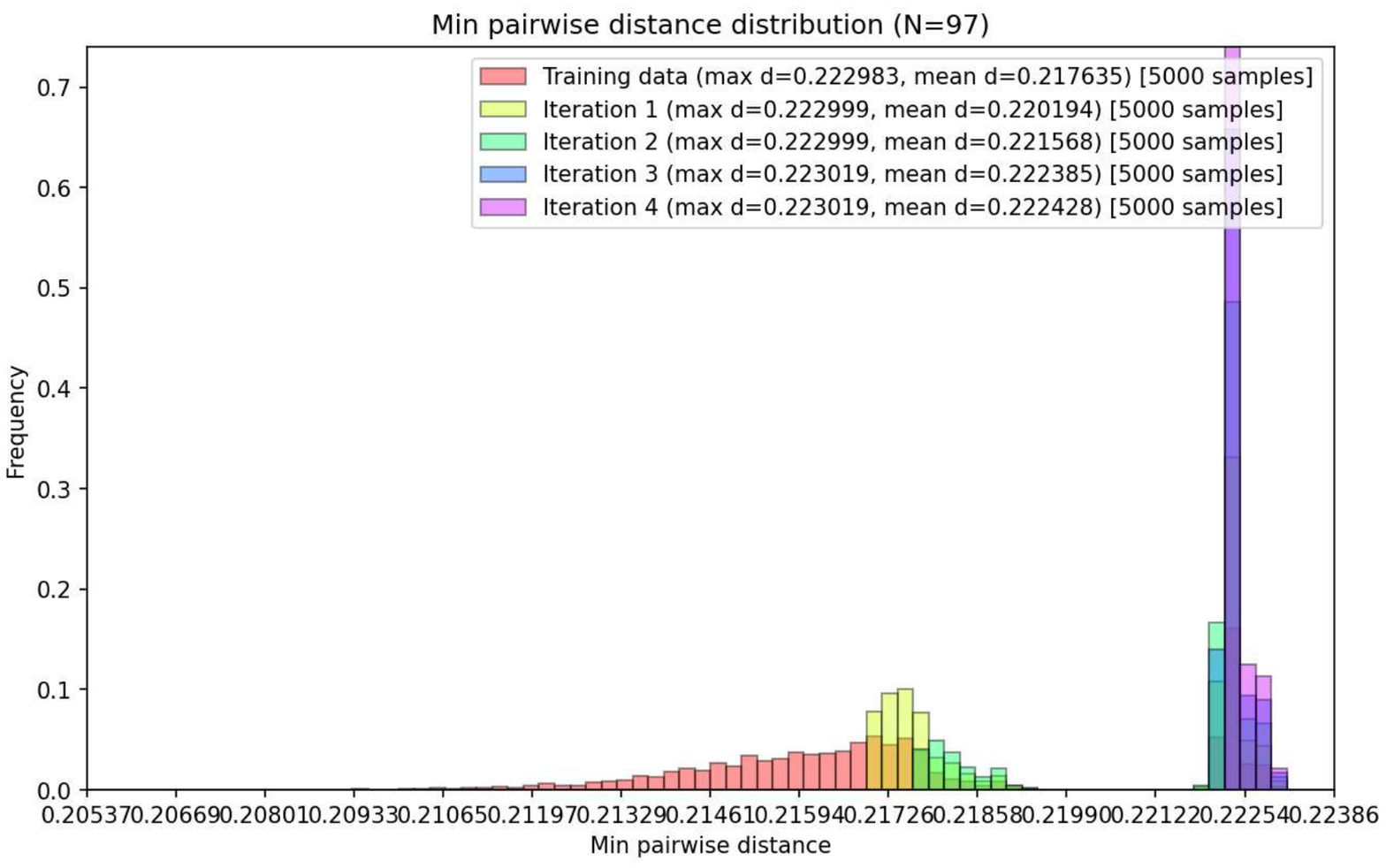}
  \caption*{$N=97$}
\end{minipage}
\begin{minipage}[t]{0.49\textwidth}
  \centering
  \includegraphics[width=\linewidth]{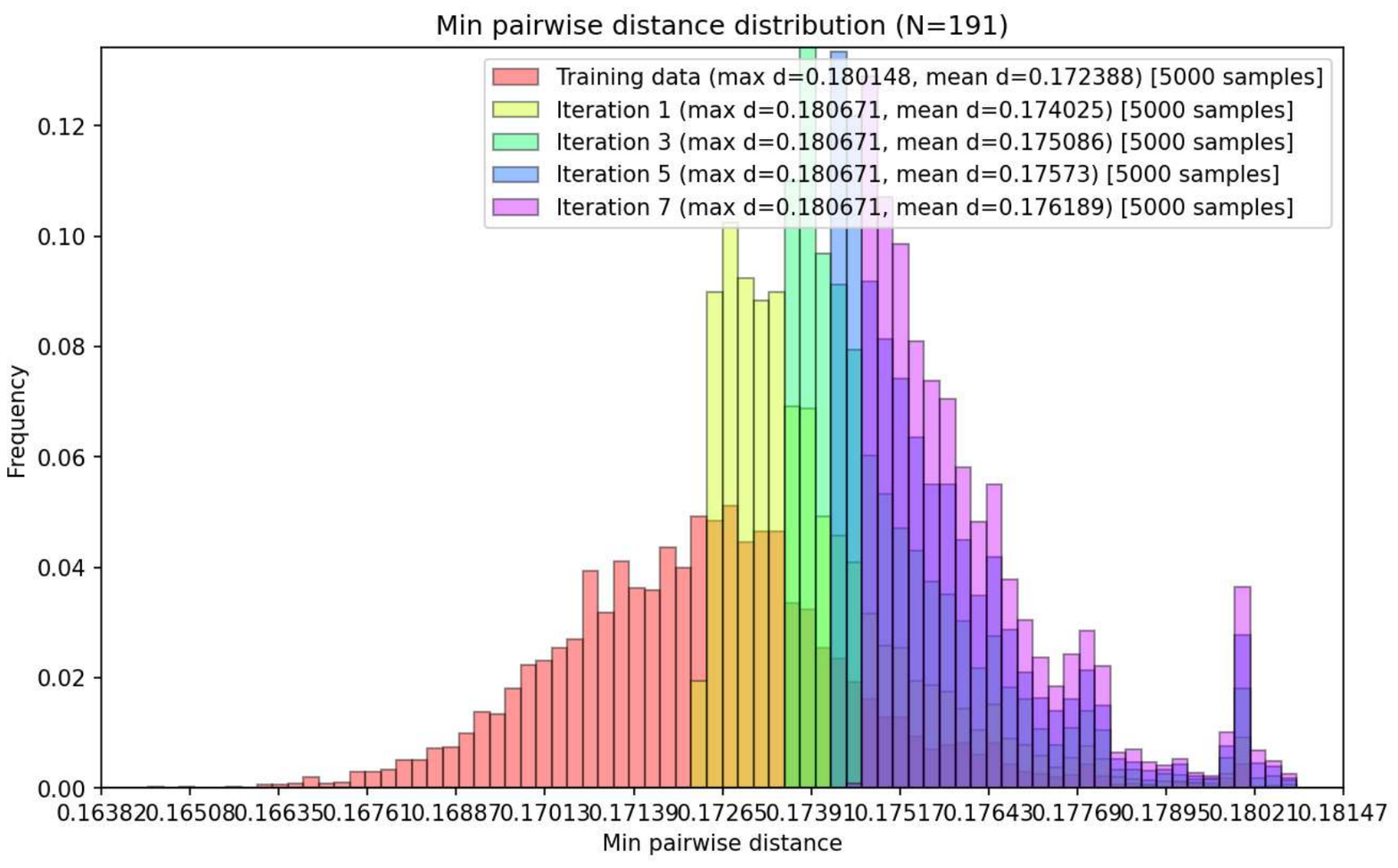}
  \caption*{$N=191$}
\end{minipage}
\caption{\textbf{3D unit-cube packing proxy: large training set; long per-iteration training.}
The distribution of $d_{\min}$ shifts right over iterations, improving both the mean and the best observed value.}
\label{fig:spheres-hists-71-191}
\end{figure}

Figure \ref{fig:spheres-hists-71-191} summarizes runs with a large SRP-generated training set and long per-iteration flow matching model training ($3000$ epochs) for $N\in\{71,73,79,97,191\}$.
In each case the main effect of boosting is a strong \emph{right-shift} and concentration of the distribution:
high-$d_{\min}$ configurations become substantially more frequent, and the mean improves monotonically, while the absolute best value typically improves only slightly (or stays essentially fixed). Concretely, 
\begin{itemize}[nosep,leftmargin=*]
\item $N=71$: mean improves from $0.239324$ (training) to $0.245881$ (iteration~4), while the best improves marginally from $0.246022$ to $0.246030$.
\item $N=73$: mean improves from $0.239102$ (training) to $0.242334$ (iteration~4), with best improves from $0.243537$ to $0.243545$.
\item $N=79$: mean improves from $0.232706$ (training) to $0.234299$ (iteration~4), with best unchanged at $0.235913$.
\item $N=97$: mean improves from $0.217635$ (training) to $0.222428$ (iteration~4), and the best improves from $0.222983$ to $0.223019$.
\item $N=191$: mean improves from $0.172388$ (training) to $0.176189$ (iteration~7), and the best improves from $0.180148$ to $0.180671$.
\end{itemize}
Thus, even if the best value is hard to move, \FMBoost{} substantially increases the probability of sampling near-record configurations after the final push.
For context, best-known values for spheres-in-a-cube are often benchmarked against the Packomania database.\footnote{\url{https://www.packomania.com/}.}

\subsubsection*{Many-iteration test for $N=191$: short training per step.} 
\begin{figure}[t]
\centering
\begin{minipage}[t]{0.64\textwidth}
  \centering
  \includegraphics[width=\linewidth]{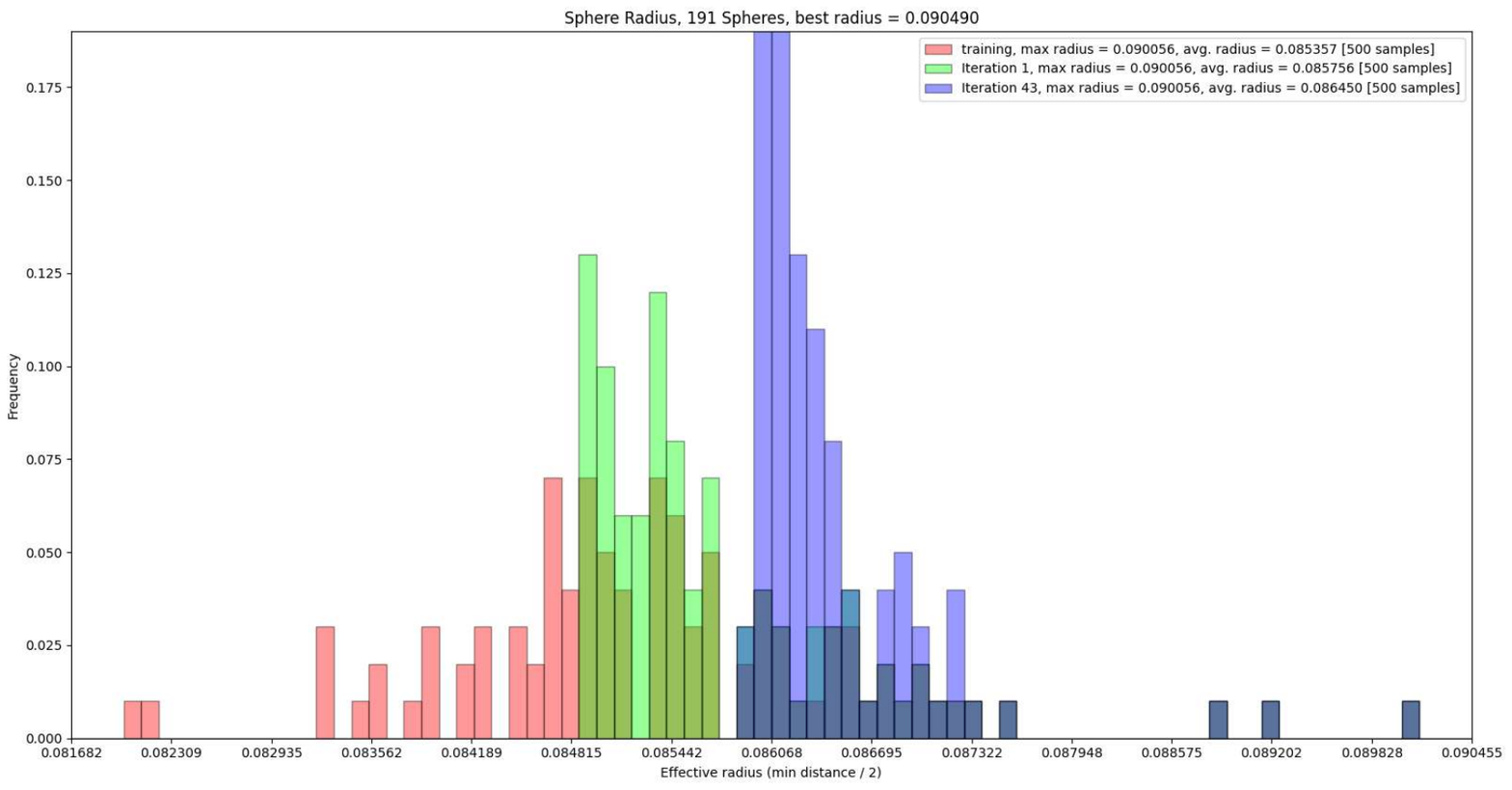}
  \caption*{\textbf{(a)} Distribution snapshots ($N=191$).}
\end{minipage}\hfill
\begin{minipage}[t]{0.34\textwidth}
  \centering
  \includegraphics[width=\linewidth]{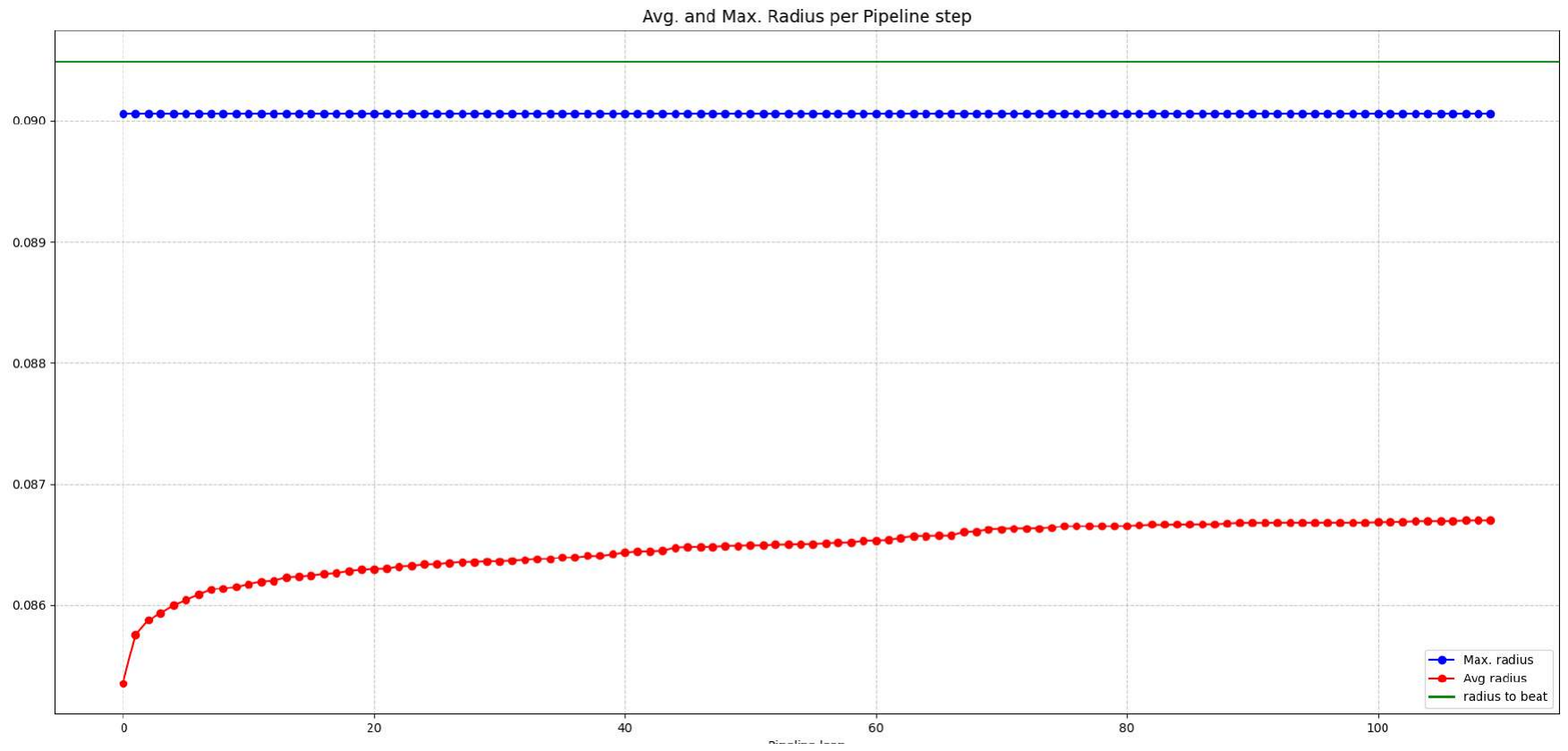}
  \caption*{\textbf{(b)} Mean/max $r_{\mathrm{eff}}$ vs.\ iteration.}
\end{minipage}

\caption{\textbf{Many-iteration regime for $N=191$: small per-iteration budget.}
 (a) shows histogram snapshots of $r_{\mathrm{eff}}=\tfrac12 d_{\min}$ at selected iterations.
 (b) tracks the average and maximum $r_{\mathrm{eff}}$ across $>100$ pipeline steps:
the maximum remains stable, while the average improves monotonically, indicating a steady shift of mass toward better basins.}
\label{fig:spheres-191-manyiters}
\end{figure}

Figure \ref{fig:spheres-191-manyiters} shows a complementary experiment at $N=191$ with a smaller per-iteration budget
($500$ samples; $300$ epochs per iteration) but more than $100$ pipeline iterations. 
Here the maximum effective radius remains essentially flat, while the average effective radius increases steadily across iterations.
This is the behavior one expects from a stable closed-loop booster under limited compute: the pipeline learns to place more mass in good basins (raising the mean) even when discovering a strictly better extreme configuration is rare.

\subsection{The Heilbronn Problem}
\label{sec:heilbronn}

We test \FMBoost{} on the classical Heilbronn triangle problem in the unit square. 

\subsubsection{Local search}

As described in Section 2, we use an SRP local search algorithm both for generating training data and for the final push applied to model samples.  For Heilbronn we use a differentiable surrogate that interpolates between a smooth global objective and the hard minimum.
We replace the non-differentiable absolute value by the smooth proxy
$|u|\approx \sqrt{u^2+\varepsilon}$ ($\varepsilon=10^{-12}$),
and define the smoothed triangle area
\[
A^{(\varepsilon)}_{ijk}(X)
:=
\frac12\sqrt{\det(p_j-p_i,p_k-p_i)^2+\varepsilon},
\qquad 1\le i<j<k\le n .
\]
To approximate the hard minimum $A_{\min}(X)=\min_{i<j<k}A_{ijk}(X)$ in a differentiable way, we use the
\emph{soft-min} (log-sum-exp) at sharpness $\beta>0$:
\begin{equation}
\label{eq:heil-softmin}
\operatorname{smin}_\beta(X)
:=
-\frac{1}{\beta}\log\!\sum_{i<j<k}\exp \bigl(-\beta\,A^{(\varepsilon)}_{ijk}(X)\bigr).
\end{equation}
Let $T=\binom{n}{3}$ and $m(X):=\min_{i<j<k}A^{(\varepsilon)}_{ijk}(X)$.  Then the standard log-sum-exp bounds give
\begin{equation}
\label{eq:heil-softmin-gap}
m(X)-\frac{\log T}{\beta}\le \operatorname{smin}_\beta(X) \le m(X),
\end{equation}
so increasing $\beta$ tightens the approximation (up to an explicit $\log T/\beta$ gap).  In practice we anneal $\beta$
geometrically from $\beta_0$ to $\beta_F$, which has the effect of starting with a smoother objective (many triangles contribute)
and gradually concentrating optimization pressure on the worst triangle(s).

We optimize the SRP surrogate loss
\begin{equation}
\label{eq:heil-loss}
\calL_\beta(X)
= w_{\mathrm{wall}} W(X)-\operatorname{smin}_\beta(X),
\end{equation}
where $W(X)$ is a quadratic wall penalty that discourages leaving $[0,1]^2$ (and we additionally clamp coordinates after each SRP step) and $w_{\mathrm{wall}}>0$ is the penalty weight that controls how strongly SRP discourages leaving the unit square:
larger $w_{\mathrm{wall}}$ makes boundary violations more expensive (a standard quadratic-penalty mechanism for constraints).
In the reported runs we use geometric annealing $\beta_0\to\beta_F$ with $\beta_0=40$ and $\beta_F=300$.

For large $\beta$ the sum in \eqref{eq:heil-softmin} is dominated by the smallest triangle areas.  To reduce computation we optionally
restrict the log-sum-exp to the $K$ smallest triangles under the current iterate $X$ (with a small tolerance), i.e.\ we replace the full index set
by an active subset of size $K$ (here $K=100$).  This keeps the dominant terms while avoiding work on clearly non-critical triangles.

After SRP we apply a deterministic L-BFGS-B polish to (locally) minimize $\calL_{\beta_F}$ under box constraints.
\texttt{L-BFGS-B} is a quasi-Newton method implemented in \texttt{scipy.optimize.minimize(method="L-BFGS-B")}, that stores only a low-rank approximation of curvature information, making it well suited to high-dimensional problems with simple bounds.
To target the true max-min objective more directly, we then solve a lifted nonlinear program in variables $(X,t)$:
\[
\max\ t\quad\text{s.t.}\quad A_{ijk}(X)\ge t\ \ \forall\, i<j<k,\qquad X\in[0,1]^{2n}.
\]
Since only few constraints are typically tight at a good solution, we enforce this system only for an \emph{active set}
consisting of the $K_{\mathrm{active}}$ currently smallest-area triangles (here $K_{\mathrm{active}}=25$), and run an SLSQP step on these constraints.\footnote{SLSQP stands for \emph{Sequential Least Squares Programming}, as implemented in \texttt{scipy.optimize.minimize(method="SLSQP")}.  It is a standard nonlinear constrained optimizer supporting both bound constraints and general equality/inequality constraints.}%
This active max-min pass is applied only to a small elite subset of candidates (top $10$ configurations) during training-set generation,
and again as the final push on flow-generated samples.

\subsubsection{Training and sampling}

We use the conditional flow matching setup from Section~2, with a permutation-equivariant set transformer velocity field.  The conditioning encodes the problem size and the target quality $c(X) = \Bigl(\frac{n}{128},A_{\min}(X)\Bigr)$, and an auxiliary penalty encourages the projected endpoint (as in Section~2) to achieve soft-min triangle area at least as large as the target.  

At sampling we follow the same geometry-aware sampling principle as GAS (Section~2), specialized to the Heilbronn objective: between ODE steps we perform short projected gradient-ascent moves on the soft-min triangle area (with an annealed temperature $\tau_{\mathrm{start}}=5\cdot10^{-3}\to\tau_{\mathrm{end}}=5\cdot10^{-4}$), and we end with a brief polishing stage.  Finally, every generated configuration is passed through the SRP final push above, and only then compared using exact $A_{\min}$.

\subsubsection{Experimental settings}

For $n=15$ we generate $1000$ SRP samples for training, keeping the top $50\%$ by $A_{\min}$ for model fitting.  SRP uses $I_{\max}=500$ outer iterations, $m=60$ inner steps, step size $0.035$ with decay $0.994$, and backtracking budget $6$; L-BFGS-B uses gtol $10^{-8}$, ftol $10^{-12}$, maxiter $2000$.  The flow model uses a set transformer of width $256$, depth $6$, and $8$ heads, trained for $500$ epochs with batch size $64$ and learning rate $10^{-4}$.  Sampling uses $30$ PCFM steps (Euler, step cap $0.05$), with $2$ projection and $2$ proximal steps per PCFM step, plus a $20$-step terminal polish.

\subsubsection{Results}

\begin{figure}[t]
\centering
\begin{minipage}[t]{0.49\textwidth}
  \centering
  \includegraphics[width=\linewidth]{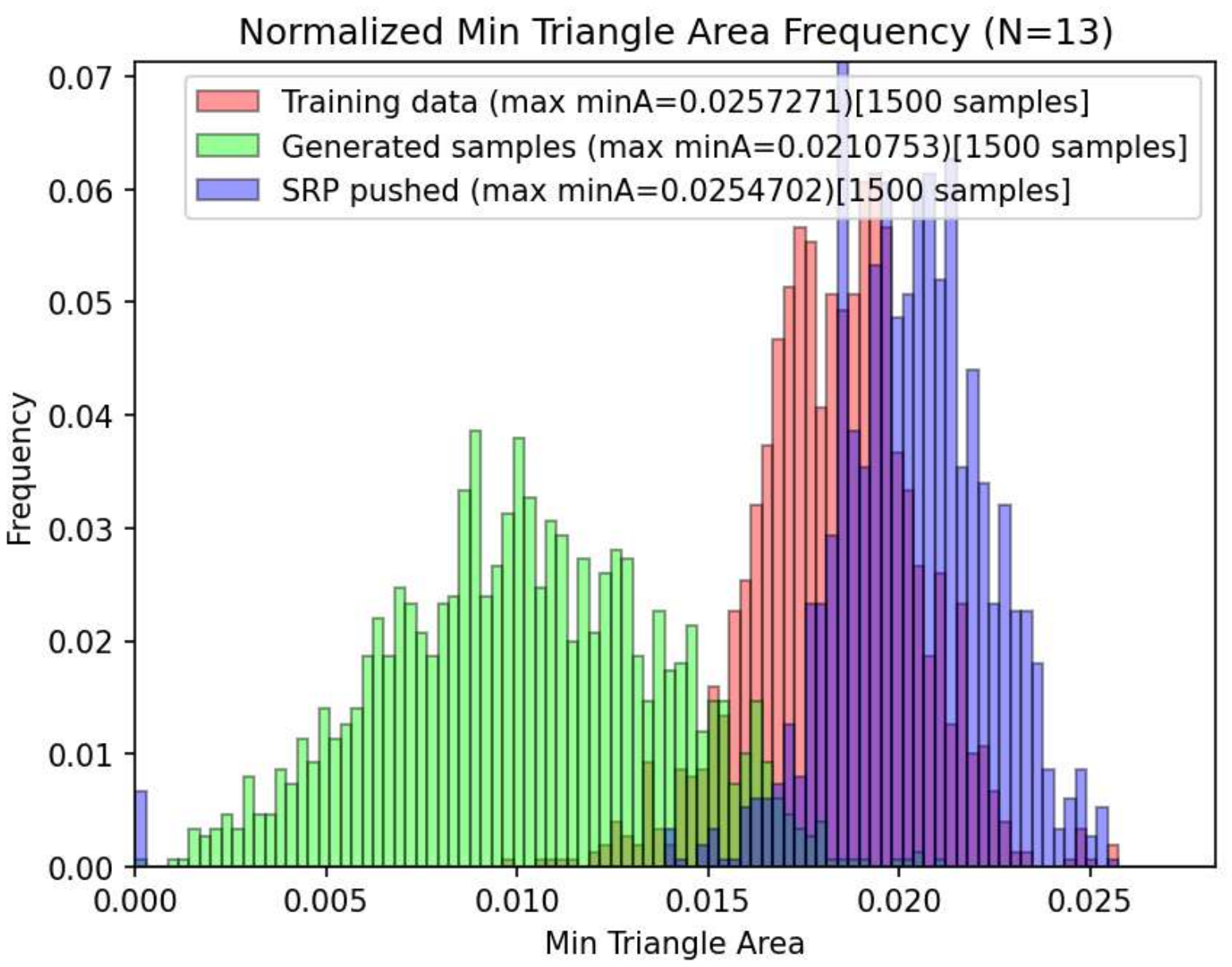}
  \caption*{\textbf{(a)} Training vs.\ generated vs.\ pushed.}
\end{minipage}\hfill
\begin{minipage}[t]{0.49\textwidth}
  \centering
  \includegraphics[width=\linewidth]{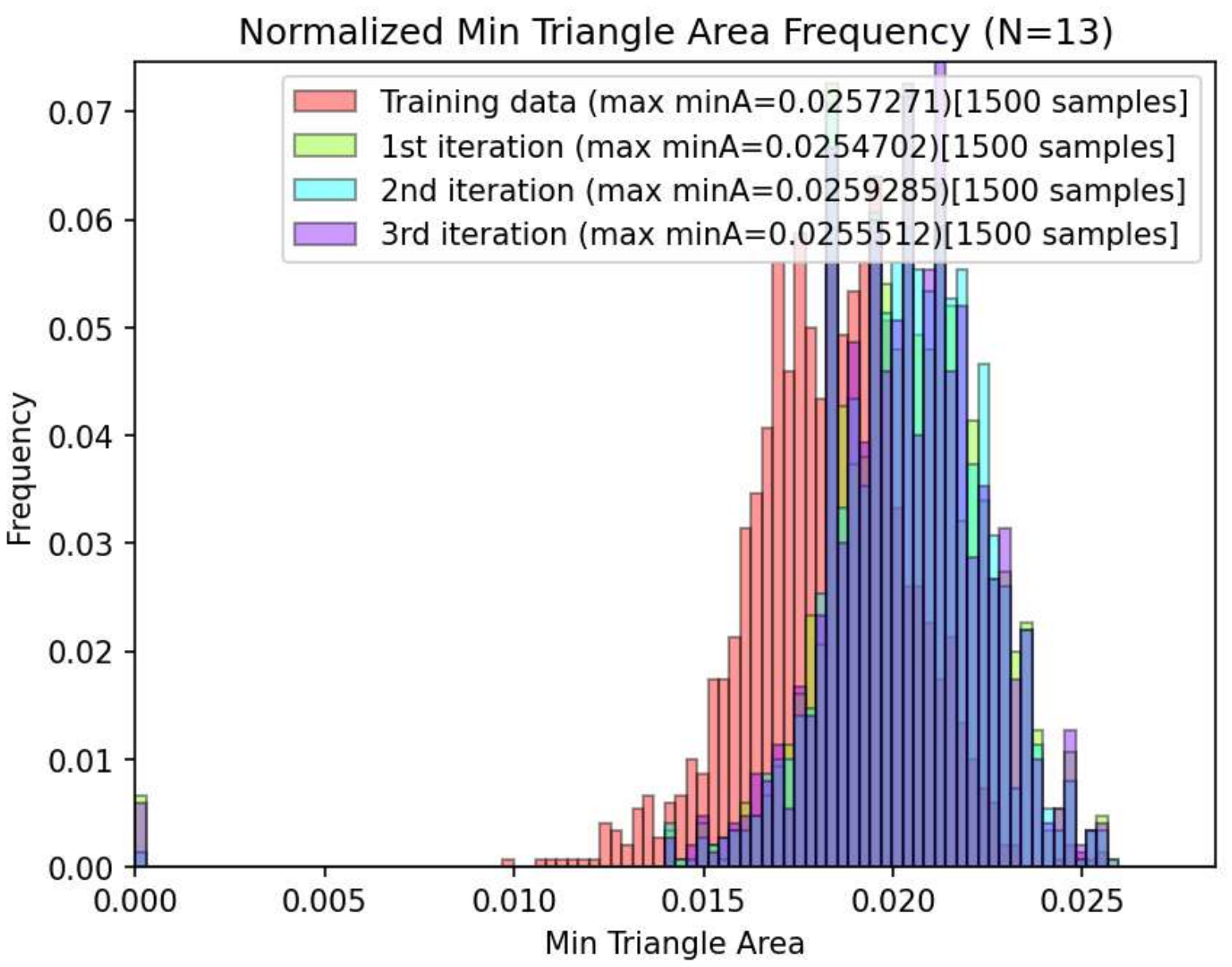}
  \caption*{\textbf{(b)} Iterations $1$--$3$.}
\end{minipage}
\caption{\textbf{Heilbronn, $n=13$.}  Normalized histograms of exact $A_{\min}$.  Legends report the best achieved $A_{\min}$ in each set.}
\label{fig:heil-n13}
\end{figure}

\begin{figure}[t]
\centering
\begin{minipage}[t]{0.49\textwidth}
  \centering
  \includegraphics[width=\linewidth]{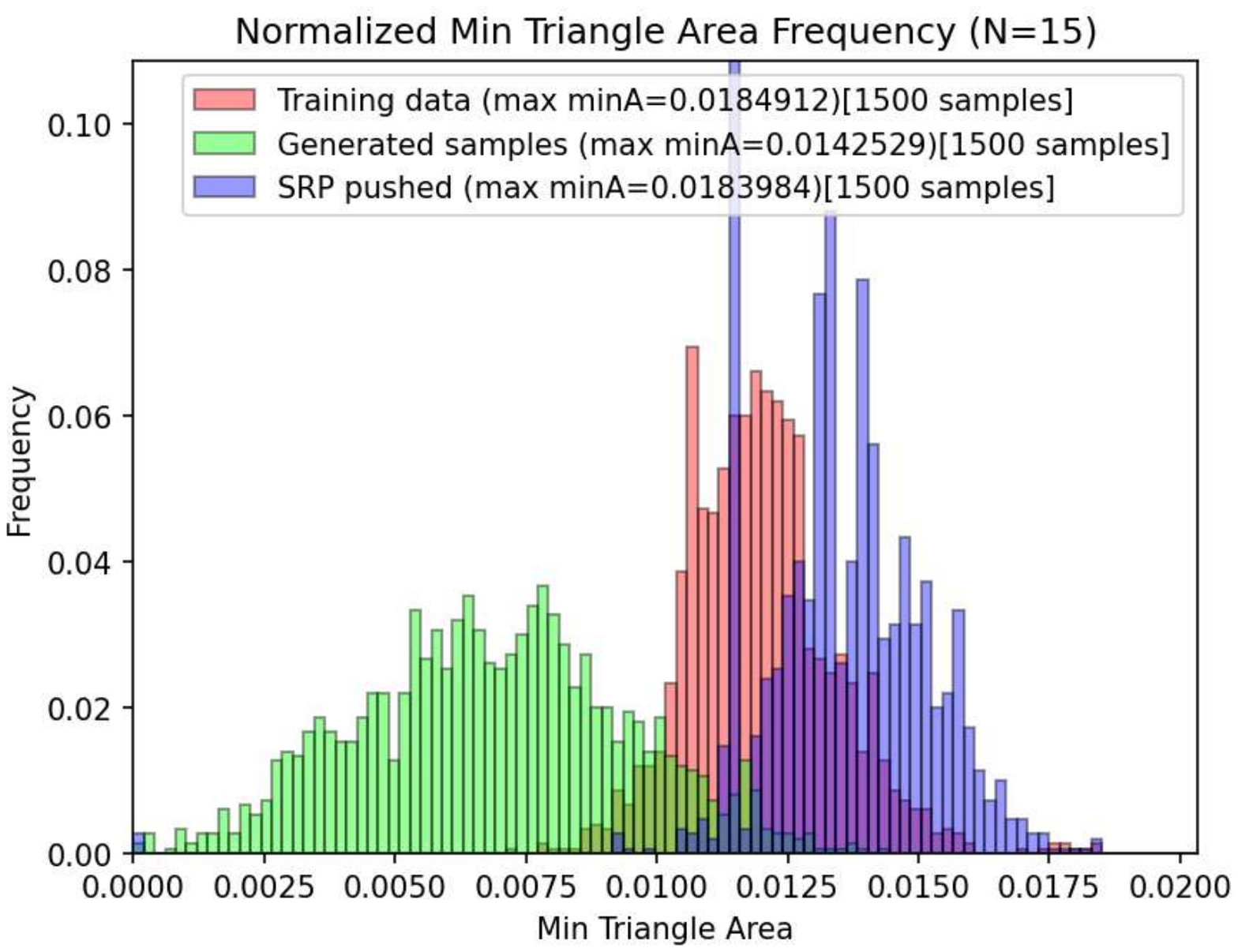}
  \caption*{\textbf{(a)} Training vs.\ generated vs.\ pushed.}
\end{minipage}\hfill
\begin{minipage}[t]{0.49\textwidth}
  \centering
  \includegraphics[width=\linewidth]{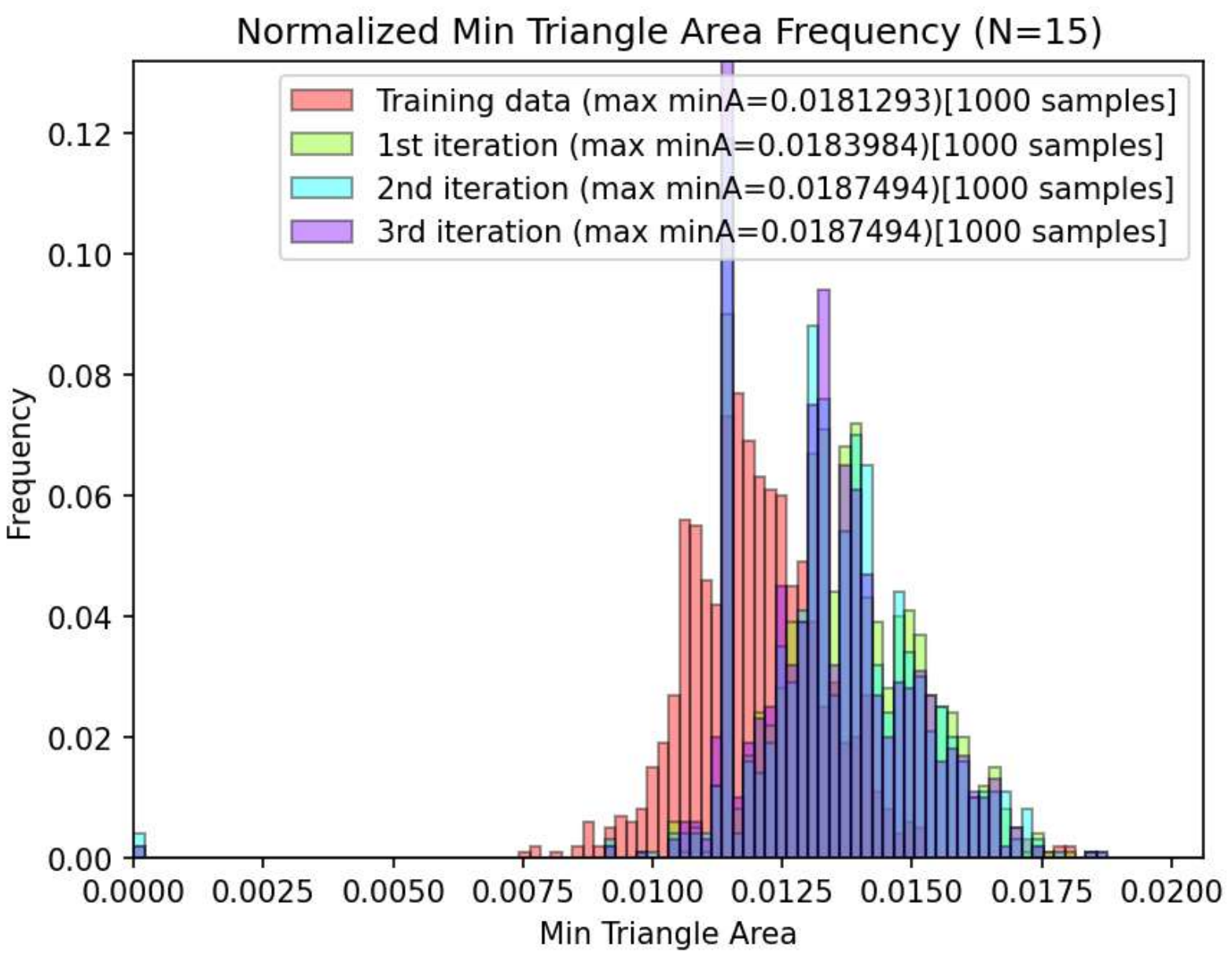}
  \caption*{\textbf{(b)} Iterations $1$--$3$.}
\end{minipage}
\caption{\textbf{Heilbronn, $n=15$.}  Normalized histograms of exact $A_{\min}$.  The final push converts "almost-good" samples into strong local optima, enabling iterative improvement.}
\label{fig:heil-n15}
\end{figure}

Figures~\ref{fig:heil-n13}--\ref{fig:heil-n15} show the empirical distributions of $A_{\min}$ for $n=13$ and $n=15$.  The main pattern is consistent across both sizes:

\begin{itemize}[nosep,leftmargin=*]
\item Raw flow samples are worse than SRP elites. Without final push (local refinement), generated samples are left-shifted (many near-degenerate triangles remain).
\item The final push repairs samples effectively and moves the distribution back toward the training regime and recovers near-elite maxima.
\item Boosting can exceed the training maximum. Across iterations, the best observed $A_{\min}$ improves and the histogram mass shifts right, indicating that the closed-loop train$\to$sample$\to$push$\to$select cycle can produce configurations that beat the previous elite set.
\end{itemize}

Concretely, for $n=13$ (Figure~\ref{fig:heil-n13}): the training max is $0.0257271$, raw generation reaches $0.0210753$, while pushing recovers $0.0254702$; over iterations the maximum improves to $0.0259285$ (iteration 2), exceeding the training maximum in that run.
For $n=15$ (Figure~\ref{fig:heil-n15}): the training max is $0.0184912$ (resp.\ $0.0181293$ in the iterative plot), raw generation is $0.0142529$, pushing recovers $0.0183984$, and the iterative run improves to $0.0187494$ by iteration 2 (matched again in iteration 3).

\begin{figure}[t]
\centering
\begin{minipage}[t]{0.5\textwidth}
  \centering
  \includegraphics[width=\linewidth]{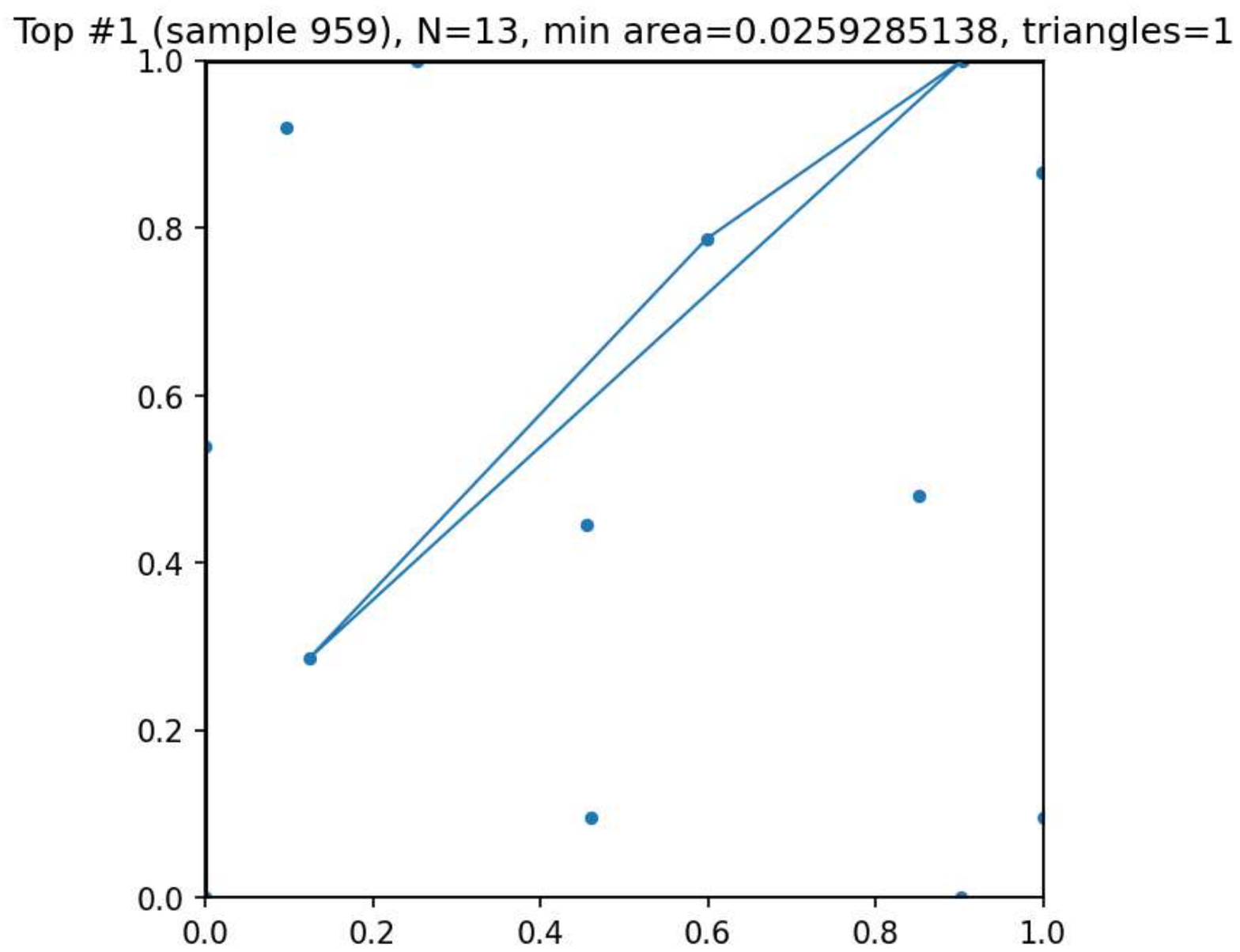}
  \caption*{$n=13$, $A_{\min}=0.0259285138$}
\end{minipage}\hfill
\begin{minipage}[t]{0.5\textwidth}
  \centering
  \includegraphics[width=\linewidth]{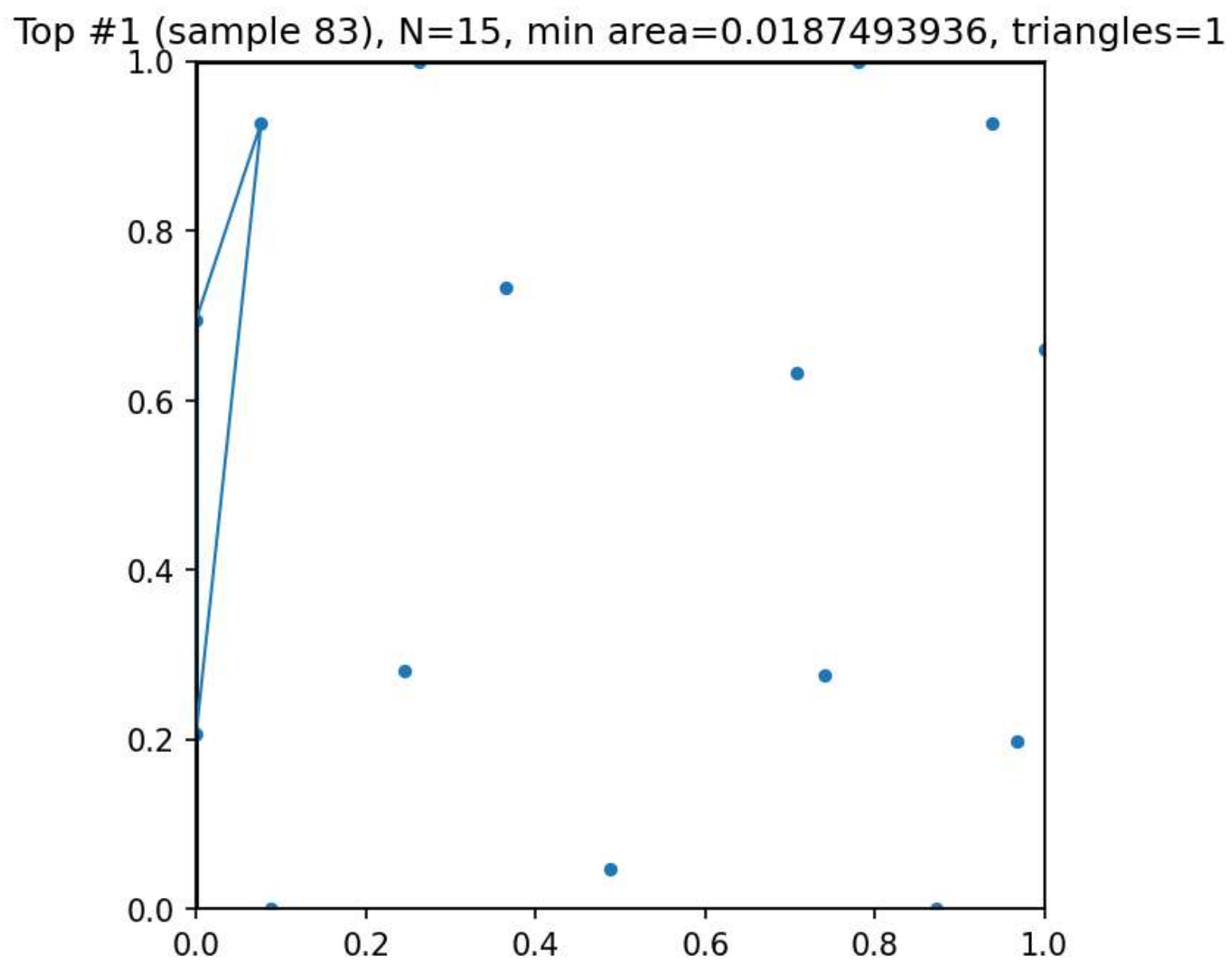}
  \caption*{$n=15$, $A_{\min}=0.0187493936$}
\end{minipage}\hfill
\caption{\textbf{Best Heilbronn configurations found by \FMBoost{} (after final push).}
The highlighted triangle(s) attain the minimum area. The best known constructions are $A_{\min} = 0.0270$ for $n=13$ and $A_{\min} = 0.0211$ for $n=15$ in \cite{FriedmanHeilbronnSquaresBenchmark}.}
\label{fig:heilbronn-best-constructions}
\end{figure}

\subsection{Circles in Unit Square with Maximal Sum of Radii}
\label{sec:sumradii}

Given $n$ circles in the unit square, the objective is to maximize the total sum of radii. Improving the sum typically requires coordinated global rearrangements of the contact graph, while feasibility is defined by a large family of hard inequalities.
A configuration consists of centers and radii
\[
X=\bigl((p_1,r_1),\dots,(p_n,r_n)\bigr),\qquad p_i=(x_i,y_i)\in[0,1]^2,\quad r_i\ge 0.
\]
We require \emph{containment} and \emph{non-overlap}:
\begin{equation}
\label{eq:sumr-feas}
r_i \le x_i, r_i \le 1-x_i, r_i \le y_i, r_i \le 1-y_i
\quad\text{and}\quad
\|p_i-p_j\| \ge r_i+r_j\ \ (i\neq j).
\end{equation}

\subsubsection{Local search}

As in Section~2, SRP is used both to generate the training set and as the final push applied to flow samples. 
In the sum-of-radii setting SRP acts on the \emph{full} variable vector
\[
X=\bigl(p_1,\dots,p_n,r_1,\dots,r_n\bigr)\in([0,1]^2)^n\times[0,\tfrac12]^n,
\qquad p_i=(x_i,y_i),
\]
and minimizes the smooth penalty objective
\begin{equation}
\label{eq:sumr-srp-loss}
\widetilde{\calL}(X)
=
w_{\mathrm{wall}}\sum_{i=1}^n \Bigl([r_i-x_i]_+^2+[r_i-(1-x_i)]_+^2+[r_i-y_i]_+^2+[r_i-(1-y_i)]_+^2\Bigr)
+
w_{\mathrm{ov}}\sum_{i<j}[r_i+r_j-\|p_i-p_j\|]_+^2
-\alpha\sum_{i=1}^n r_i,
\end{equation}
followed by a deterministic \texttt{L-BFGS-B} polish (SciPy) under the box bounds $p_i\in[0,1]^2$ and $r_i\in[0,1/2]$.
Heuristically, the first two (quadratic) terms drive constraint violations to zero (walls and overlaps), while the last term encourages
large total radius.

After SRP+\texttt{L-BFGS-B} we apply a post-processing step: we freeze the centers $p_1,\dots,p_n$ produced by the local search and then recompute radii by solving the
\emph{exact} max-sum-radii linear program:
\begin{equation}
\label{eq:sumr-lp}
\max_{r\in\R_{\ge 0}^n}\ \sum_{i=1}^n r_i
\quad\text{s.t.}\quad
r_i \le \min\{x_i,1-x_i,y_i,1-y_i\},\ \ 
r_i+r_j \le \|p_i-p_j\|\ (i<j).
\end{equation}
This LP is solved with the HiGHS backend via \texttt{scipy.optimize.linprog(method="highs")}.  The LP~\eqref{eq:sumr-lp} returns a vector $r^\star$ that is \emph{(weakly) feasible by construction} for the frozen centers
$p_1,\dots,p_n$: the circles are allowed to be tangent to each other or to the boundary. Since the objective is $\sum_i r_i$, this $r^\star$ maximizes the total radius among all radii compatible with the given centers.
In practice we then apply a tiny \emph{safety shrink} before saving/plotting:
\[
r_i \leftarrow \max\{0, r_i^\star-\delta\}\qquad(\delta \approx 10^{-9}),
\]
so that strict inequalities $r_i+r_j<\|p_i-p_j\|$ and $r_i<\min\{x_i,1-x_i,y_i,1-y_i\}$ hold numerically.
This does \emph{not} change the geometry in any meaningful way (it perturbs $\sum_i r_i$ by at most $n\delta$), but it prevents borderline tangencies from
appearing as tiny overlaps due to rounding error in downstream evaluation and visualization.

Note that the LP is \emph{not} used instead of SRP; rather, SRP (and the polish) searches for a good geometry of centers, while the LP performs an
exact hard projection of the radii to the best feasible values compatible with that geometry.  This plays the same conceptual role as a projection step in GAS:
it enforces the hard constraints and optimizes the objective in a subproblem that becomes convex/linear once the centers are fixed.

For pushed flow samples we use a slightly more robust variant: keeping the flow-generated centers fixed, we run SRP+polish+LP twice, once initialized from the flow radii and once from small random radii, and keep the better outcome by $\sum r_i$.  This mitigates occasional poor radius initialization from the generator.

\subsubsection{Training and sampling}

Our flow model generates \emph{centers only}.  Concretely, the set-transformer velocity field is trained on $(x,y)$-tokens (two channels), while radii are treated as auxiliary data used only to compute conditioning statistics and to evaluate/push samples.  This decoupling is deliberate: radii are highly constrained by the contact graph and are cheaply recovered by the LP~\eqref{eq:sumr-lp} once good centers are found.
Each training sample carries the summary statistics
\[
c(X)=\Bigl(\frac{n}{128},\ \sum_i r_i,\ \min_{i<j}\bigl(\|p_i-p_j\|-(r_i+r_j)\bigr),\ \min_i\min\{x_i-r_i,1-x_i-r_i,y_i-r_i,1-y_i-r_i\}\Bigr),
\]
so the model sees both target quality ($\sum r_i$) and slack information (pair and wall clearances).

Sampling follows the PCFM-style loop used throughout (Section~2): we integrate the learned ODE for centers and interleave with simple box projection for feasibility of centers.  Radii are then set to zero at generation time and recovered by the final push (SRP+LP) before scoring.

\subsubsection{Experimental settings}

In the runs shown below we use $3000$ samples per iteration for each $n\in\{26,30,32\}$.  A representative configuration (used for $n=26$ centers-only retrain-and-sample) is:
batch size $64$, learning rate $2\cdot 10^{-4}$, $1000$ epochs, and \texttt{train\_top\_fraction}$=0.5$; sampling uses \texttt{pcfm\_steps}$=40$ with midpoint integration and a small projection/proximal schedule.
\subsubsection{Results}

\begin{figure}[t]
\centering
\includegraphics[width=0.92\linewidth]{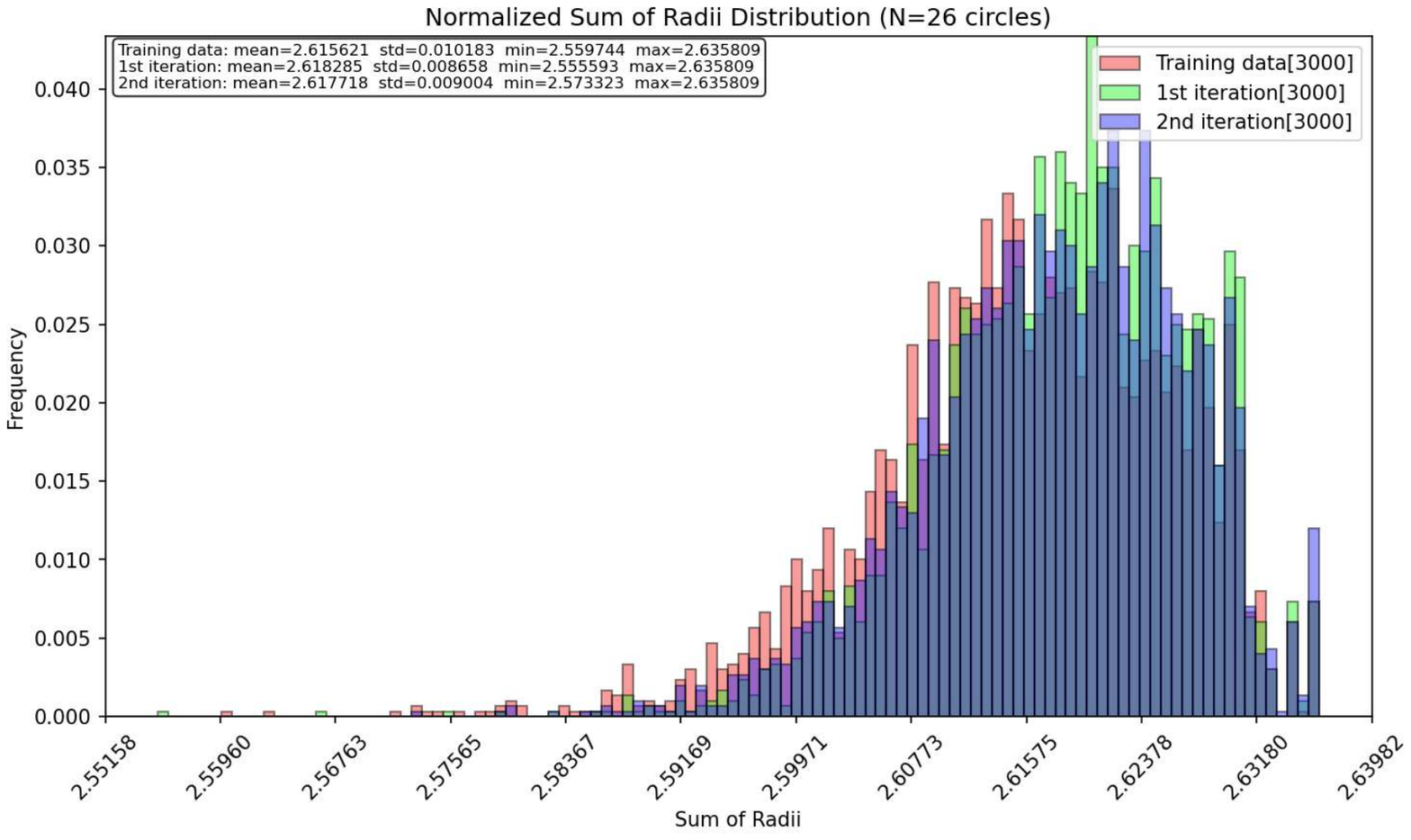}
\caption{\textbf{Sum-of-radii circle packing, $n=26$.} Normalized histograms of $\sum_i r_i$ over the training set and two \FMBoost{} iterations (each with $3000$ samples, after final push).}
\label{fig:sumr-26}
\end{figure}

\begin{figure}[t]
\centering
\includegraphics[width=0.92\linewidth]{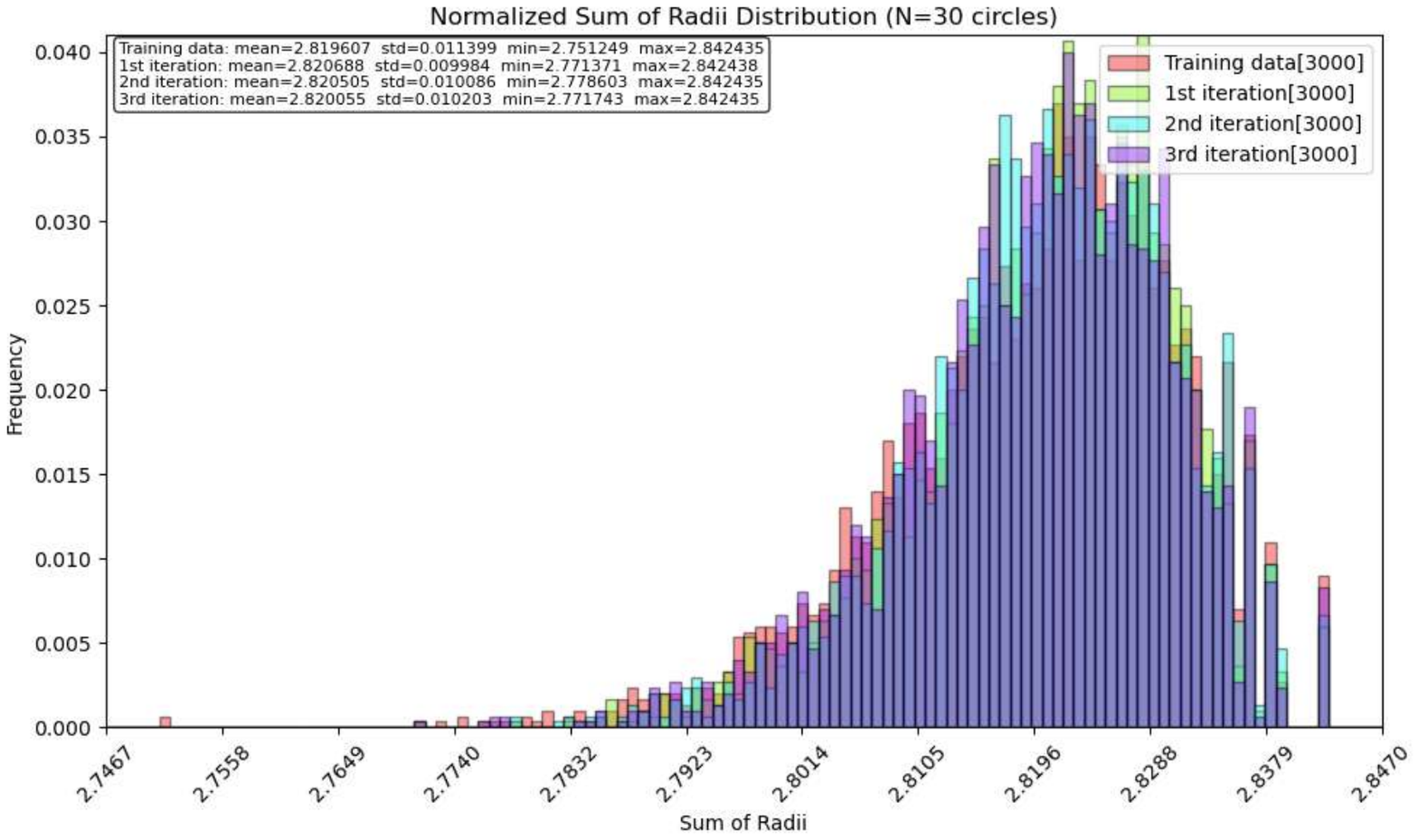}
\caption{\textbf{Sum-of-radii circle packing, $n=30$.} Normalized histograms of $\sum_i r_i$ over the training set and three \FMBoost{} iterations (each with $3000$ samples, after final push).}
\label{fig:sumr-30}
\end{figure}

\begin{figure}[t]
\centering
\includegraphics[width=0.92\linewidth]{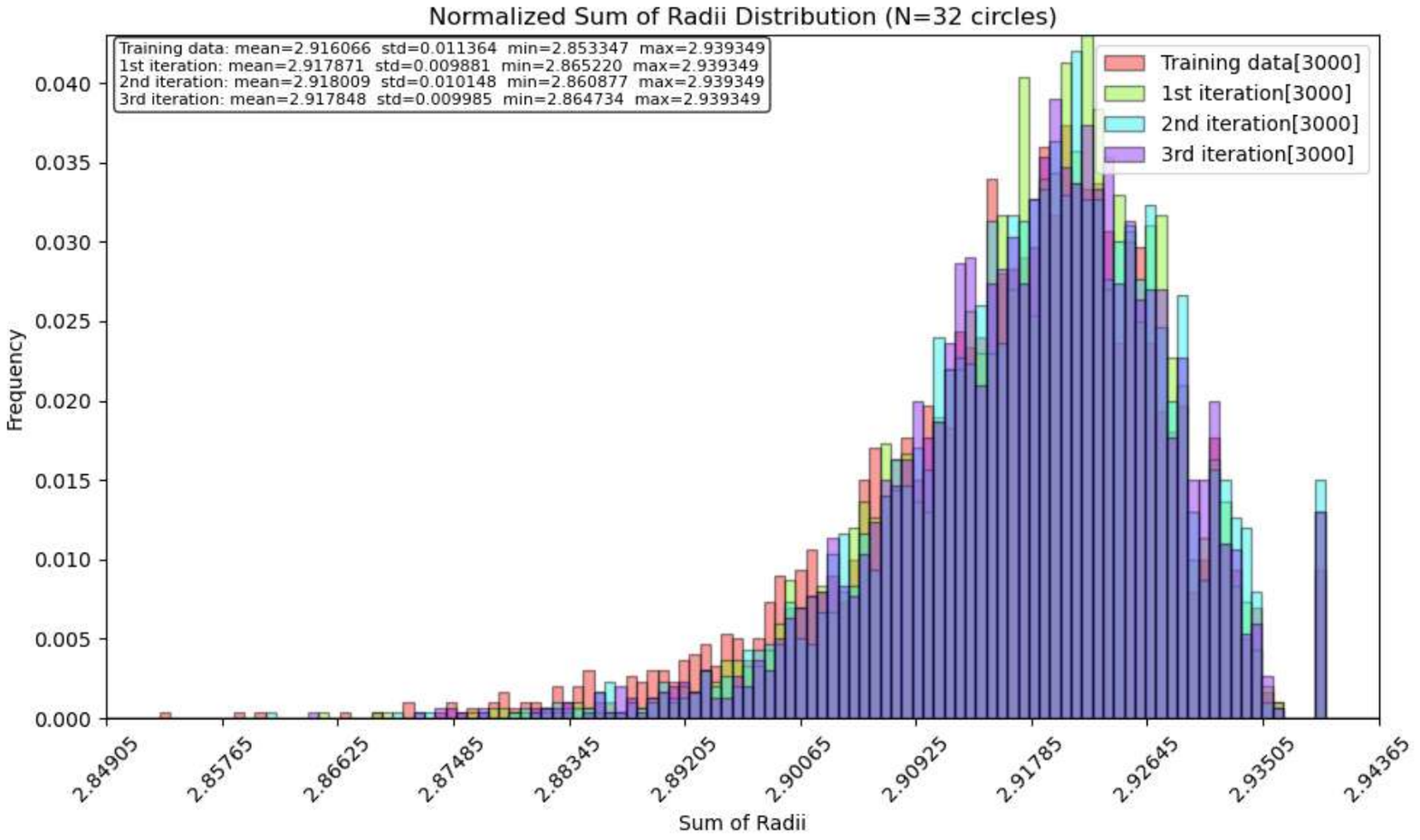}
\caption{\textbf{Sum-of-radii circle packing, $n=32$.} Normalized histograms of $\sum_i r_i$ over the training set and three \FMBoost{} iterations (each with $3000$ samples, after final push).}
\label{fig:sumr-32}
\end{figure}

Figures \ref{fig:sumr-26}-\ref{fig:sumr-32} show the empirical distributions of $\sum_i r_i$ over successive \FMBoost{} iterations.  Across all three $n$, the main effect of boosting is to shift probability mass toward the right tail: the mean increases slightly and the low-score mass is reduced, while the best observed configurations remain at the extreme right edge of the distribution.

For $n=26$ the training distribution has mean $\approx 2.6156$ and maximum $2.635809$.  After one \FMBoost{} iteration the mean increases to $\approx 2.6183$ (and remains $\approx 2.6177$ in the second), with a visible right-shift and reduced variance, while the maximum is maintained at $2.635809$.  In particular, the pipeline reliably reproduces near-record configurations and increases their frequency.

For $n=30$ the distributions are already tight around $\sum r_i\approx 2.82$.  Boosting primarily trims the left tail; the mean is stable ($\approx 2.820$) and the maximum stays at $2.842435$ across three iterations.

For $n=32$ we observe the clearest sustained right-shift at $n=32$: the mean increases from $\approx 2.9161$ (training) to $\approx 2.9180$ (iteration 2), and the upper tail becomes denser while maintaining the maximum $2.939349$.  This regime is particularly sensitive to global rearrangements of the contact graph, and the "centers-only" generator plus LP/SRP push appears to be an effective division of labor.

\remark
The AlphaEvolve report states $\sum r_i\ge 2.635$ for $n=26$ and $\sum r_i\ge 2.937$ for $n=32$ (rounded/thresholded in the paper figure).  Our best values in Figures~\ref{fig:sumr-26} and~\ref{fig:sumr-32} beat these reported thresholds, and the main benefit of \FMBoost{} is that such near-extremal solutions occur with substantially higher frequency in the generated-and-pushed distribution.

\begin{figure}[t]
\centering
\begin{minipage}[t]{0.32\textwidth}
  \centering
  \includegraphics[width=\linewidth]{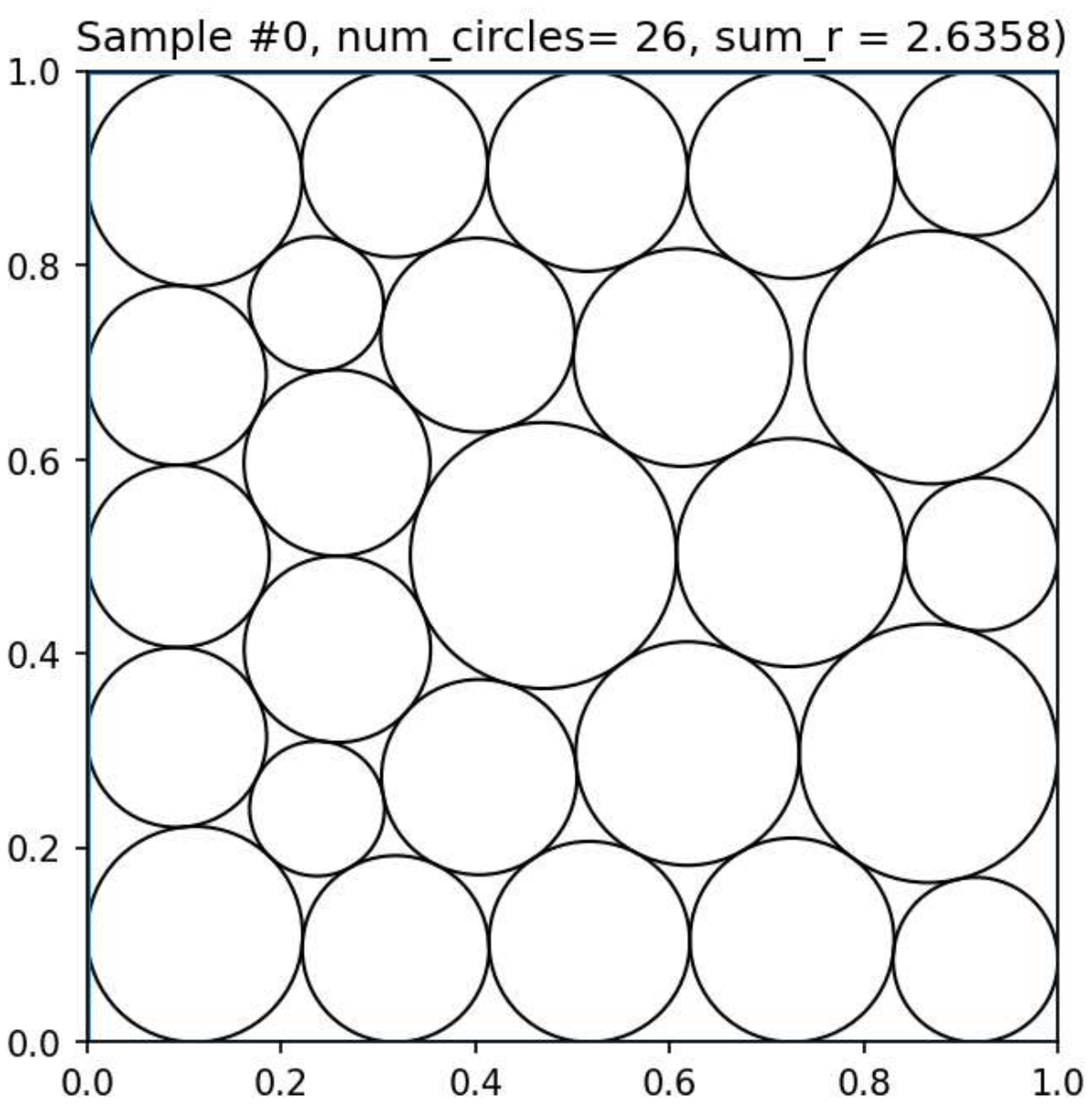}
  \caption*{$n=26$, $\sum r_i=2.6358$ (A)}
\end{minipage}\hfill
\begin{minipage}[t]{0.32\textwidth}
  \centering
  \includegraphics[width=\linewidth]{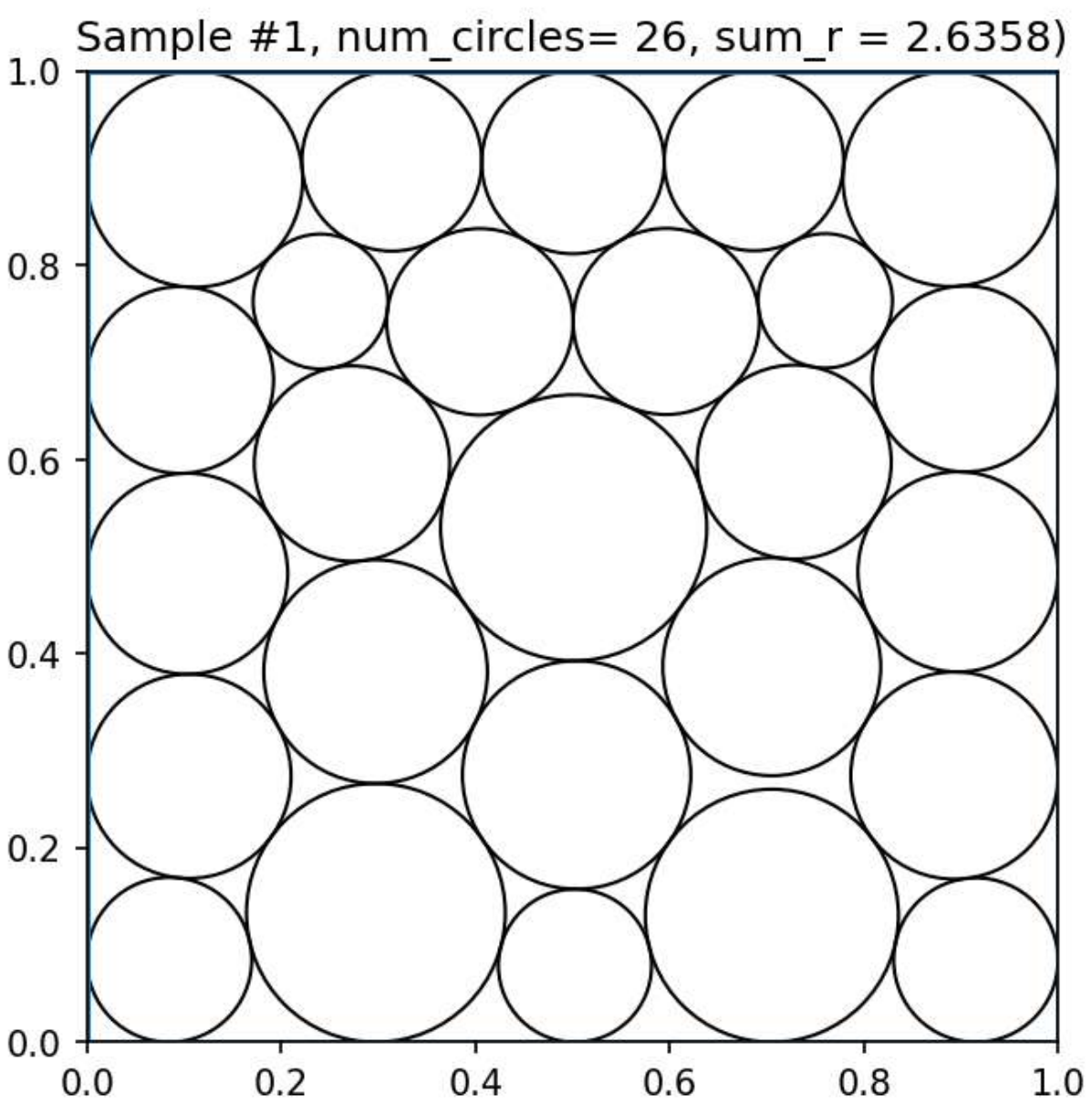}
  \caption*{$n=26$, $\sum r_i=2.6358$ (B)}
\end{minipage}\hfill
\begin{minipage}[t]{0.32\textwidth}
  \centering
  \includegraphics[width=\linewidth]{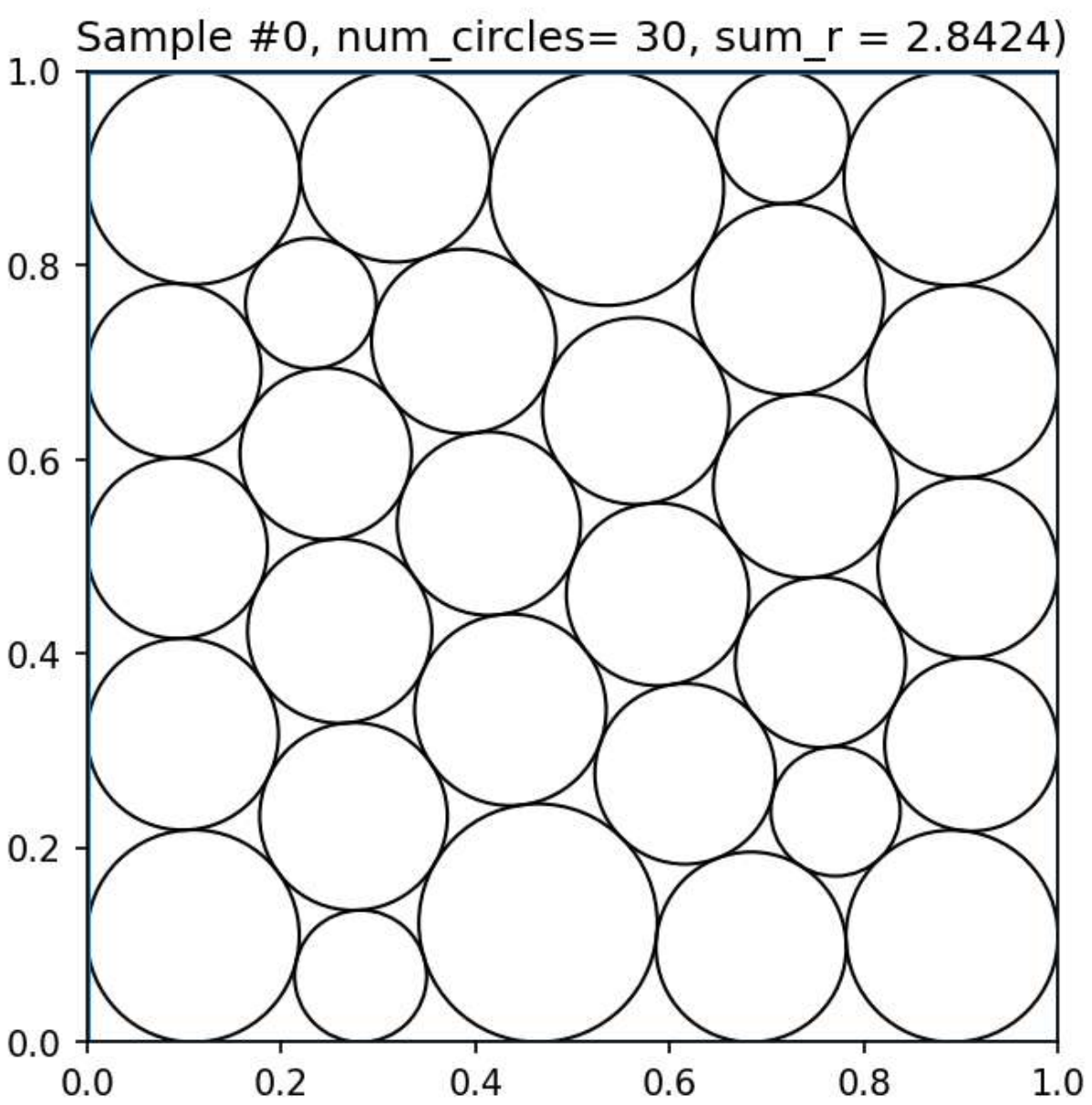}
  \caption*{$n=30$, $\sum r_i=2.8424$ (A)}
\end{minipage}

\vspace{0.8em}

\begin{minipage}[t]{0.32\textwidth}
  \centering
  \includegraphics[width=\linewidth]{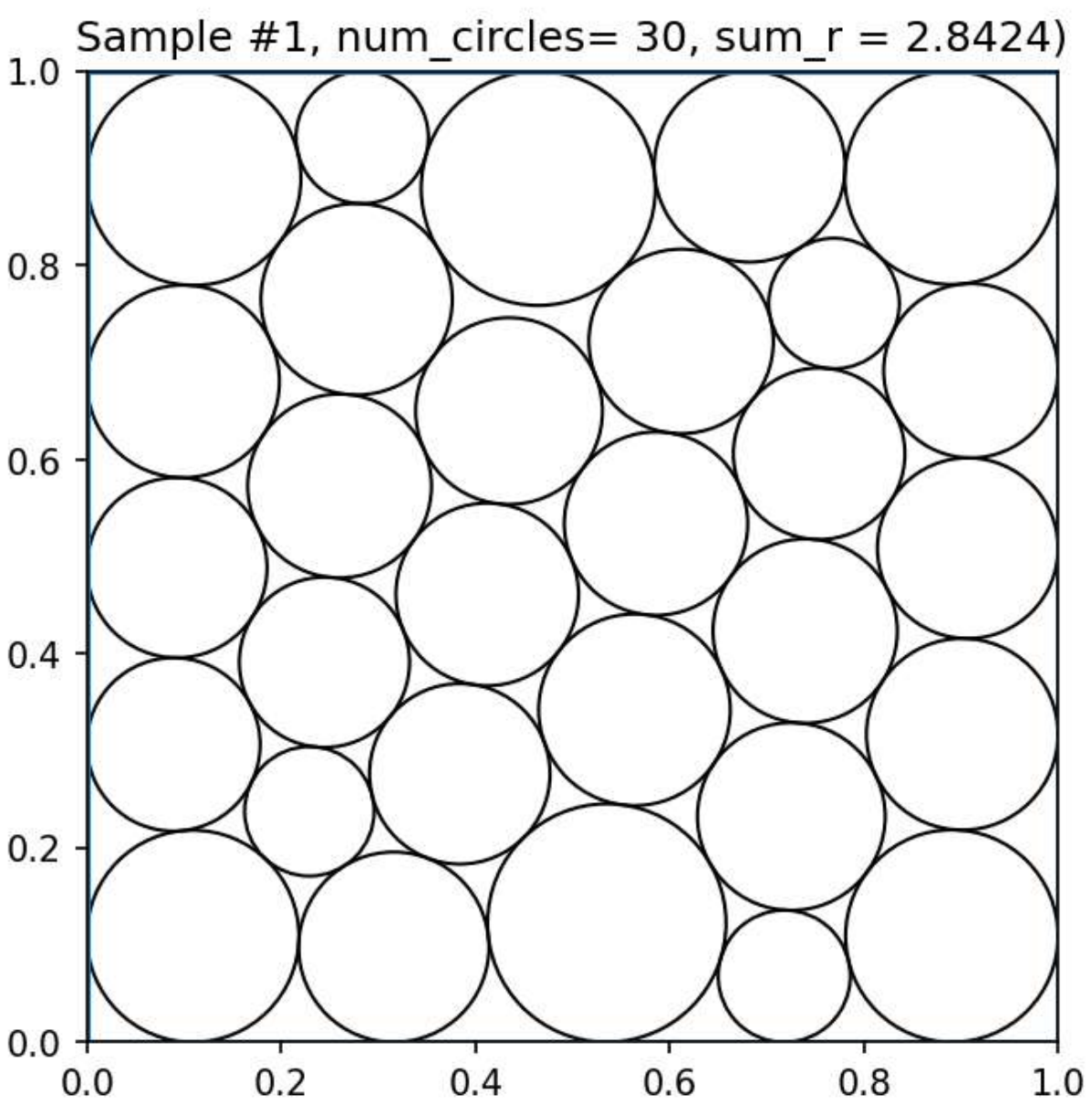}
  \caption*{$n=30$, $\sum r_i=2.8424$ (B)}
\end{minipage}\hfill
\begin{minipage}[t]{0.32\textwidth}
  \centering
  \includegraphics[width=\linewidth]{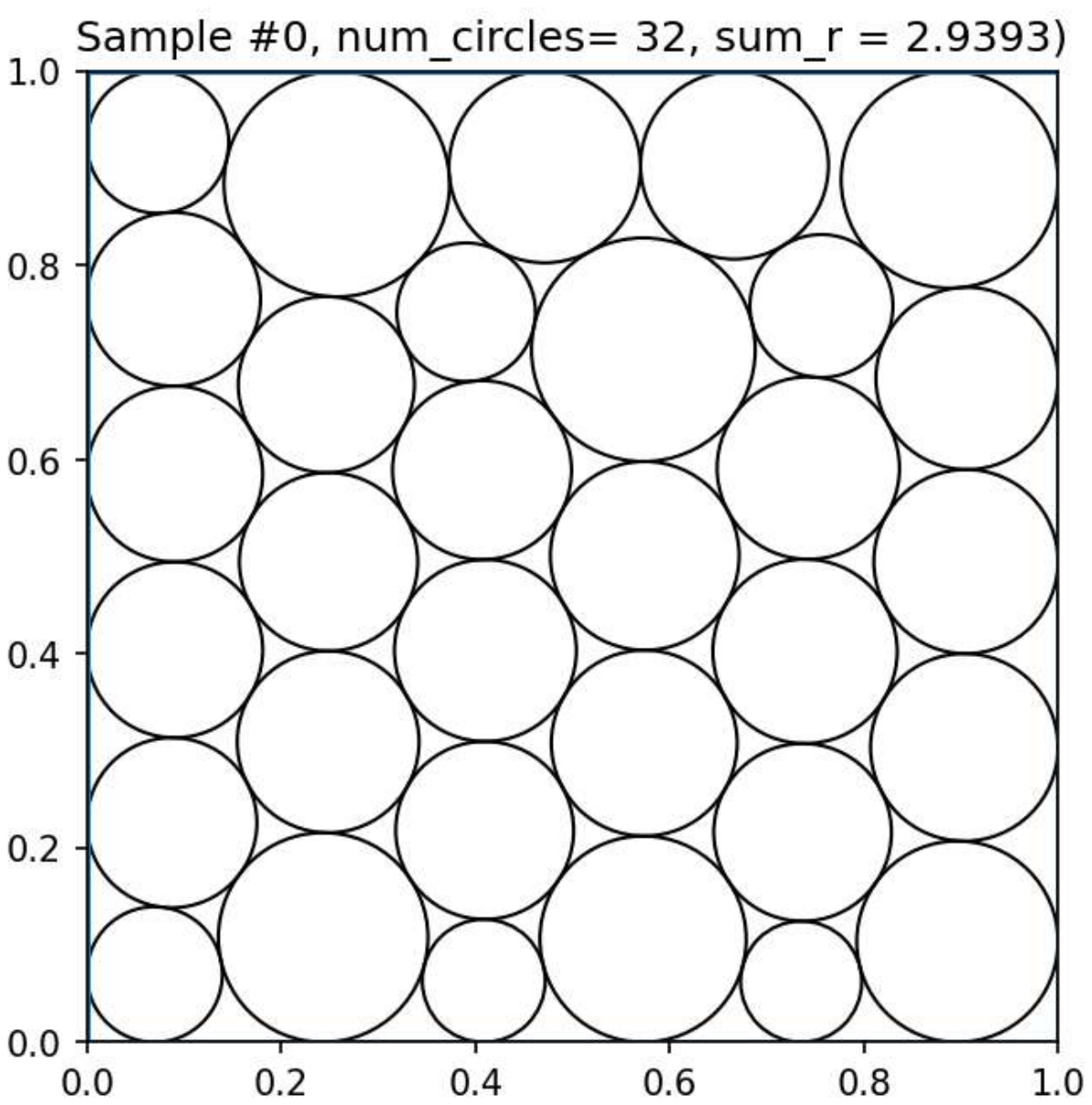}
  \caption*{$n=32$, $\sum r_i=2.9393$ (A)}
\end{minipage}\hfill
\begin{minipage}[t]{0.32\textwidth}
  \centering
  \includegraphics[width=\linewidth]{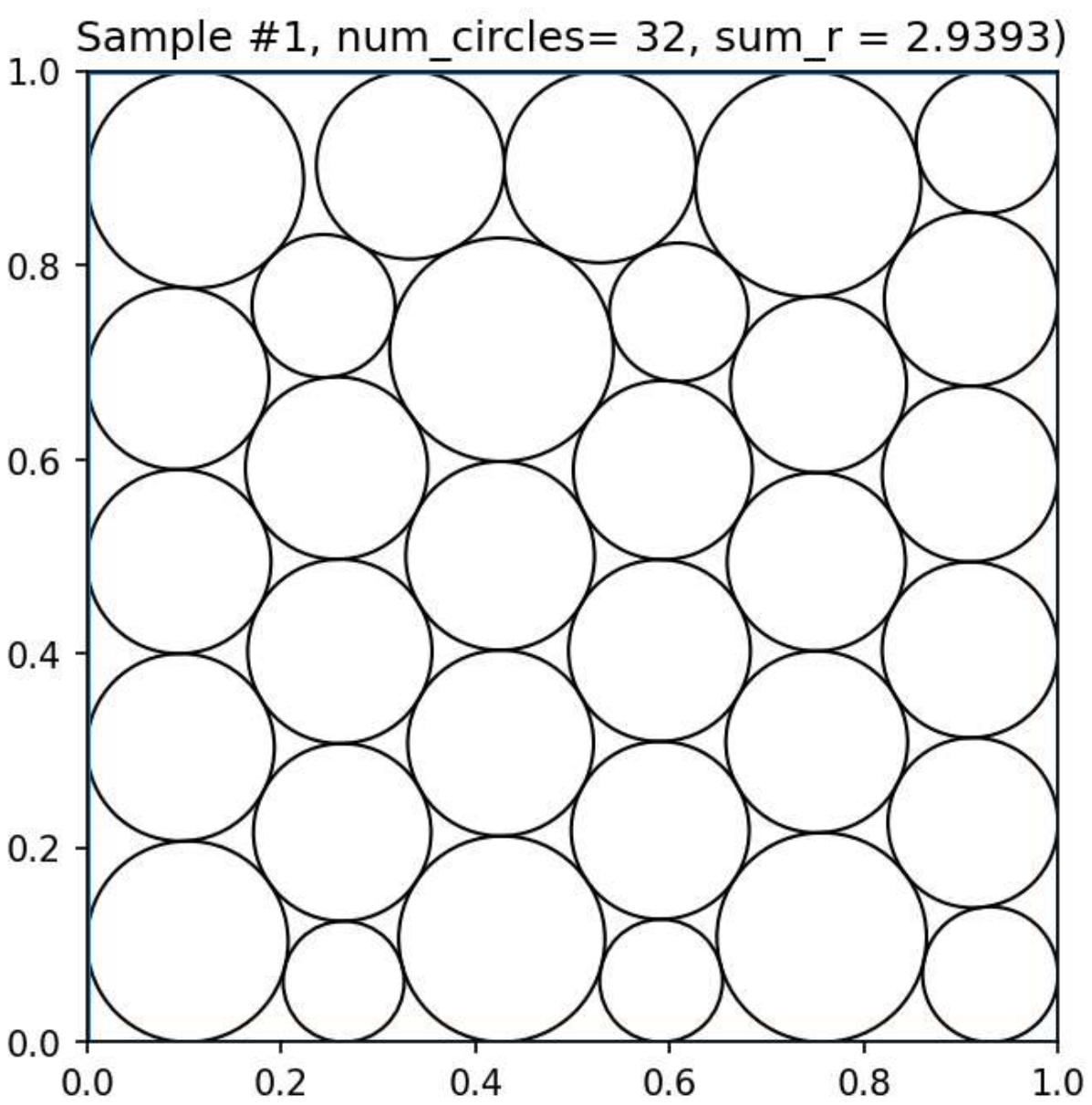}
  \caption*{$n=32$, $\sum r_i=2.9393$ (B)}
\end{minipage}

\caption{\textbf{Best sum-of-radii circle packings found by \FMBoost{} (after final push).}
For each $n\in\{26,30,32\}$ we show two distinct arrangements (A/B) attaining the same best objective value
$\sum_{i=1}^n r_i$. For $n=26$ and $n=32$ these beat the best known result found by AlphaEvolve \cite{novikov2025alphaevolve}.}
\label{fig:sumradii-best-packings}
\end{figure}

\subsection{Star Discrepancy Minimization}
\label{sec:discrepancy} 
The star discrepancy problem asks to place $n$ points in the unit square so that they are as uniformly distributed as possible, i.e. minimizing the worst-case deviation between the fraction of points in any axis-aligned box anchored at the origin and that box's area.

\subsubsection{Local search}

As in the previous problems, SRP serves both as a training-set generator and as the final push (local search) to repair flow samples (Section~2).  
Here the local search is driven by a smooth approximation of the minimax objective in~\eqref{def:star_discrepancy}, obtained by (i) smoothing the box indicators and (ii) replacing the supremum over all anchored boxes by a soft-max over a finite set of candidate boxes.

An anchored box is parametrized by its upper-right corner $q=(a,b)\in[0,1]^2$, i.e.\ $[0,a)\times[0,b)$.  
In the smooth surrogate we do not optimize over a continuum of corners $q$; instead we choose two finite grids
\[
A_x\subset(0,1],\qquad A_y\subset(0,1],
\]
and restrict attention to the discrete family of boxes with corners $(a,b)\in A_x\times A_y$.
We refer to $A_x$ and $A_y$ as the evaluation grids (or anchored grids) for the surrogate: they specify which anchored boxes are probed when approximating the worst-case discrepancy.
In our implementation $A_x,A_y$ are either (a) a fixed uniform grid (for stable comparisons), or (b) a critical grid obtained from the current point coordinates (unique $x$- and $y$-coordinates, with $1$ included), refreshed periodically during SRP.

For an anchored box $(a,b)$ we approximate the indicator $\mathbf{1}_{x_i<a}\mathbf{1}_{y_i<b}$ by a product of sigmoids
\[
\mathbf{1}_{x_i<a}\approx \sigma\!\left(\frac{a-x_i}{\tau}\right),
\qquad
\mathbf{1}_{y_i<b}\approx \sigma\!\left(\frac{b-y_i}{\tau}\right),
\]
with gate temperature $\tau=\texttt{tau\_sigmoid}$.
For $(a,b)\in A_x\times A_y$ we then define the smoothed discrepancy at that box by
\[
\Delta(a,b)
:=
\frac{1}{N}\sum_{i=1}^N 
\sigma\!\left(\frac{a-x_i}{\tau}\right)\sigma\!\left(\frac{b-y_i}{\tau}\right)
-ab,
\qquad
|\Delta|\approx\sqrt{\Delta^2+\varepsilon},
\]
where the $\sqrt{\Delta^2+\varepsilon}$ term (with small $\varepsilon$) provides a differentiable proxy for $|\Delta|$.
The SRP objective then aggregates these values over the grid (via a log-sum-exp soft-max, as described next), so that improving the surrogate corresponds to reducing the worst anchored-box error \emph{among the boxes on the chosen grid}.

We then approximate the supremum in~\eqref{def:star_discrepancy} by a log-sum-exp (softmax) at sharpness $\beta$:
\begin{equation}
\label{eq:star-softmax}
\widetilde D_{\beta,\tau}(P)
:=
\frac{1}{\beta}\log\!\sum_{(a,b)\in A_x\times A_y}\exp\!\Big(\beta\,\sqrt{\Delta(a,b)^2+\varepsilon}\Big).
\end{equation}
SRP anneals $\beta$ from \texttt{beta\_softmax\_start} to \texttt{beta\_softmax\_final}, so early iterations smooth over many boxes, while later iterations concentrate on near-worst boxes.  For efficiency we optionally restrict to the top-$K$ boxes with largest $|\Delta|$ (Top-$K$ mode), and we update the "critical" grid periodically during SRP.

The push stage is SRP on $\widetilde D_{\beta,\tau}$ (with clipping to $[0,1]^2$), followed by L-BFGS-B on the same surrogate.  The exact $D^\ast$ is then recomputed and used for selection and sorting.

\subsubsection{Training and sampling}

The flow model is the same permutation-equivariant set-transformer velocity field as in Section~2, now conditioned by $c(P)=\Bigl(\frac{N}{128},\ D^\ast(P)\Bigr),
$ where $D^\ast(P)$ is the exact discrepancy of the training sample.  In addition to the standard flow-matching regression, we include an auxiliary penalty (again as in Section~2) that discourages the projected endpoint from having larger smooth discrepancy than its target, i.e.\ it penalizes $(\widetilde D_{\beta,\tau}-D^\ast(P))_+$.

Sampling uses a PCFM-style loop: we integrate the learned ODE and interleave short gradient steps that \emph{decrease} the differentiable surrogate~\eqref{eq:star-softmax} (projection) plus a proximal relaxation toward the path blend.  The initial noise distribution $x_1$ is taken to be Latin hypercube sampling (LHS), which provides a stronger geometric baseline than i.i.d.\ uniform points.  As with other tasks, all samples are evaluated only after the SRP final push.

\subsubsection{Results}

\begin{figure}[t]
\centering
\begin{minipage}[t]{0.49\textwidth}
  \centering
  \includegraphics[width=\linewidth]{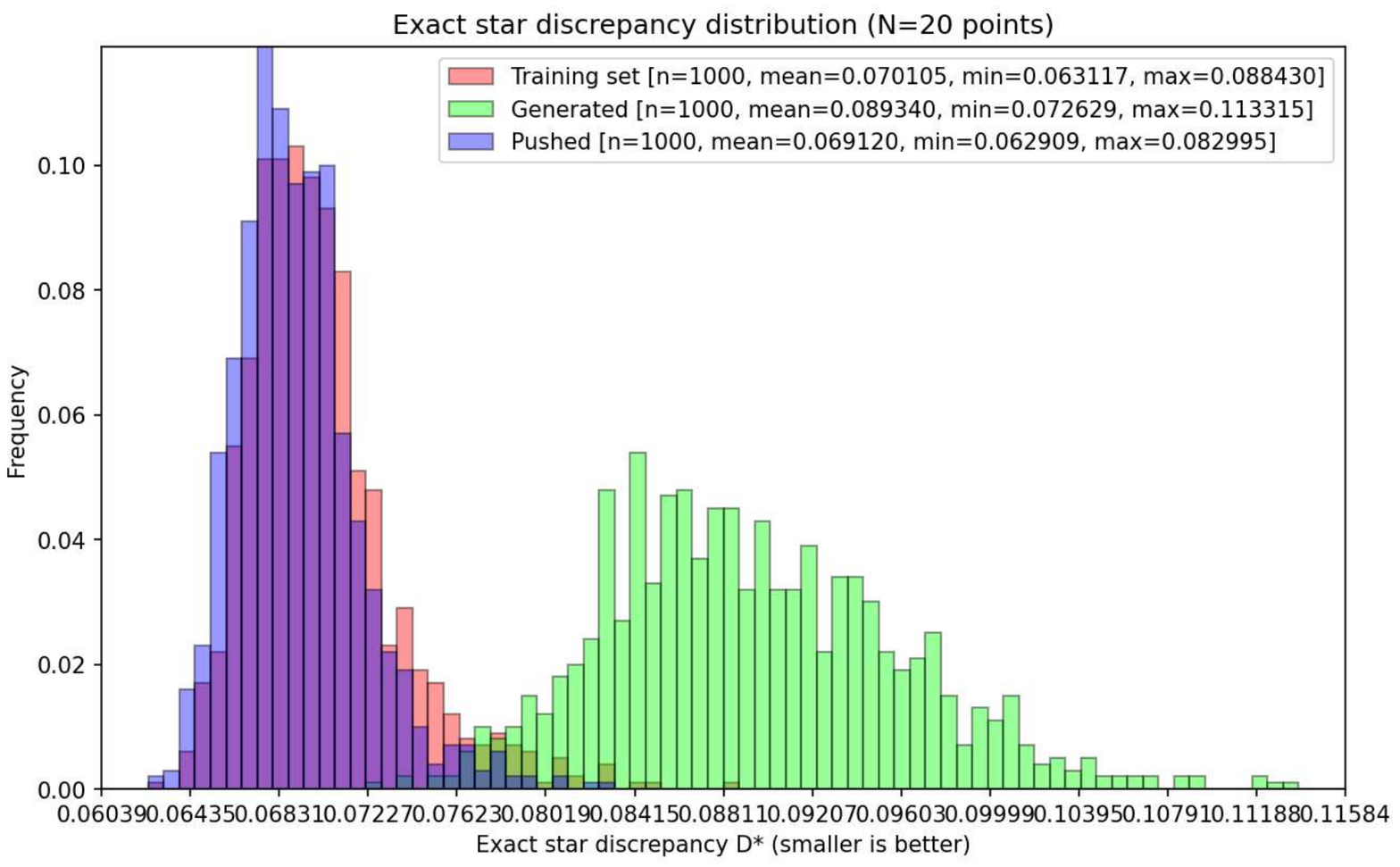}
  \caption*{\textbf{(a)} Training vs.\ generated vs.\ pushed.}
\end{minipage}\hfill
\begin{minipage}[t]{0.49\textwidth}
  \centering
  \includegraphics[width=\linewidth]{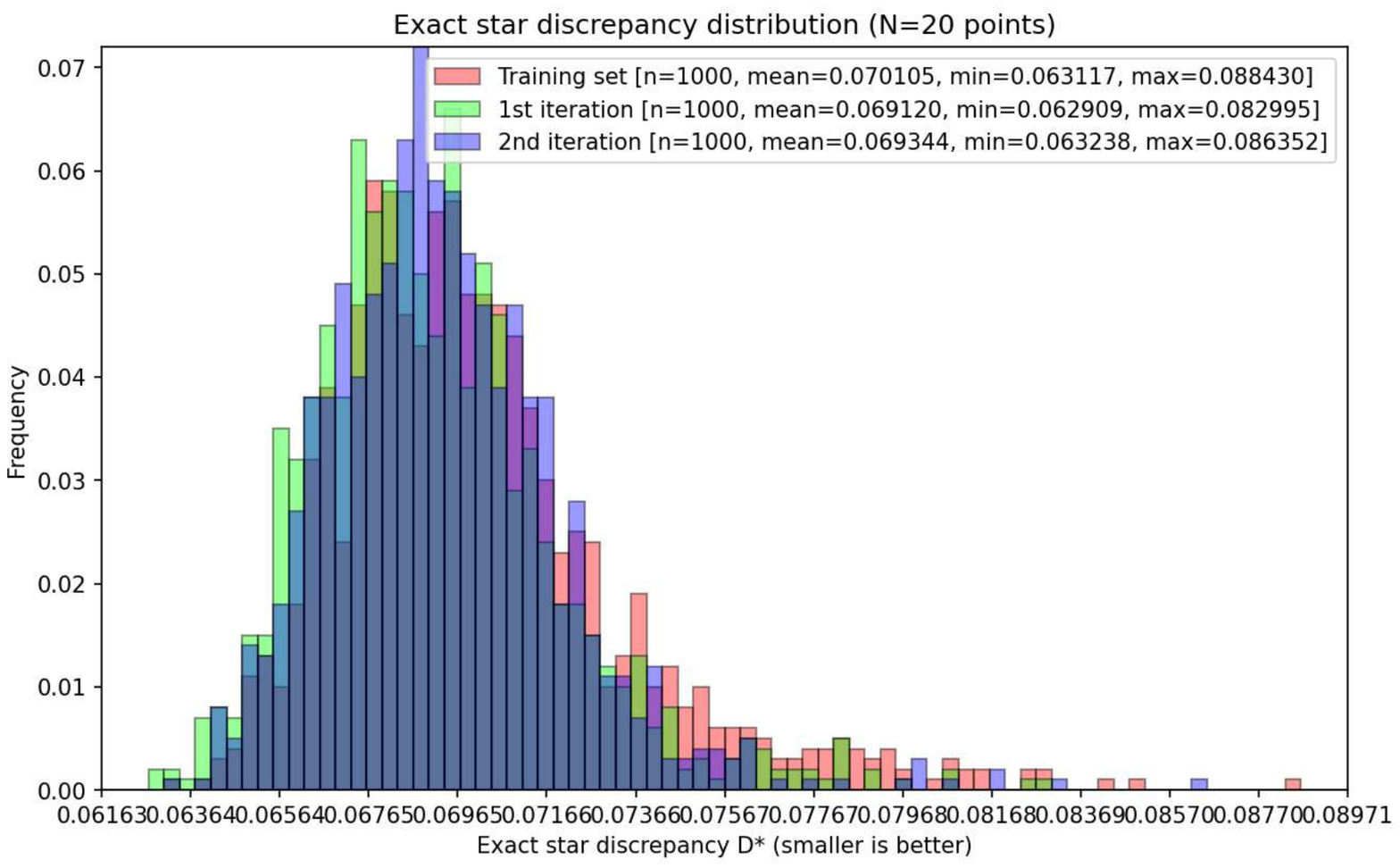}
  \caption*{\textbf{(b)} Training vs.\ pushed iterations 1--2.}
\end{minipage}
\caption{\textbf{Exact star discrepancy, $N=20$.}  Normalized histograms of $D^\ast$ (smaller is better).}
\label{fig:star-n20}
\end{figure}

\begin{figure}[t]
\centering
\begin{minipage}[t]{0.49\textwidth}
  \centering
  \includegraphics[width=\linewidth]{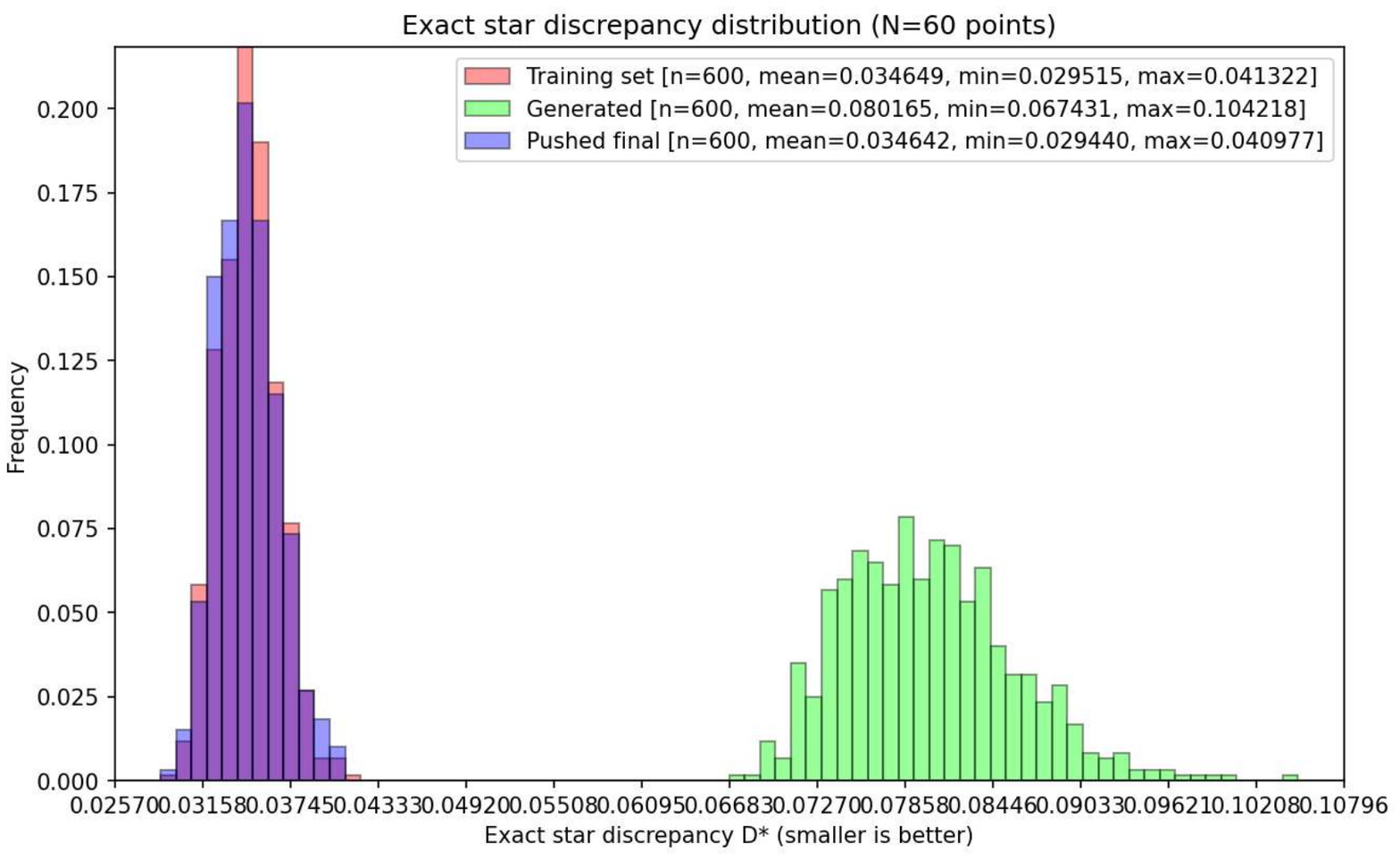}
  \caption*{\textbf{(a)} Training vs.\ generated vs.\ pushed.}
\end{minipage}\hfill
\begin{minipage}[t]{0.49\textwidth}
  \centering
  \includegraphics[width=\linewidth]{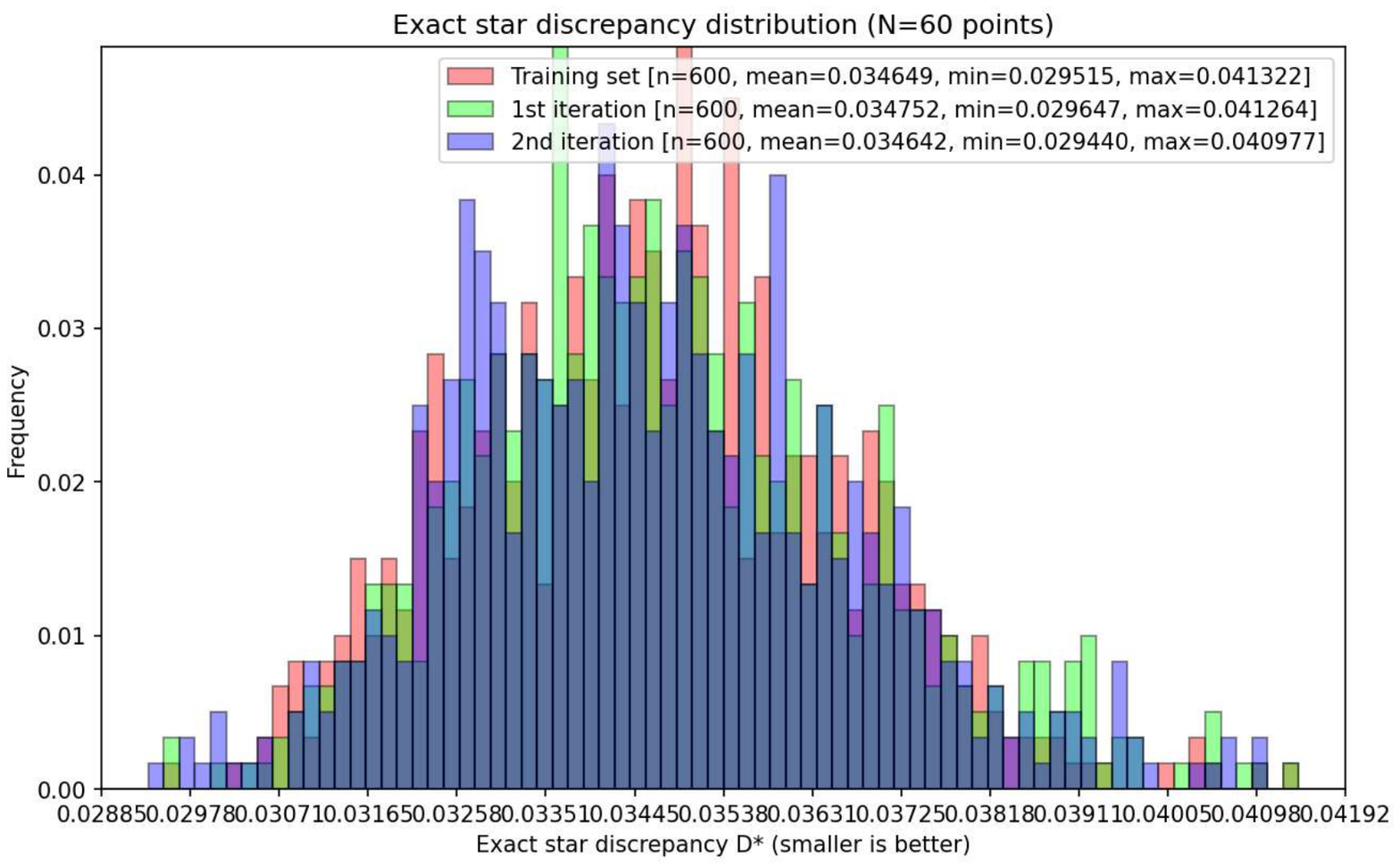}
  \caption*{\textbf{(b)} Training vs.\ pushed iterations 1--2.}
\end{minipage}
\caption{\textbf{Exact star discrepancy, $N=60$.}  The final push is essential: raw generation is far off-scale, but pushing recovers and slightly improves the best tail.}
\label{fig:star-n60}
\end{figure}

Figures~\ref{fig:star-n20} and~\ref{fig:star-n60} summarize experiments for $N=20$ and $N=60$ points (each histogram uses the sample counts stated in the legends).  Two consistent effects stand out like for the other tested problems:

\begin{itemize}[nosep,leftmargin=*]
\item Generation alone is not enough.  The raw flow samples (green) have substantially worse discrepancy (larger $D^\ast$), reflecting that minimax structure is hard to hit without local correction.
\item The final push closes the gap (and can slightly improve records).  After SRP+L-BFGS-B, the pushed distribution (blue) essentially matches the training regime and modestly improves the best observed values.
\end{itemize}

For $N=20$, the training set has mean $0.070105$ with best (minimum) $0.063117$.  Raw generation degrades badly (mean $0.089340$).  After the final push, the mean improves to $0.069120$ and the best value improves slightly to $0.062909$ (Figure~\ref{fig:star-n20}a).  In the iterative view (Figure~\ref{fig:star-n20}b), the first iteration produces the best tail; the second iteration remains better on average than training but does not further improve the minimum.

For $N=60$ the same phenomenon is more dramatic: raw generation produces discrepancies around $0.08$ (mean $0.080165$), while the pushed distribution returns to the training scale $\approx 0.035$ (Figure~\ref{fig:star-n60}a).  Across iterations (Figure~\ref{fig:star-n60}b), improvements are necessarily small because the training distribution is already tight; nevertheless, iteration~2 improves the best observed value from $0.029515$ to $0.029440$ and slightly reduces the worst tail (max drops from $0.041322$ to $0.040977$).

\begin{figure}[t]
\centering
\begin{minipage}[t]{0.49\textwidth}
  \centering
  \includegraphics[width=\linewidth]{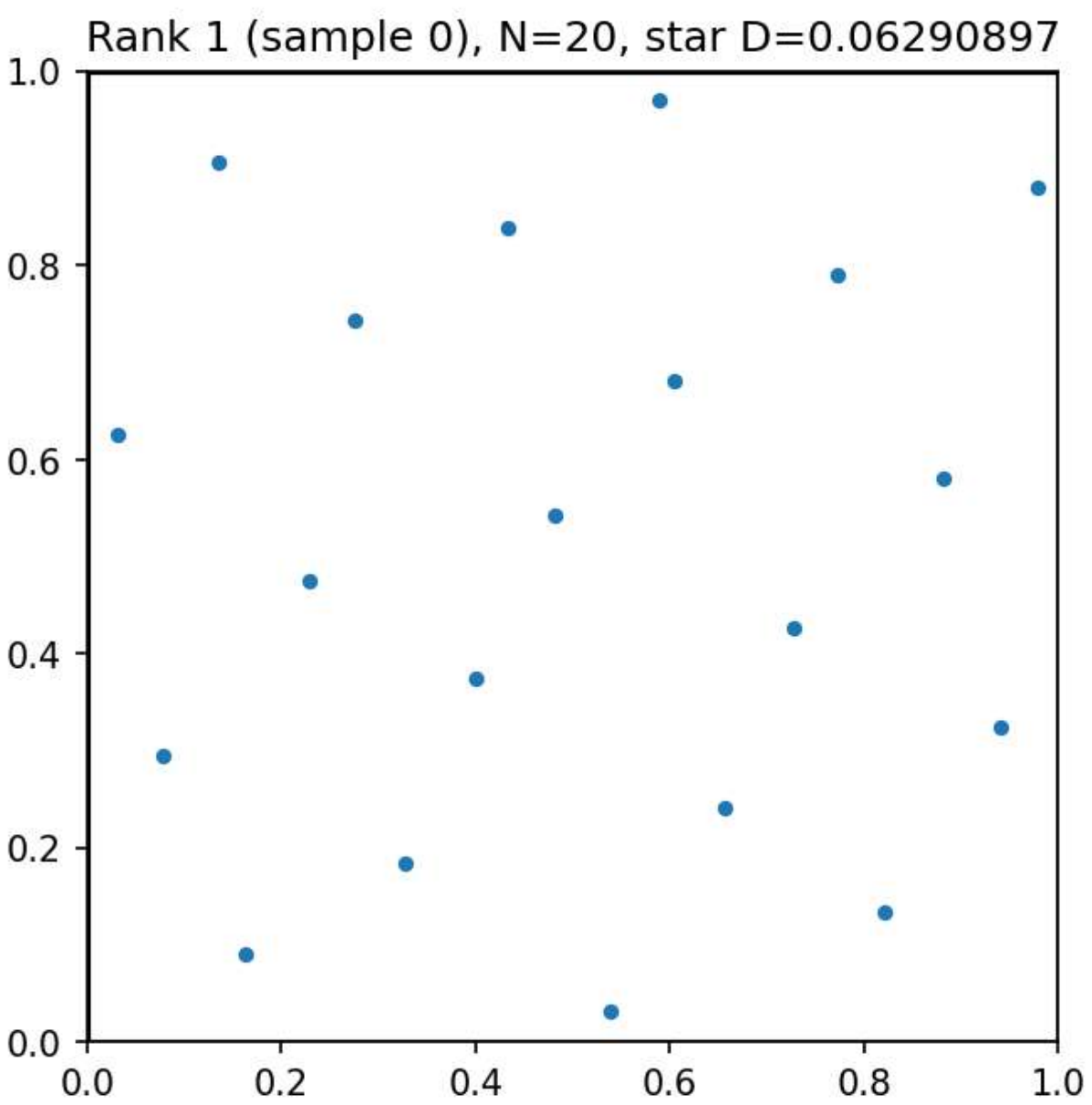}
  \caption*{$N=20$, $D^\ast=0.06290897$}
\end{minipage}\hfill
\begin{minipage}[t]{0.49\textwidth}
  \centering
  \includegraphics[width=\linewidth]{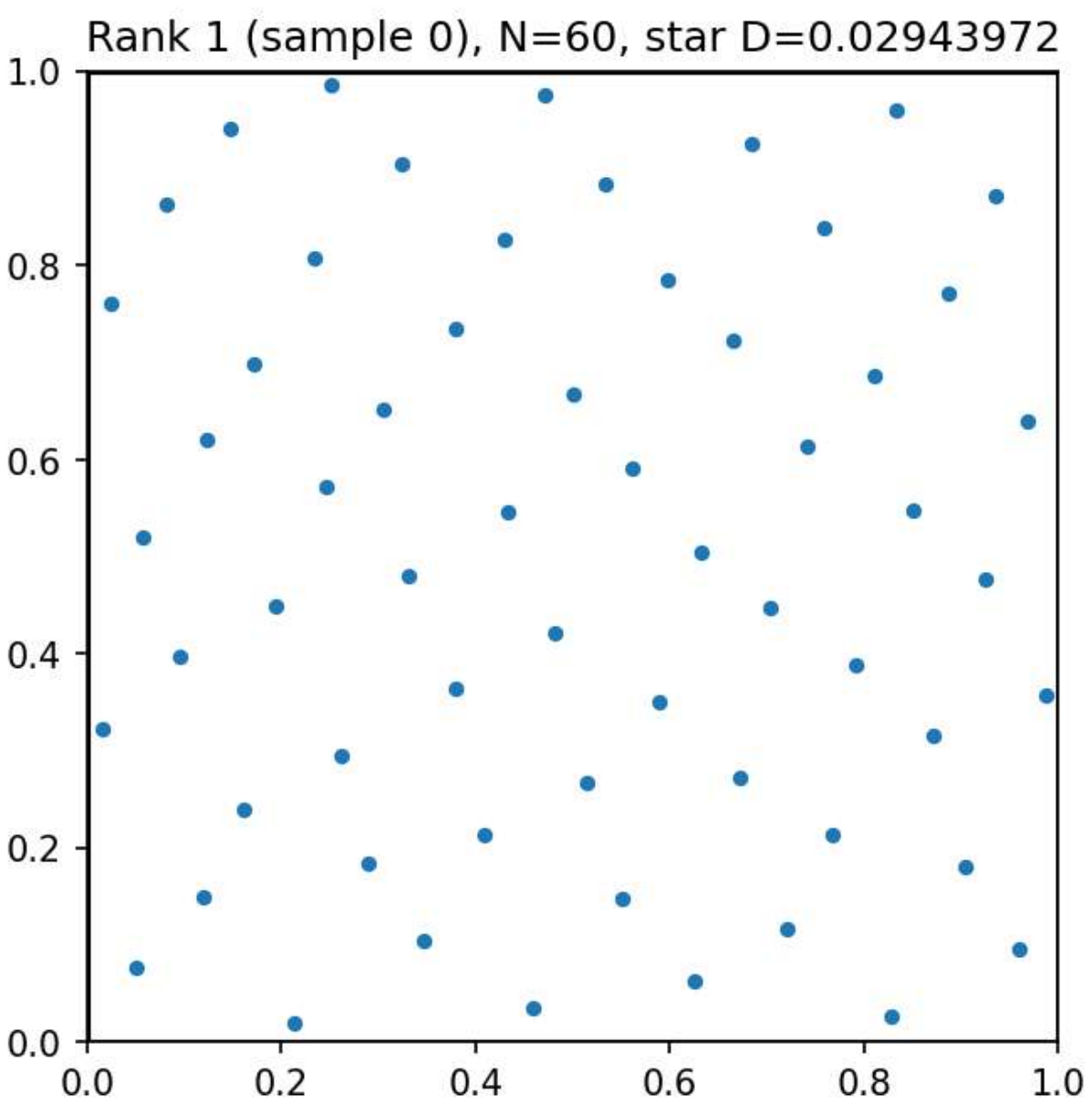}
  \caption*{$N=60$, $D^\ast=0.02943972$}
\end{minipage}

\caption{Best star-discrepancy point sets found by \FMBoost{}: they beat the best known constructions $0.073086$ ($N=20$) and $0.032772$ ($N=60$) in \cite{ALGOLabStarDiscrepancyBenchmark} and for $N=20$ it approaches the optimal $0.0604$ proved in  \cite{ClementDoerrKlamrothPaquete2025OptimalStarDiscrepancy}.} 
\label{fig:star-best-pointsets}
\end{figure}

\section{Conclusion}
\label{sec:discussion} 
We have introduced \FMBoost{}, a reward-guided closed-loop generative framework for continuous simulation-based optimization that synthesizes ideas from flow-based generative modeling, self-supervised representation learning, and reinforcement learning. The core insight is that generative models need not be passive density estimators trained on static datasets; instead, they can be active participants in the optimization process, receiving direct feedback from the objective function and adapting their sampling distribution accordingly.

The framework rests on three components. First, \emph{Conditional Flow Matching} learns a velocity field that transports a source distribution to the empirical distribution of high-quality solutions, providing a smooth, low-dimensional parameterisation of the solution manifold. Second, \emph{Geometry-Aware Sampling} integrates this vector field while enforcing hard constraints through interleaved projection and proximal relaxation, ensuring that generated samples are geometrically valid without sacrificing fidelity to the learned distribution. Third, \emph{Reward-Guided CFM} fine-tunes the generative model online using importance-weighted reward, with a consistency regulariser, inspired by teacher--student architectures in self-supervised learning, that prevents the failure mode of generative collapse.

The transition from open-loop to closed-loop optimization is the key conceptual advance. In previous open-loop methods, the generative model approximates the distribution of current best solutions, and improvement relies on a mere stochastic hope that sampling occasionally escapes to better basins. In our closed-loop method, the model is fine-tuned online with direct gradient signal from the reward function, systematically guiding sampling toward higher-quality regions while the consistency term maintains the diversity necessary for continued exploration. This transforms the boosting loop from an indirect, variance-driven process into a principled policy-optimization algorithm with controllable convergence behaviour.

We have demonstrated the effectiveness of \FMBoost{} on four geometric optimization problems: sphere packing in hypercubes, circle packing maximizing sum of radii, the Heilbronn triangle problem, and star discrepancy minimization. In several cases, \FMBoost{} discovers configurations that match or exceed the best known results. Crucially, the closed-loop variant achieves these results in remarkably few iterations, beating the available data within a single boosting round, whereas open-loop methods require orders of magnitude more samples to reach comparable quality.

\subsection{Discussion and Future Directions}
Several directions for future work emerge from this study. \FMBoost{} is suitable to study a broad family of extremal problems in algebra, geometry, combinatorics and number theory, such as  covering problems, discrepancy minimization in higher dimensions, energy minimization on manifolds, optimal transport, constraint-satisfaction problems with continuous variables (e.g. counter-example mining). Many of these share the structure that makes \FMBoost{} effective, a smooth objective landscape with many local optima, hard constraints that define a complex feasible region, and a lack of exploitable algebraic structure. We anticipate that the same pipeline, with problem-specific adaptations to the projection operator and reward function, will prove effective across this family.

The current framework treats the reward function as a black box, querying it only at generated samples. Incorporating gradient information from differentiable objectives, or learning a surrogate reward model for expensive-to-evaluate functions, could substantially accelerate convergence. Similarly, the teacher-student architecture admits natural extensions, namely an exponential moving average (EMA) teacher, as in DINO \cite{caron2021emerging} and BYOL \cite{grill2020bootstrap}, might provide smoother consistency targets than a frozen checkpoint, and curriculum strategies that progressively tighten the consistency constraint could enable more aggressive exploration in early stages.

From a broader perspective, our experiments suggest that flow-based generative models are not merely tools for image or signal synthesis, but can serve as robust and flexible components in computational and experimental mathematics. The key enabling insight is that the same architectures and training objectives that capture complex data distributions can, with appropriate closed-loop feedback, be repurposed to \emph{improve} upon those distributions toward extremal configurations. We hope that \FMBoost{} contributes to a growing paradigm in which machine learning and experimental mathematics are not separate disciplines but mutually reinforcing tools for discovery.

\section*{Acknowledgment}
The authors would like to thank the Center of Mathematical Sciences and Applications (CMSA) at Harvard University for running the 2024 \emph{Mathematics and Machine Learning} program, the organizers of the \emph{Machine Learning and Mathematics} workshop at the Korea Institute for Advanced Study (KIAS), and the organizers of the 2025 MATRIX--MFO Tandem Workshop \emph{Machine Learning and AI for Mathematics}. Our work benefited greatly from the inspiration and discussions at these events.

G.B.\ gratefully acknowledges help and long discussions with Fran\c{c}ois Charton, Jordan Ellenberg, Geordie Williamson, Adam Zsolt Wagner, Mike Douglas and Amaury Hayat. G.B.\ and J.K.\ were supported by the Independent Research Fund Denmark. 

B.H.\ was supported by the Excellence Cluster ORIGINS, funded by the Deutsche Forschungsgemeinschaft (DFG, German Research Foundation) under Germany’s Excellence Strategy -- EXC-2094-390783311. B.H.\ extends his gratitude to Lukas Heinrich.

\bibliographystyle{unsrt}
\bibliography{references}
\end{document}